\documentclass{amsart}

\usepackage{version}
\includeversion{arxiv}
\excludeversion{duke}

\usepackage{amssymb,amsmath}
\usepackage{epic,eepic}
\usepackage{url}
\usepackage[all]{xypic}
\xyoption{dvips}
\usepackage[neveradjust]{paralist}
\usepackage{slashbox}

\theoremstyle{definition}
\newtheorem{ntn}{Notation}[section]
\newtheorem{dfn}[ntn]{Definition}
\theoremstyle{plain}
\newtheorem{lem}[ntn]{Lemma}
\newtheorem{prp}[ntn]{Proposition}
\newtheorem{thm}[ntn]{Theorem}
\newtheorem{cor}[ntn]{Corollary}
\newtheorem{cnj}[ntn]{Conjecture}
\theoremstyle{remark}

\newtheorem{rmk}[ntn]{Remark}
\newtheorem*{rmk*}{Remark}
\newtheorem{exa}[ntn]{Example}

\newcommand{\ol}[1]{{\overline{#1}}}
\newcommand{\del}{\partial}
\newcommand{\eps}{\varepsilon}
\newcommand{\ideal}[1]{{\langle#1\rangle}}
\newcommand{\xymat}{\SelectTips{cm}{}\xymatrix}

\newcommand{\boldzero}{{\mathbf{0}}}
\newcommand{\boldone}{{\mathbf{1}}}
\newcommand{\bolda}{{\mathbf{a}}}
\newcommand{\boldb}{{\mathbf{b}}}
\newcommand{\boldc}{{\mathbf{c}}}
\newcommand{\boldu}{{\mathbf{u}}}
\newcommand{\boldv}{{\mathbf{v}}}

\newcommand{\calD}{\mathcal{D}}
\newcommand{\calE}{\mathcal{E}}
\newcommand{\calF}{\mathcal{F}}
\newcommand{\calH}{\mathcal{H}}
\newcommand{\calI}{\mathcal{I}}
\newcommand{\calK}{\mathcal{K}}

\newcommand{\calM}{\mathcal{M}}
\newcommand{\calO}{\mathcal{O}}
\newcommand{\calP}{\mathcal{P}}
\newcommand{\calR}{\mathcal{R}}

\newcommand{\calW}{\mathcal{W}}

\newcommand{\frakm}{\mathfrak{m}}
\newcommand{\frakP}{\mathfrak{P}}
\newcommand{\frakT}{\mathfrak{T}}
\newcommand{\frakV}{\mathfrak{V}}

\newcommand{\CC}{\mathbb{C}}

\newcommand{\NN}{\mathbb{N}}

\newcommand{\PP}{\mathbb{P}}
\newcommand{\QQ}{\mathbb{Q}}

\newcommand{\TT}{\mathbb{T}}

\newcommand{\ZZ}{\mathbb{Z}}

\DeclareMathOperator{\ch}{Ch}
\DeclareMathOperator{\cc}{CC}
\DeclareMathOperator{\codim}{codim}

\DeclareMathOperator{\conv}{conv}

\DeclareMathOperator{\GL}{Gl}

\DeclareMathOperator{\gr}{gr}

\DeclareMathOperator{\Hom}{Hom}
\DeclareMathOperator{\image}{im}

\DeclareMathOperator{\Irr}{{\mathcal I}\it{rr}}

\DeclareMathOperator{\N}{N}
\DeclareMathOperator{\modulo}{\,\,mod\,\,}

\DeclareMathOperator{\Orb}{Orb}

\DeclareMathOperator{\Proj}{Proj}
\DeclareMathOperator{\Rees}{Rees}
\DeclareMathOperator{\rk}{rk}
\DeclareMathOperator{\SHom}{{\mathcal Hom}}
\DeclareMathOperator{\sign}{sign}
\DeclareMathOperator{\Spec}{Spec}
\DeclareMathOperator{\supp}{supp}
\DeclareMathOperator{\Tor}{Tor}
\DeclareMathOperator{\vol}{vol}
\DeclareMathOperator{\Var}{Var}

\makeindex

\begin{document}

\title[Slopes of hypergeometric systems]{Irregularity of hypergeometric
systems\\ via slopes along coordinate subspaces}

\author{Mathias Schulze}
\address{
M. Schulze\\
Oklahoma State University\\
Dept. of Mathematics\\
401 MSCS\\
Stillwater, OK 74078\\
USA}
\email{mschulze@math.okstate.edu}
\thanks{MS was supported by the Humboldt foundation.}

\author{Uli Walther}
\address{
U. Walther\\
Purdue University\\
Dept. of Mathematics\\
150 N. University St.\\
West Lafayette, IN 47907\\
USA}
\email{walther@math.purdue.edu}
\thanks{UW was supported by the NSF under grant DMS~0555319, and by the NSA under grant H98230-06-1-0012.}

\begin{abstract}
We study the irregularity sheaves attached to the $A$-hypergeometric $D$-module $M_A(\beta)$ introduced by I.M.~Gel'fand et al.~\cite{GGZ87,GKZ89}, where $A\in\mathbb{Z}^{d\times n}$ is pointed of full rank and $\beta\in\mathbb{C}^d$.
More precisely, we investigate the slopes of this module along coordinate subspaces.
 
In the process we describe the associated graded ring to a positive semigroup ring for a filtration defined by an arbitrary weight vector $L $ on torus-equivariant generators. 
To this end we introduce the $(A,L)$-umbrella, a cell complex determined by $A$ and $L$, and identify its facets with the components of the associated graded ring.
 
We then establish a correspondence between the full $(A,L)$-umbrella and the components of the $L$-characteristic variety of $M_A(\beta)$. 
We compute in combinatorial terms the multiplicities of these components in the $L$-characteristic cycle of the associated Euler--Koszul complex, identifying them with certain intersection multiplicities.
 
We deduce from this that slopes of $M_A(\beta)$ are combinatorial, independent of $\beta$, and in one-to-one correspondence with jumps of the $(A,L)$-umbrella. 
This confirms a conjecture of Sturmfels and gives a converse of a theorem of R.~Hotta~\cite[Ch.~II, \S6.2, Thm.]{Hot98}: 
$M_A(\beta)$ is regular if and only if $A$ defines a projective variety.
\end{abstract}



\maketitle
\newpage
\tableofcontents
\numberwithin{equation}{subsection}

\section{Introduction and overview}\label{intro}

The solutions $f$ of an ordinary linear differential equation $P\bullet f(x)=0$ where
\[
P=p_m(x)\del_x^m+p_{m-1}(x)\del_x^{m-1}+\cdots+p_1(x)\del_x+p_0(x)
\]
form a $\CC$-vector bundle of dimension $m$ away from the zero locus of $p_m$. 
At $p_m(x_0)=0$, two types of singularities may occur: at a regular singular point $x_0$, the (multivalued) solutions have polynomial growth for $x\rightarrow x_0$ while in all other cases
$x_0$ is called irregular. 
By Fuchs' Theorem, $P$ is regular at the origin if and only if the Newton polygon $\N(P)$ of $P$ is a quadrant (see \cite{Inc44}). 
The slopes (or critical indices) of $P$ are just the slopes of $\N(P)$; they represent a refined notion of irregularity (see \cite{Lau85}) and indicate the growth of solutions near the critical point. 
Regularity at $x_0$ is equivalent to equality of formal and convergent solutions at $x_0$.

The concept of regularity in higher dimension is considerably more involved. 
Denote by $\calO_X$ the structure sheaf of the complex manifold $X$ and by $\calO_{\widehat{X|Y}}$ the completion of $\calO_X$ along the submanifold $Y$. 
For any coherent $\calD_X$-module $\calM$, Z.~Mebkhout introduced the irregularity complex $\Irr_Y(M)={\bf R}\SHom_{\calD_X}(\calM,\calO_{\widehat{X|Y}}/\calO_{X|Y})$, $\calO_{X|Y}$ being the restriction of $\calO_X$ to $Y$ \cite{Meb89}. 
This intrinsically analytic notion is inspired by earlier work of B.~Malgrange on regularity in the univariate case \cite{Mal74}, but quite difficult to use.

On the other hand, Y.~Laurent \cite{Lau87} gave a generalization for
the concept of a slope to the multivariate case based on more algebraic
methods. With $X,Y$ as above, let $V$ be the Kashiwara--Malgrange
filtration along $Y$, let $F$ be the order filtration, and put
$L=pF+qV$.  Then the rational number $p/q>0$ is a slope of $\calM$
along $Y$ if the $L$-characteristic variety
$\ch^L(\calM)=\supp(\gr^L(\calM))\subseteq T^*X$ jumps (is not locally
constant) at $p/q$.

The theorem of Fuchs generalizes to the multivariate case: by the
analytic-algebraic comparison theorem for slopes
\cite[Thm.~2.4.2]{LM99}, $\calM$ has no slopes along $Y$ precisely if
$\Irr_Y(\calM)$ is exact. In fact, the slopes agree with the
jumps of the Gevrey filtration on the irregularity complex which
provide a measure of growth for the solutions of $\calM$ near $Y$.

All slopes of $\calM$ along $Y$ are rational and there are only a
finite number of them, (see \cite{Lau87}). 
A.~Assi, F.~Castro, and M.~Granger \cite{ACG96} developed a Gr\"obner basis algorithm to compute slopes of algebraic $D$-modules using the algebraic counterparts of $F$, $V$, and $L$ on the Weyl algebra $D$. 
In the process they proved a comparison theorem: slopes of modules over the Weyl algebra can be computed without leaving the algebraic category, where rationality and finiteness follow from the existence of the Gr\"obner fan (see \cite{ACG00}). 
Explicit formul\ae\ for slopes of $D$-modules are very rare. 
The purpose of this article is to describe the slopes of $A$-hypergeometric $D$-modules.

I.M.~Gel'fand, M.I.~Graev, M.M.~Kapranov and A.V.~Zelevinski{\u\i} \cite{GGZ87,GKZ89} defined a class of $D$-modules that includes as particular cases the differential systems satisfied by the classical hypergeometric functions of Gau\ss, Appell, and others.
These $A$-hypergeometric, or GKZ (after Gel'fand, Kapranov, and Zelevinski\u\i), systems are special cases of the equivariant $D$-modules of R.~Hotta and M.~Kashiwara \cite{HK84,Hot98}.
A $d\times n$ integer matrix $A$ defines an action of a $d$-torus $\TT:=(\CC^*)^d$ on $T^*_0\CC^n=\CC^n$.
Our general hypothesis is that $\NN A$ is a positive semigroup with $\ZZ A=\ZZ^d$.
The closure of the orbit through $(1,\dots,1)$ is defined by the toric ideal $I_A\subseteq\CC[\del]=:R$ where $\del:=\del_1,\dots,\del_n$ and $S_A:=R/I_A=\CC[\NN A]$ is the associated semigroup ring.
The Euler vector fields $E=(E_1,\dots,E_d)$ are the pushforwards to this orbit of the Lie algebra generators $t_1\del_1,\dots,t_d\del_d$ of $\TT$.
The $A$-hypergeometric system $M_A(\beta)$, depending on the Lie algebra character $\beta\in\CC^d$, is the $D$-module defined by $I_A$ and the Euler operators $E-\beta$.
It arises in various situations in algebraic geometry such as in the theory of toric residues (see \cite{CDS01}), the study of hyperplane arrangements (see \cite{OT01}), and in the Picard--Fuchs equations governing the variation of Hodge structures for Calabi--Yau toric hypersurfaces (see \cite{CK99}).

By a theorem of R.~Hotta \cite[Ch.~II, \S6.2, Thm.]{Hot98}, homogeneous $A$-hypergeometric systems are regular and hence have no slopes.
In dimension one and in codimension one, slopes of $M_A(\beta)$ were studied by F.~Castro and N.~Takayama \cite{CT03} and M.~Hartillo \cite{Har03,Har05}.
Cohen--Macaulayness of the toric rings in question makes these cases comparatively tractable.
In our general situation, a key tool is a (generalization of a) computation by A.~Adolphson \cite{Ado94} identifying candidate components of the $F$-characteristic variety with the set of faces not containing $0$ of the convex hull $\Delta^F_A$ of $0$ and the columns $\bolda_1,\dots,\bolda_n$ of $A$.
In Section~\ref{spectoric}, $L$ is the filtration on $R$ defined by an
arbitrary  weight vector on $\del_1,\dots,\del_n$.
We introduce  in Definition~\ref{35} the $(A,L)$-polyhedron $\Delta_A^L$ as the convex hull in projective space of $0$ and all $\bolda_i^L:=\bolda_i/\deg^L(\del_i)$ for $i=1,\dots,n$. 
Its faces not containing zero form the $(A,L)$-umbrella $\Phi^L_A$, a
combinatorial object independent of $\beta$ that encodes all
information regarding the $L$-characteristic variety of $M_A(\beta)$.

\begin{figure}[ht]
\caption{Some $(A,L)$-umbrellas. (Shaded $\Delta^L_A$ with fat boundary
$\Phi^L_A$.)}\label{70}
\begin{center}
\setlength{\unitlength}{0.8mm}
\begin{picture}(150,40)(-5,-5)
\put(0,0){\shade\path(0,0)(0,30)(40,10)(10,0)
\put(0,0){\thicklines\path(0,30)(40,10)(10,0)}
\put(10,20){\circle*{1}\makebox(0,-1)[t]{$\bolda_3$}}
\put(5,35){\makebox(0,0)[tl]{$L=(1,1,1,1)$}}}
\put(50,0){\shade\path(0,0)(0,30)(10,20)(20,5)(10,0)
\put(0,0){\thicklines\path(0,30)(10,20)(20,5)(10,0)}
\put(10,20){\circle*{1}\makebox(0,0)[bl]{$\bolda_3$}}
\put(20,5){\circle*{1}\put(-1,1){\makebox(0,0)[bl]{$\bolda_4^L$}}}
\put(5,35){\makebox(0,0)[tl]{$L=(1,1,1,2)$}}}
\put(100,0){\shade\path(0,0)(0,30)(10,20)(10,0)(0,0)
\put(0,0){\thicklines\path(0,30)(10,20)(10,0)}
\put(10,20){\circle*{1}\makebox(0,0)[bl]{$\bolda_3$}}
\put(8,2){\circle*{1}\put(0,0){\makebox(0,0)[br]{$\bolda_4^L$}}}
\put(5,35){\makebox(0,0)[tl]{$L=(1,1,1,5)$}}}
\multiput(0,0)(50,0){3}{
\put(0,0){\vector(1,0){45}}
\put(0,0){\vector(0,1){35}}
\multiput(0,0)(10,0){5}{\path(0,-1)(0,1)}
\multiput(0,0)(0,10){4}{\path(-1,0)(1,0)}
\put(0,30){\circle*{1}\put(-1,0){\makebox(0,0)[r]{$\bolda_1$}}}
\put(10,0){\circle*{1}\put(0,-1){\makebox(0,0)[t]{$\bolda_2$}}}
\put(40,10){\circle*{1}\put(0,-1){\makebox(0,0)[lt]{$\bolda_4$}}}
\dottedline{1}(0,0)(45,11.25)}
\end{picture}
\end{center}
\end{figure}

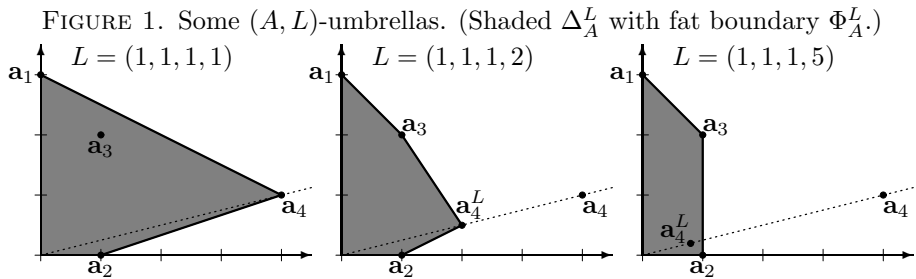
Figure~\ref{70} shows the $(A,L)$-umbrella in the example
\[
A=\begin{pmatrix}0&1&1&4\\3&0&2&1\end{pmatrix}
\]
for the family of filtrations $L=(1,1,1,t)$ for three parameters resulting in
combinatorially different umbrellas.
In Theorem~\ref{13}, we identify $\Phi^L_A$ with both the
$\ZZ^d$-graded spectrum of the $L$-graded toric ring
$S_A^L:=\gr^L(S_A)$ and the collection of torus orbits in $\Spec(S^L_A)$.
In Theorem~\ref{36}, we characterize Cohen--Macaulayness of $S_A^L$ by the corresponding property of its subrings generated by facet cones.

In Section~\ref{charvar}, we extend $L$ to the Weyl algebra
$D=\CC[x]\ideal{\del}$ in the variables $x:=x_1,\dots,x_n$ in such a
way that $E$ is $L$-homogeneous and $W:=\gr^L(D)$ is a polynomial ring.
The family of filtrations in Figure~\ref{70} may be viewed as the
restriction of $L=F+qV$ to $R$ for the $V$-filtration along $x_4$ for
$q=0,1,4$.  In Proposition~\ref{55} we show that all components of the
$L$-characteristic variety of $M_A(\beta)$ correspond to faces
$\tau\in\Phi_A^L$: each is the closure $\bar C_A^\tau$ of the conormal
space $C_A^\tau$ to a torus orbit in $\Spec(S_A^L)$.  The facet components
correspond to orbits in the smooth, and hence Cohen--Macaulay, locus
of $S_A^L$ on which the $L$-symbols of the Euler vector fields form a
regular sequence.  From this we conclude in Proposition~\ref{10} that
the facet components actually occur and their multiplicity is given by
an index formula independent of the parameter $\beta$.  This shows in
particular that $p/q$ is a slope of $M_A(\beta)$ at the origin
whenever $\Phi_A^{pF+qV}$ jumps at $p/q$.  We thus obtain a converse
to Hotta's theorem in Corollary~\ref{38}: regular $A$-hypergeometric
systems are homogeneous.

Since orbits to nonfacets may be outside the Cohen--Macaulay locus of
$S_A^L$, we consider in Section~\ref{cycle} the full Euler--Koszul
complex $K_{A,\bullet}(S_A;\beta)$ from \cite{MMW05}.  In order to
apply methods of homological algebra, we discuss the basics of good
$L$-filtrations on $R$- and $D$-modules.  Using the spectral sequence
for the $L$-filtration on $K_{A,\bullet}(S_A;\beta)$, we identify its
$L$-characteristic cycle with the intersection cycle between $S_A^L$
and the $L$-graded Euler ideal $\gr^L(\ideal{E})$.  In
Theorem~\ref{50}, we use Serre's Intersection Theorem~\cite[Ch.~V, \S C.1, Thm.~1]{Ser65} to
show that the $L$-characteristic variety of $K_{A,\bullet}(S_A;\beta)$
contains all candidate components.  To show that this also holds for
$M_A(\beta)=H_0(K_{A,\bullet}(S_A;\beta))$ we use in Theorem~\ref{60}
results from \cite{MMW05} and an induction argument on toric modules.
In particular (see Corollary~\ref{23}), the components of the $L$-characteristic variety of $M_A(\beta)$ are in one-to-one correspondence to the faces in the $(A,L)$-umbrella $\Phi_A^L$. 
It follows (see Corollary~\ref{39}) that the slopes of $M_A(\beta)$ along coordinate subspaces at the origin correspond exactly to jump parameters $p/q$ of $\Phi_A^{pF+qV}$, confirming and extending a conjecture by B.~Sturmfels.

The holonomic rank $\rk(M):=\dim_{\CC(x)}(M\otimes_{\CC[x]}\CC(x))$ of
a $D$-module $M$ is the dimension of the $\CC$-vector space of its
solution space near a regular point. 
It is a classical result (see \cite{GKZ89,Ado94}) that for generic $\beta$ the rank of $M_A(\beta)$
equals the volume of the convex hull $\Delta^F_A$ of $0$ and
$\bolda_1,\dots,\bolda_n$ where the volume of the unit simplex is
normalized to $1$.  The exceptional set
$\calE(A)$\index{EA@$\calE(A)$} of $A$ is the set of rank-jumping
parameters: $\beta\in\calE(A)$ precisely if the rank of $M_A(\beta)$
exceeds the volume.  By \cite{MMW05}, $\calE(A)$ is a finite subspace
arrangement.  In Theorem~\ref{26} we give a general index/volume
formula for the multiplicity $\mu^{L,\tau}_A$ of $\bar C_A^\tau$,
$\tau\in\Phi^L_A$, in the $L$-characteristic cycle of
$K_{A,\bullet}(S_A;\beta)$. We show then that for non-rank-jumping
parameter $\beta$, or if $\tau$ is a facet, the number
$\mu^{L,\tau}_A$ equals the multiplicity $\mu^{L,\tau}_{A,0}(\beta)$
of $\bar C_A^\tau$ in the $L$-characteristic cycle of $M_A(\beta)$.
From \cite{MMW05} it is known that the rank
$\mu^{F,\emptyset}_{A,0}(\beta)$ of $M_A(\beta)$ is upper
semicontinuous in the parameter $\beta$.  Theorem~\ref{47}
generalizes a weaker statement to all $L$ and $\tau\in\Phi_A^L$: the
multiplicity $\mu^{L,\tau}_{A,0}(\beta)$, $\tau\in\Phi^L_A$, is always
minimal at generic $\beta$.  We conjecture that
$\mu^{L,\tau}_{A,0}(\beta)$ is upper semicontinuous in $\beta$.

\begin{arxiv}
In the last section, we generalize results of the previous sections to a natural extension of $M_A(\beta)$ to the projective closure $(\PP_\CC^1)^n$.
\end{arxiv}

\begin{duke}
\begin{rmk*}\label{95}
In \cite{SW08} we introduce a natural extension of $M_A(\beta)$
to a $\calD$-module $\calM_A(\beta)$ on the product
$\prod_{j=1}^n\PP^1_\CC$ by noting that the generators $E$ and $I_A$
of the hypergeometric ideal extend to global differential
operators. On each standard affine patch one obtains a $D$-module that
may be viewed as a twisted $A$-hypergeometric system. 

Based on results in this article we describe in \cite{SW08}
the $L$-characteristic variety of $\calM_A(\beta)$ in terms of a
combinatorial gadget generalizing the $(A,L)$-umbrella, and we
determine the multiplicities in the $L$-characteristic Euler--Koszul
cycle for generic $\beta$. In order to discuss the slopes of
$\calM_A(\beta)$ (i.e., the slopes of $M_A(\beta)$ at infinity) one
then needs to understand for each point of $\prod_{j=1}^n\PP^1_\CC$
precisely which components of the $L$-characteristic variety pass
through its cotangent space; we state a conjecture to this account in
\cite[Conj.~5.18]{SW08}. The core of this issue seems to be related to the
question whether one can relax the basic condition on $A$ in this
article: that it be pointed. Without pointedness, thus placing the
origin in $T^*_0X$ outside the corresponding toric variety, three
issues come up: a) the definition of $\Delta^L_A$ must be changed; b)
Lemma~\ref{78} fails; c) the arguments in Section~\ref{cycle} need
adjustment. While we have some ideas for a) and our conjecture stated
in \cite{SW08} addresses b), c) is perhaps more involved, but
should have interesting answers.
\end{rmk*}
\end{duke}

\section{Filtrations on the toric ring}\label{spectoric}

\subsection{Torus action and toric ring}

\begin{ntn}\label{80}
By $\QQ_+$\index{Q@$\QQ_+$} we mean the nonnegative rational numbers
and we include $0$ in $\NN$.

Let $A=(a_{i,j})\in\ZZ^{d\times n}$\index{A@$A$} be an integer matrix
of rank $d$ whose columns
$\bolda_1,\dots,\bolda_n\in\ZZ^d$\index{a@$\bolda_i$} are nonzero.
We assume that $\NN A$ is a {\em positive} (see \cite[\S6.1]{BH93}) (or {\em pointed}, see \cite[\S1]{MMW05}) semigroup with $\ZZ A=\ZZ^d$.  
We write $\tau\subseteq A$ if $\tau$ is a subset of the
column set of $A$ and consider it both as a submatrix of $A$ and a
subset of the set of column indices $\{1,\dots,n\}$.  
Then the \emph{dimension} $\dim(\tau)$ is $\dim_\QQ(\QQ\tau)-1$. 
For a vector or
collection with index set $\{1,\dots,n\}$, a lower index $\tau$
denotes the subvector or subcollection with indices in
$\tau$\index{xt@$x_\tau$}\index{dt@$\del_\tau$}\index{tt@$t_\tau$}.
We abbreviate
$\ol\tau:=\{1,\dots,n\}\smallsetminus\tau$\index{t@$\ol\tau$}.
For any set $\tau$, its cardinality is denoted $|\tau|$. 

We shall frequently denote by $\bar C$ the Zariski closure of a set $C$.

For $\boldu\in\ZZ^n$ define $\boldu_+$ by\index{u@$\boldu_+,\boldu_-$}
$(\boldu_+)_j=\max(0,\boldu_j)$ and put $\boldu_-=\boldu_+-\boldu$.
For $\boldu,\boldv\in\NN^n$ write $\min(\boldu,\boldv)$ for the vector
whose $j$-th entry is $ \min(\boldu_j,\boldv_j)$. For any vector
$\boldu$ we mean by $\boldu>0$ that $\boldu$ is componentwise
positive: $u_i>0$ for all $i$.
\end{ntn}

The base space in this note is $X:=\Spec(\CC[x])=\CC^n$ where
$x:=x_1,\ldots,x_n$.  Let $R:=\CC[\del]$\index{R@$R$} be the
polynomial ring in $n$ variables
$\del:=\del_1,\dots,\del_n$\index{d@$\del$}.  Identifying $\del_i$
with the partial derivation $\del/\del x_i$, $\Spec(R)$ becomes the
conormal space $T^*_0X$ of $X$ at $0$.  The $d$-torus
$\TT:=(\CC^*)^d=\Spec(\CC[t_1^{\pm1},\dots,t_d^{\pm1}])$\index{T@$\TT$}
with coordinates $t:=t_1,\dots,t_d$\index{t@$t$} acts on $\Spec(R)$
\begin{equation}\label{16}
\xymat{(t,\xi)\ar@{|->}[r]&t\cdot\xi:=(t^{\bolda_1}\xi_1,\dots,t^{\bolda_n}\xi_n)}.
\end{equation}
This induces a $\ZZ^d$-grading on $R$ by $\deg(\del_i)=\bolda_i$\index{deg@$\deg$}.
For a $\ZZ^d$-graded $R$-module $N$, we denote by $\deg(N)$\index{degN@$\deg(N)$} the set of its $\ZZ^d$-degrees. 

\begin{dfn}\label{17}
We denote the orbit $\TT\cdot\xi$ through $\xi\in T^*_0X$ by $\Orb(\xi)$\index{Orb@$\Orb$}.
Let $\tau$ be a subset of columns of $A$.
We define $\boldone^\tau_A\in\{0,1\}^n$\index{1t@$\boldone^\tau_A$} by
\[
(\boldone^\tau_A)_j:=
\begin{cases}
1&\text{if}\quad\bolda_j\in\tau,\\
0&\text{if}\quad\bolda_j\notin\tau,
\end{cases}
\]
and denote by $O_A^\tau$\index{Ot@$O_A^\tau$} the orbit of $\boldone^\tau_A$.

The {\em toric ideal $I_\tau\subseteq R_\tau:=\CC[\del_\tau]$ of $\tau$}\index{It@$I_\tau$}\index{Rt@$R_\tau$} is the $\ZZ^d$-graded prime ideal generated by all $\square_\boldu=\del^{\boldu_+}-\del^{\boldu_-}$\index{$\square_\boldu$} where $\boldu\in\ZZ^{|\tau|}$ such that $\tau\cdot\boldu=0$. 
We set $I_A^\tau:=RI_\tau+J_\tau$\index{IAt@$I_A^\tau$} where $J_\tau$ is the $R$-ideal generated by $\{\del_i\mid i\not\in\tau\}$\index{Jt@$J_\tau$}.
Based on the following lemma, the {\em semigroup ring of $\tau$} is
\[
S_\tau:=\CC[\NN\tau]=\bigoplus_{\bolda\in\NN\tau}\CC\cdot t^\bolda=\sum_{\boldu\in\NN^{|\tau|}}\CC\cdot(\del_\tau^\boldu\modulo I_\tau)\subseteq S_A\index{St@$S_\tau$},
\]
The normalization of $S_A$ is the Cohen--Macaulay ring $\tilde
S_A=\bigoplus_{\bolda\in (\QQ_+ A\cap\ZZ^d)}\CC\cdot
t^\bolda$.\index{SA@$\tilde S_A$}
\end{dfn}

\begin{lem}\label{62}
For $\tau\subseteq A$, $I_A^\tau=\{f\in R\mid f(O_A^\tau)=0\}$ and $S_\tau=R_\tau/I_\tau=R/I_A^\tau$.\qed
\end{lem}

\subsection{$L$-filtration on the toric ring}

Let $L=(L_{\del_1},\ldots,L_{\del_n})\in\QQ^n$\index{L@$L$} be a {\em
weight vector}.  It induces an increasing
filtration $L$ of $\CC$-vector spaces on $R$ via $[\del^\boldu\in
L_kR]\Leftrightarrow[L\cdot\boldu\le k]$.  Note that $L$ has a
rational discrete index set. If $f\in L_kR\smallsetminus
L_{<k}R$ then $k=:\deg^L(f)$ is the \emph{$L$-degree} of $f$.  Let
$\sigma^L:R\to\gr^L(R)$\index{sigmaL@$\sigma^L$} be the {\em
$L$-symbol map} defined by $\sigma^L(f)=f\modulo L_{<k}R$ if
$\deg^L(f)=k$.  An element $f\in R$ is \emph{$L$-homogeneous} if 
$f=\sum_{L\cdot \boldu=k}f_\boldu \del^\boldu$ for some $k\in\QQ$
where $f_\boldu\in\CC$. 
By abuse of notation, we identify $R$ and
$\gr^L(R)$ via the $\CC$-linear isomorphism induced by $f\mapsto
\sigma^L(f)$ for $L$-homogeneous $f$ in right-normal form.

For any $\tau\subseteq A$, $L$ induces a filtration on $S_\tau$ by 
\[
L_k S_\tau:=\sum_{\deg^L(\del_\tau^\boldu)\le
  k}\CC\cdot(\del_\tau^\boldu\modulo I_\tau).
\]
With $I_\tau^L$ denoting the $\ZZ^d$-graded ideal $\gr^L(I_\tau)$\index{ItL@$I_\tau^L$}, we abbreviate
\[
S_\tau^L:=\gr^L(S_\tau)\cong R/I_\tau^L\index{StL@$S_\tau^L$}.
\]

The following is a mild extension of \cite[Cor.~4.4]{Stu96} (see also \cite[Lem.~4.11]{Ado94}).

\begin{lem}\label{5}
One has the identity $I^L_A=R\ideal{\sigma^L(\square_\boldu)\mid\boldu\in\ZZ^n,A\cdot\boldu=0}$.
\end{lem}

\begin{duke}
\begin{proof}
Form a matrix $B$ by adding the row $L$ on top, and then a
column $(1,0,\ldots,0)$ on the left of $A$. Then $I_B$ is the
$L$-homogenization of $I_A$ relative to the new variable $\del_0$. By
\cite[Cor.~4.4]{Stu96}, for any term order that eliminates
$\del_0$, $I_B$ has a Gr\"obner basis that consists of
$L$-homogenizations of binomials $\del^{\boldu_+}-\del^{\boldu_-}$ with
$A\cdot\boldu=0$. Dehomogenization leads to an $L$-Gr\"obner basis for
$I_A$ and the claim follows.
\end{proof}
\end{duke}

\begin{arxiv}
\begin{proof}
Generators of $I^L_A$ are the $L$-symbols of a reduced Gr\"obner basis
of $I_A$ and can be computed by Buchberger's algorithm.
But each S-pair and reduction step of the algorithm preserves the
generators $\{\square_\boldu\mid\boldu\in\ZZ^n\cap\ker(A)\}$:
If $\boldu'+\boldu_+=\boldv'+\boldv_+$ where
$\boldu,\boldv\in\ZZ^n\cap\ker(A)$ and $\boldu',\boldv'\in\NN^n$ then
\[
A(\boldv'+\boldv_-)=A(\boldv'+\boldv_+)=A(\boldu'+\boldu_+)=A(\boldu'+\boldu_-).
\]
and
\[
\del^{\boldu'}\square_\boldu-\del^{\boldv'}\square_\boldv=
\del^{\boldv'+\boldv_-}-\del^{\boldu'+\boldu_-}=
\del^{\min(\boldu'+\boldu_-,\boldv'+\boldv_-)}
\square_{\boldv'+\boldv_--\boldu'-\boldu_-}.
\]
Note that this argument is valid also for $L\not>0$ where Buchberger's
algorithm with (de-)homogenization is used.
\end{proof}
\end{arxiv}

In this section we study the geometry of $S^L_A$ in terms of $A$ and $L$. 
Since the Gr\"obner fan of any $R$-ideal is defined over $\QQ$, the study of real weights can be reduced to the case $L\in\QQ^n$ (see \cite[\S7.4]{MS05} for a discussion in the case where $L>0$ defines a term order).

\subsection{The $(A,L)$-umbrella}

We consider the embedding of the affine space $\QQ^d\supseteq A$ into the rational projective $d$-space 
\[
\PP^d_\QQ=\PP_\QQ(\QQ^d\times \QQ)
\]
via the map $q\mapsto (q:1)$. Denote
$\infty:=\PP_\QQ^d\smallsetminus\QQ^d$ the hyperplane at infinity.

In $\PP^d_\QQ$, any two distinct points $\bolda,\boldb\in\PP^d_\QQ$
are joined by two line segments.  If $H$ is a hyperplane in
$\PP^d_\QQ$ containing neither $\bolda$ nor $\boldb$ then there is a
unique line segment joining the points and not meeting $H$.  This is
exactly the convex hull of $\bolda$ and $\boldb$ in the affine space
$\PP^d_\QQ\smallsetminus H$.

\begin{dfn}
Let $H\subseteq \PP^d_\QQ$ be a hyperplane and let
$U_H:=\PP^d_\QQ\smallsetminus H$ be its
complement.  For $B\subseteq U_H$, the {\em
convex hull} of $B$ relative to $H$ is the set
$\conv_H(B)$ defined as the convex hull of $B$ in the affine space
$U_H$.
\end{dfn}

Note that, for varying $H$, $\conv_H(B)$ changes exactly when $H$ is
moved through a point of $B$. 
Within $\QQ^d$, elements of convex hulls are linear combinations with nonnegative coefficients that add to unity. 
Convex hulls relative to $H\neq \infty$, with coordinates from $\QQ^d\subseteq\PP_\QQ^d$, obey slightly different rules.
Let $h\in\Hom_\QQ(\QQ^d,\QQ)$ be a linear form such that $H$ is the closure of $h^{-1}(0)$ in $\PP_\QQ^d$.
The line through $0$ and $\bolda\in\QQ^d\smallsetminus\{0\}$ meets
$\infty$ in a point that we denote $\bolda/0$.  For $q\in \QQ$ let
$\sign(q)$ be the usual signum function:
\[
\sign(q)=\begin{cases} -1&\text{if $q<0$;}\\
                        0&\text{if $q=0$;}\\
                        1&\text{if $q>0$.}\end{cases}
\]

\begin{lem}\label{25}
Let $B=\{\boldb_1,\ldots,\boldb_n\}\subseteq U_H$ with $B\cap\infty=\{\boldb_{m+1},\ldots,\boldb_n\}$ and pick $\{\boldb'_{m+1},\ldots,\boldb'_n\}\subseteq\QQ^d\cap U_H$ such that $\boldb_j=\boldb'_j/0$ for $m<j\le n$. 
If $\boldb\in\conv_H(B)\cap\QQ^d$ then in coordinates of $\QQ^d$ there is an equation $\boldb=\sum_{j=1}^m\eps_j\boldb_j+\sum_{j=m+1}^n\eps'_j\boldb'_j$ where $\sum_{j=1}^m\eps_j=1$, $\sign(\eps_jh(\boldb_j))=\sign(h(\boldb))$ for all $j$ with $\eps_j\not=0$, and $\sign(\eps'_jh(\boldb'_j))=\sign(h(\boldb))$ for all $j$ with $\eps'_j\not=0$.
\end{lem}

\begin{proof}
Let $B_+$ (resp.~$B_-$) be the subsets of $B\smallsetminus\infty$ on
which $h$ evaluates positively (resp.~negatively), and put
$B_\infty=B\cap\infty$, $B'=\{\boldb'_{m+1},\ldots,\boldb'_n\}$. A
general element $\boldb\in\conv_H(B)$ is the convex combination of three
points: $\boldb_+\in\conv_H(B_+)$, $\boldb_-\in \conv_H(B_-)$, and
$\boldb_\infty\in \conv_H(B_\infty)$.

Clearly, $\boldb_+=\sum_{\boldb_j\in B_+}\eps_j\boldb_j$ where
$\eps_j\geq 0$ and $\sum_{\boldb_j\in B_+}\eps_j=1$. A similar statement
holds for $\boldb_-$. Now
$\conv_H(\boldb_+,\boldb_-)\cap \QQ^d$ is the union of rays
$\{\lambda_+\boldb_+ +\lambda_-\boldb_-\mid \lambda_++\lambda_-=1,\,
\lambda_+\lambda_-\leq 0\}$. Thus, if
$\boldb_0\in\conv_H(\boldb_+,\boldb_-)$ then
$\sign(\lambda_+h(\boldb_+))=\sign(\boldb_0)=\sign(\lambda_-h(\boldb_-))$. 

It suffices to show the lemma if $h(B')>0$. 
Then $\boldb_\infty$ is of the form $\boldb'/0$ with $\boldb'\in\conv_H(B')$, and points of $\conv_H(B')$ are of the form $\sum_{\boldb'_j\in B'}\eps_j\boldb'_j$ with $\sum_{\boldb'_j\in B'}\eps_j=1$ and all $\eps_j\geq 0$.
We are thus reduced to considering $\conv_H(\{\boldb_0,\boldb_\infty\})$ with $\boldb_\infty=\boldb'/0$ and $h(\boldb')>0$. If $h(\boldb_0)>0$ then $\conv_H(\{\boldb_0,\boldb_\infty\})$ is the ray $\{\boldb_0+\lambda'\boldb'\mid\lambda'\geq 0\}$, while if $h(\boldb_0)<0$ then it is the ray $\{\boldb_0+\lambda'\boldb'\mid\lambda\leq 0\}$. 
The condition on $\lambda'$ can, if $\lambda'\ne 0$, be packaged as $\sign(\lambda'h(\boldb'))=\sign(h(\boldb_0))$.
The lemma follows.
\end{proof}

We view the columns $\bolda_1,\dots,\bolda_n\in\ZZ^d$ of $A$ as points
in $\QQ^d=\PP_\QQ^d\setminus\infty$.  By assumption, $\NN A$ is positive and
hence $h\in\Hom_\QQ(\QQ^d,\QQ)$ can be chosen such that
$h(\bolda_j)>0$ for all $j$.  For any $\lambda\in\QQ$, set
$H_\lambda:=h^{-1}(-\lambda)$ and $U_\lambda:=\PP_\QQ^d\smallsetminus
H_\lambda$.

\begin{dfn}\label{35}
Choose $h\in\Hom_\QQ(\QQ^d,\QQ)$ such that $h(\bolda_j)>0$ for all $j$, and let $\eps$ be such that $0<\eps<|h(\bolda_j)/L_{\del_j}|$ whenever $L_{\del_j}\ne 0$.
We set $\bolda_j^L:=\bolda_j/L_{\del_j}$\index{aL@$\bolda_j^L$} and call $\Delta^L_A:=\conv_{H_\eps}(\{0,\bolda_1^L,\dots,\bolda_n^L\})\subseteq\PP^d_\QQ$\index{DLA@$\Delta^L_A$} the {\em $(A,L)$-polyhedron}.

Let the {\em $(A,L)$-umbrella} be the set $\Phi_A^L$\index{PAL@$\Phi_A^L$} of faces of $\Delta^L_A$ which do not contain $0$. 
In particular, $\Phi^L_A$ contains the empty face.
Whenever it suits us, we identify $\tau\in\Phi_A^L$\index{t@$\tau$} with $\{j\mid\bolda^L_j\in\tau\}$, or with $\{\bolda_j\mid\bolda^L_j\in\tau\}$, or with the corresponding submatrix of $A$.
By $\Phi_A^{L,k}\subseteq\Phi_A^L$\index{PALk@$\Phi_A^{L,k}$}, we denote the subset of faces of dimension $k$.

By $\Gamma_A^L:=\bigcup\Phi_A^L$\index{GAL@$\Gamma_A^L$} (resp.\
$\Gamma_A^{L,k}:=\bigcup\Phi_A^{L,k}$\index{GALk@$\Gamma_A^{L,k}$}),
we denote the underlying point set of $\Phi_A^L$ (resp.\
$\Phi_A^{L,k}$). Note that $\Gamma^L_A$ is a piecewise linear manifold
with boundary, homeomorphic to the $(d-1)$-disk. If $\bolda\in
\Delta^L_A\smallsetminus\{0\}$ then the line through $0$ and $\bolda$
meets $\Gamma^L_A$ in $\Gamma^L_A(\bolda)$.\index{Ga@$\Gamma(\bolda)$}

The matrix $A$ is called {\em $L$-homogeneous} if all $\bolda^L_j$ lie
on a common hyperplane of $\PP_\QQ^d$.  Every $A$ is
$\boldzero$-homogeneous and we call
$\Phi_A^\boldzero$\index{PA0@$\Phi_A^\boldzero$} the \emph{$A$-umbrella}.
If $A$ is $L$-homogeneous then $I^L_A=I_A$ under the identification of
$\gr^L(R)$ with $R$. Note that $\Phi^\boldzero_A$ can be
identified with the lattice of nonempty faces of the polyhedral cone
$\QQ_+A$ via the face lattice of a cross-section $\QQ_+A\cap
h^{-1}(1)$.
\end{dfn}

\begin{figure}[ht]
\caption{More $(A,L)$-umbrellas. (Shaded $\Delta^L_A$ with fat boundary
$\Phi^L_A$.)}\label{90}
\begin{center}
\setlength{\unitlength}{0.8mm}
\begin{picture}(150,65)(0,-25)
\put(0,0){\shade\path(35,38.75)(0,30)(0,0)(10,0)(40,7.5)
\put(0,0){\thicklines\path(0,30)(35,38.75)\dottedline{2}(35,38.75)(40,40)}
\put(0,0){\thicklines\path(10,0)(40,7.5)\dottedline{2}(40,7.5)(45,8.625)}
\put(5,-15){\makebox(0,0)[bl]{$L=(1,1,1,0)$}}}
\put(105,0){\shade\path(5,35)(0,30)(0,0)(10,0)(40,6)
\put(0,0){\thicklines\path(0,30)(5,35)\dottedline{2}(5,35)(10,40)}
\put(0,0){\thicklines\path(10,0)(40,6)\dottedline{2}(40,6)(45,7)}
\shade\path(-50,-12)(-40,-10)(-45,-15)
\put(0,0){\thicklines\path(-40,-10)(-50,-12)\dottedline{2}(-50,-12)(-55,-13.25)}
\put(0,0){\thicklines\path(-40,-10)(-45,-15)\dottedline{2}(-45,-15)(-50,-20)}
\dottedline{1}(-40,-10)(0,30)
\dottedline{1}(-40,-10)(10,0)
\dottedline{1}(-40,-10)(00,0)
\put(0,0){\path(0,0)(-55,0)}
\put(0,0){\path(0,0)(0,-15)}
\multiput(0,0)(-10,0){6}{\path(0,-1)(0,1)}
\multiput(0,0)(0,-10){2}{\path(-1,0)(1,0)}
\put(-40,-10){\circle*{1}\put(0,-1){\makebox(0,0)[t]{$\bolda_4^L$}}}
\put(5,-15){\makebox(0,0)[bl]{$L=(1,1,1,-1)$}}}
\multiput(0,0)(105,0){2}{
\put(0,0){\vector(1,0){45}}
\put(0,0){\vector(0,1){35}}
\multiput(0,0)(10,0){5}{\path(0,-1)(0,1)}
\multiput(0,0)(0,10){4}{\path(-1,0)(1,0)}
\put(0,30){\circle*{1}\put(1,-1){\makebox(0,0)[tl]{$\bolda_1$}}}
\put(10,0){\circle*{1}\put(0,-1){\makebox(0,0)[t]{$\bolda_2$}}}
\put(10,20){\circle*{1}\put(1,0){\makebox(0,0)[l]{$\bolda_3$}}}
\put(40,10){\circle*{1}\put(0,1){\makebox(0,0)[b]{$\bolda_4$}}}
\dottedline{1}(0,0)(45,11.25)}
\end{picture}
\end{center}
\end{figure}
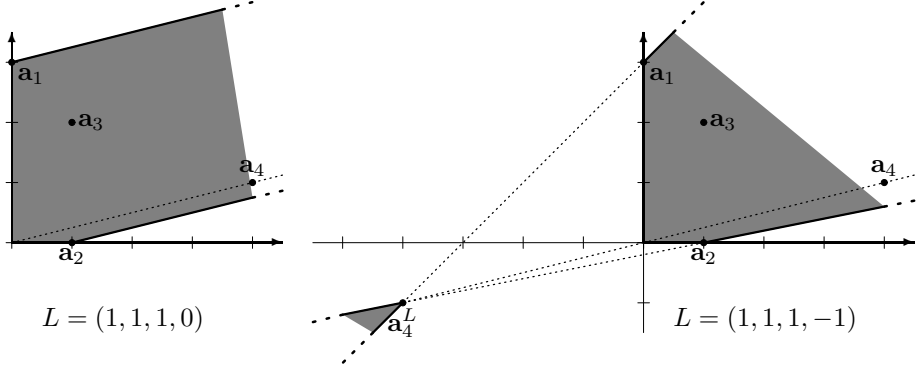

Figure~\ref{90} shows the intersection with $\QQ^d$ of the $(A,L)$-umbrella for $L=(1,1,1,t)$ and $t\le 0$ in the example from the introduction, $A=\begin{pmatrix}0&1&1&4\\3&0&2&1\end{pmatrix}$.

\subsection{Monomials in the graded toric ideal}

The following result generalizes Lemma 3.1 and 3.2 in \cite{Ado94}.

\begin{lem}\label{1}
If $\bolda_{k_1}^L,\dots,\bolda_{k_m}^L\in\Gamma_A^L$ do not lie in a common $\tau\in\Phi_A^{L,d-1}$ then $\del_{k_1}\cdots\del_{k_m}\in\sqrt{I^L_A}$.
\end{lem}

\begin{proof}
We write $v_j$ for $L_{\del_j}$.

Assume first that $v_{k_1}\ne0$ and $\Gamma_A^L(\bolda_{k_1})=\bolda_{k_1}/0\in\infty$.
In particular, $v_{k_1}>0$ by Definition~\ref{35}.
The polyhedron $\Delta^L_A\cap\infty$ is the convex hull of the points $\bolda^L_i=\bolda_i/0$ with $v_i=0$, and of the intersection points $(\bolda^L_i-\bolda^L_j)/0$ with $\infty$ of line segments from $\bolda^L_i$ to $\bolda^L_j$ with $v_i>0>v_j$. 
Therefore
\[
\bolda_{k_1}^L=\sum_{v_i=0}\eta_i\bolda_i+\sum_{v_i>0>v_j}\eta_{i,j}\left(\bolda^L_i-\bolda^L_j\right)
\]
for some $0\le\eta_i,\eta_{i,j}\in\QQ$.
This equality gives rise to an element
\[
\square=\del_{k_1}^{s/v_k}-\prod_{v_i=0}\del_i^{s\eta_i}\prod_{v_i>0>v_j}\del_i^{s\eta_{i,j}/v_i}\del_j^{-s\eta_{i,j}/v_j}\in I_A
\]
where $s\in\NN$ is chosen to clear all denominators.
The $L$-degree of the left monomial is positive while that of the right one is zero.
Thus $\sigma^L(\square)=\del_{k_1}^{s/v_k}\in\gr^L(I_A)$ and the claim follows in this case.

We now keep the assumption $v_{k_1}\ne0$ but assume that $\Gamma_A^L(\bolda_{k_1})\not\in\infty$.
By hypothesis there is $\bolda\in\conv(\bolda_{k_1}^L,\dots,\bolda_{k_m}^L)$ in the interior of $\Delta_A^L$.
By continuity of the function $\bolda\mapsto \Gamma_A^L(\bolda)$, choosing $\bolda$ sufficiently close to $\bolda_{k_1}^L$ implies that $\Gamma_A^L(\bolda)\not\in\infty$.
As $\bolda$ lies in the interior of $\Delta_A^L$ and since $\Gamma^L_A(\bolda)\notin \infty$, there is $0\ne t\in\QQ$ with $\Gamma^L_A(\bolda)=t\bolda\in\Gamma_A^L$ and hence either $\sign(t h(\bolda))=1$ and $t>1$, or $\sign(t h(\bolda))=-1$ and $0<t<1$. 
In either case, $t\sign(th(\bolda))>\sign(th(\bolda))$.
Using Lemma~\ref{25}, we can write
\begin{equation}\label{103}
\bolda=\sum_{v_{k_j}\ne0}\eps_{k_j}\bolda_{k_j}^L+\sum_{v_{k_j}=0}\eps_{k_j}\bolda_{k_j}
\end{equation}
for some $\eps_{k_j}\in\QQ$ with $\sum_{v_{k_j}\ne0}\eps_{k_j}=1$, $\sign(\eps_{k_j}/v_{k_j})=\sign(h(\bolda))$ if $\eps_{k_j}v_{k_j}\ne0$, and $\sign(\eps_{k_j})=\sign(h(\bolda))$ if $\eps_{k_j}\ne0=v_{k_j}$.
Again by Lemma~\ref{25} and Definition~\ref{35}, we can write
\begin{equation}\label{104}
t\bolda=\sum_{v_i\ne0}\eta_i\bolda_i^L+\sum_{v_i=0}\eta_i\bolda_i
\end{equation}
for some $\eta_i\in\QQ$ where $\sum_{v_i\ne0}\eta_i=1$, $\sign(\eta_i/v_i)=\sign(th(\bolda))$ if $\eta_iv_i\ne0$, and $\sign(\eta_i)=\sign(th(\bolda))$ if $\eta_i\ne0=v_i$.
Combining equations \eqref{103} and \eqref{104} we find an element
\begin{equation}\label{105}
\square=\prod_{v_{k_j}\ne0}\del_{k_j}^{s|t\eps_{k_j}/v_{k_j}|}\prod_{v_{k_j}=0}\del_{k_j}^{s|t\eps_{k_j}|}-\prod_{v_i\ne0}\del_i^{s|\eta_i/v_i|}\prod_{v_i=0}\del_i^{s|\eta_i|}\in I_A
\end{equation}
where $s\in\NN$ chosen to clear all denominators.
From
\begin{gather}\label{109}
\deg^L\bigl(\prod_{v_{k_j}\ne0}\del_{k_j}^{s|t\eps_{k_j}/v_{k_j}|}\prod_{v_{k_j}=0}\del_{k_j}^{s|t\eps_{k_j}|}\bigr)
=s|t|\sum_{v_{k_j}\ne0}\eps_{k_j}\sign(\eps_{k_j}/v_{k_j})\\
\nonumber=st\sum_{v_{k_j}\ne0}\eps_{k_j}\sign(th(\bolda))
=st\sign(th(\bolda))
>s\sign(th(\bolda))\\
\nonumber=s\sum_{v_i\ne0}\sign(\eta_i/v_i)\eta_i
=s\sum_{v_i\ne0}|\eta_i/v_i|v_i
=\deg^L\Bigl(\prod_{v_i\ne0}\del_i^{s|\eta_i/v_i|}\prod_{v_i=0}\del_i^{s|\eta_i|}\Bigr),
\end{gather}
we conclude that the $L$-leading term of $\square$ is the left of the two monomials in \eqref{105}. 
The claim follows in the case where at least one $v_{k_j}$ is nonzero.

Suppose finally that $v_{k_j}=0$ and hence $\bolda_{k_j}^L\in\infty$ for all $j=1,\dots,m$.
By assumption and Definition~\ref{35} we can pick an element $\bolda'\in\conv(\bolda_{k_1}^L,\ldots,\bolda_{k_m}^L)\smallsetminus\partial\Delta_A^L\subseteq\conv(\bolda_{k_1}^L,\ldots,\bolda_{k_m}^L)\smallsetminus\Gamma_A^L$.
Then $\bolda'=\bolda/0$ with $\bolda\in\conv(\bolda_{k_1},\ldots,\bolda_{k_m})$. 
It follows that there is an equation of type \eqref{103} with conditions as indicated there.
By construction, $\Gamma^L_A(\bolda')\neq\bolda'\in\infty$ and so $\Gamma^L_A(\bolda')\not\in\infty$. 
Hence, there is a positive $t\in\QQ$ with $\Gamma^L_A(\bolda')=t\bolda$ and so there is an equation of type \eqref{104} with conditions as indicated there.
As $t\sign(th(\bolda))=-t>0>-1=\sign(th(\bolda))$ we get an equation of type \eqref{109}, and the claim follows as in the previous case.
\end{proof}

\subsection{Homogeneity in the graded toric ideal}

\begin{dfn}
Let $\tau$ be a set of columns of $A$.
For $\boldu\in\ZZ^n$, we write
$\supp(\boldu)\subseteq\tau$\index{supp@$\supp$} if $\boldu_i\not =0$ implies $\bolda_i\in\tau$. 
For $f\in R$, we write $\supp(f)\subseteq\tau$ if $\supp(\boldu)\subseteq\tau$ for all monomials $\del^\boldu$ of $f$.
\end{dfn}

\begin{lem}\label{6}
Let $\tau\in\Phi_A^L$ and pick $\boldu\in\ZZ^n$ such that $A\cdot\boldu=0$.  
\begin{enumerate}
\item\label{6a} If $\supp(\square_\boldu)\subseteq\tau$ then $\sigma^L(\square_\boldu)=\square_\boldu$. In particular, the toric ideal $I_\tau$ is $L$-homogeneous.
\item\label{6b} If $\supp(\del^{\boldu_\pm})\subseteq\tau$ and $\supp(\del^{\boldu_\mp})\not\subseteq\tau$ then $\sigma^L(\square_\boldu)=\mp\del^{\boldu_\mp}$.
\end{enumerate}
\end{lem}

\begin{proof}
We write $v_j$ for $L_{\del_j}$.

Consider first the case where the facet $\tau\in\Phi^{L,d-1}_A$ lies
entirely in $\infty$. Then $\bolda_i^L\in\tau$ implies $v_i=0$ and
hence $\supp(\square_\boldu)\subseteq\tau$ implies that
$\square_\boldu$ is $L$-homogeneous of degree zero.  Suppose
$\supp(\boldu_+)\subseteq\tau$ but $\supp(\boldu_-)\not\subseteq\tau$.
As $\tau\subseteq\infty$ is a facet of $\Delta^L_A$ and by
Definition~\ref{35}, the interior of $\Delta^L_A$ meets neither
$\infty$ nor $H_\eps$. Hence $\Delta^L_A$ is completely contained in
the $\PP^d_\QQ$-closure of one of the regions $\{t\in\QQ^d\mid
h(t)\geq-\eps\}$ and $\{t\in\QQ^d\mid h(t)\leq-\eps\}$. Since
$0\in\Delta^L_A$, it must be the former. Moreover, by definition of
$\eps$, $h(\bolda^L_j) \not\in [-\eps,0]$ for all $v_j\not =0$, whence
$v_j\geq0$ in all cases.  Thus $\deg^L(\del^{\boldu_+})=0$ and
$\deg^L(\del^{\boldu_-})>0$, hence
$\sigma^L(\square_\boldu)=-\del^{\boldu_-}$.  The case
$\supp(\boldu_-)\subseteq\tau\not\supseteq\supp(\boldu_+)$ is similar.

Now let $\tau\in\Phi^L_A$ be not contained in $\infty$, or let
$\tau\in\Phi^L_A\smallsetminus \Phi^{L,d-1}_A$ be a nonfacet face
contained in $\infty$.  Then there is a linear form
$h_\tau\in\Hom_\QQ(\QQ^d,\QQ)$ such that the closure of
$h_\tau^{-1}(1)$ in $\PP^d_\QQ$ meets $\Delta^L_A$ precisely in
$\tau$, and for which 
$\sign(v_i)h_\tau(\bolda^L_i)\le\sign(v_i)$ whenever
$v_i\ne0$.

For all $\bolda^L_i\in\tau$ with $v_i=0$ the line through $0$ and
$\bolda_i$ meets $\tau$ in $\bolda^L_i\in\infty$.  This means that
$\bolda_i$ is parallel to $h_\tau^{-1}(1)$ and so $h_\tau(\bolda_i)=0$.  Thus,
$A\cdot\boldu=0$ implies that
\[
0=h_\tau(A\cdot\boldu)=h_\tau\left(\sum_{i=1}^nu_i\bolda_i\right)
=\sum_{u_i>0\ne v_i}u_iv_ih_\tau(\bolda^L_i)-\sum_{u_i<0\ne v_i}(-u_i)v_ih_\tau(\bolda^L_i).
\]
If $\supp(\square_\boldu)\subseteq\tau$ and $u_iv_i\ne0$ then $h_\tau(\bolda_i^L)=1$ and hence $\deg^L(\del^{\boldu_+})=\deg^L(\del^{\boldu_-})$.

Now suppose $\supp(\del^{\boldu_+})\subseteq\tau\not\supseteq\supp(\del^{\boldu_-})$.  Hence $h_\tau(\bolda_i^L)=1$ if $u_i>0\ne v_i$, $\sign(v_i)h_\tau(\bolda_i^L)\le\sign(v_i)$ if $u_i<0\ne v_i$, and $\sign(v_i)h_\tau(\bolda_i^L)<\sign(v_i)$ for at least one $u_i<0$.  
Thus,
\begin{align*}
\deg^L(\del^{\boldu_+})=
\sum_{u_i>0}v_iu_i&=\sum_{u_i>0\ne v_i}v_iu_ih_\tau(\bolda_i^L)\\
&=\sum_{u_i<0\ne v_i}v_i(-u_i)h_\tau(\bolda_i^L)<\sum_{u_i<0\ne v_i}v_i(-u_i)
=\deg^L(\del^{\boldu_-}).
\end{align*}
The case $\supp(\boldu_-)\subseteq\tau\not\supseteq\supp(\boldu_+)$ is similar.
\end{proof}

\subsection{Minimal associated prime ideals}

We now identify the components of the $L$-graded toric ring $S_A^L$.

\begin{dfn}
For $\boldu\in\ZZ^n$ let 
\[
\tau^L_\boldu=\bigcap_{\supp(\boldu)\subseteq\tau\in\Phi^L_A}\tau
\]
be the smallest element of $\Phi^L_A$ containing $\supp(\boldu)$ and put $\tau^L_\boldu=A$ if there is none.\index{tu@$\tau_\boldu$}
Note that $[\del^\boldu\in J_\tau]\Rightarrow[\tau^L_\boldu\not\subseteq\tau]$ for $\tau\subseteq A$ and $\boldu\in\NN^n$. 
\end{dfn}

Recall that $I_A$, $I_A^\tau$, and $J_\tau$ are $\ZZ^d$-graded prime ideals.

\begin{lem}\label{7}
Let $\tilde I_A^L\subseteq R$ be generated by all elements of the following types:
\begin{enumerate}
\item\label{7a} $\del_{k_1}\cdots\del_{k_m}$ where $\bolda_{k_1}^L,\dots,\bolda_{k_m}^L$ do not lie in a common facet $\tau\in\Phi_A^{L,d-1}$; 
\item\label{7b} $\square_\boldu$ where $\boldu\in\ZZ^n$, $A\cdot\boldu=0$, and $\tau^L_\boldu\not=A$.
\end{enumerate}
Then $\tilde I_A^L=\bigcap_{\tau\in\Phi_A^{L,d-1}}I_A^\tau$.
\end{lem}

\begin{proof}
Let $I=\bigcap_{\tau\in\Phi_A^{L,d-1}}I_A^\tau$. 
We show first that $\tilde I^L_A\subseteq I$. 
Let $I_0$ denote the ideal of $R$ generated by the elements from \eqref{7a}.
Then clearly $I_0\subseteq J_\tau\subseteq  I^\tau_A$ for any $\tau\in\Phi^{L,d-1}_A$.
Pick $\square_\boldu$ as in \eqref{7b}. 
By part~\eqref{6a} of Lemma~\ref{6}, $\sigma^L(\square_\boldu)=\square_\boldu$.
Let $\tau\in\Phi_A^{L,d-1}$ and suppose $\square_\boldu\notin J_\tau$. 
Then (without loss of generality) $\supp(\boldu_+)\subseteq\tau$. 
By part~\eqref{6b} of Lemma~\ref{6}, since $\sigma^L(\square_\boldu)=\square_\boldu$, $\supp(\boldu_-)\subseteq\tau$ as well. 
Hence $\square_\boldu\in I_\tau\subseteq I^\tau_A$ and so $\tilde I_A^L\subseteq I$.

For the converse inclusion suppose now $m\in I$.
We write
$m=m_0+\sum_{\tau\in\Phi^L_A}m_\tau$ where $m_0\in I_0$, and 
\[
m_\tau=\sum_{\tau^L_\boldu=\tau}c_\boldu\del^\boldu\in\CC[\del_\tau]
\]
collects the monomials in $m$ minimally supported in $\tau$. 
Since $I_0\subseteq \tilde I^L_A\subseteq I$ we may assume that $m_0=0$. 
Now pick any $\tau\in\Phi^{L}_A$.
Then $m=\hat m_\tau+m_{\ol\tau}$ where $\hat m_\tau=\sum_{\tau'\subseteq\tau}m_{\tau'}$ are the terms in $m$ supported in $\tau$. 
Since $m_{\ol\tau}\in J_\tau\subseteq I^\tau_A$ and $m\in I\subseteq I_A^\tau$, $\hat m_\tau=m-m_{\ol\tau}\in I_A^\tau\cap\CC[\del_\tau]=I_\tau\subseteq\tilde I^L_A$ in view of \eqref{7b}. 
Since every $\hat m_\tau$ is in $\tilde I^L_A$, so is every $m_\tau$, and hence $m\in\tilde I^L_A$ as well.
\end{proof}

\begin{lem}\label{3}
The radical ideal $\tilde I^L_A$ from Lemma~\ref{7} equals $\sqrt{I^L_A}$.
\end{lem}

\begin{proof}
By Lemma~\ref{1}, the elements in \ref{7}.\eqref{7a} are in $\sqrt{I^L_A}$.
By part~\eqref{6a} of Lemma~\ref{6}, all elements from \ref{7}.\eqref{7b} are in $I^L_A$. 
Hence $\tilde I^L_A\subseteq\sqrt{I^L_A}$.
 
By Lemmas~\ref{5} and \ref{7}, it suffices to show conversely that, for any $\boldu\in\ZZ^n$ with $A\cdot\boldu=0$, $\sigma^L(\square_\boldu)\in I^\tau_A$ for all $\tau\in\Phi^{L,d-1}_A$. 
Pick such $\boldu$ and let $\tau\in\Phi^{L,d-1}_A$. 
If $\tau^L_\boldu\subseteq\tau$ then by part~\eqref{6a} of Lemma~\ref{6}, $\sigma^L(\square_\boldu)=\square_\boldu\in I_\tau\subseteq I^\tau_A$. 
If $\tau^L_{\boldu_+}\subseteq\tau$ but $\tau^L_{\boldu_-}\not\subseteq\tau$ then by part~\eqref{6b} of Lemma~\ref{6}, $\sigma^L(\square_\boldu)=-\del^{\boldu_-}\in J_\tau\subseteq I^\tau_A$.
Similarly, $\sigma^L(\square_\boldu)\in I^\tau_A$ if $\tau^L_{\boldu_-}\subseteq\tau$ but $\tau^L_{\boldu_+}\not\subseteq\tau$. 
Finally, if $\tau^L_{\boldu_+}\not\subseteq\tau\not\supseteq\tau^L_{\boldu_-}$ then $\sigma^L(\square_\boldu)\in J_\tau\subseteq I^\tau_A$.
\end{proof}

The following consequence of Lemmas~\ref{62}, \ref{7}, and \ref{3} generalizes \cite[Lem.~3.2]{Ado94}.

\begin{thm}\label{13}
The set of $\ZZ^d$-graded prime ideals of $R$ containing $I^L_A$ equals $\{I_A^\tau\mid\tau\in\Phi_A^L\}$ and hence the $(A,L)$-umbrella encodes the geometry of $S^L_A$:
\[
\Spec(S^L_A)=\Var(I_A^L)=\bigcup_{\tau\in\Phi_A^{L,d-1}}\bar O_A^\tau=\bigsqcup_{\tau\in\Phi_A^L}O_A^\tau.
\]
Adjacencies of orbit strata correspond to inclusions in the $(A,L)$-umbrella:
\[\pushQED{\qed}
\left[ O_A^{\tau'}\subseteq\bar O_A^\tau\right]\quad\Leftrightarrow\quad\left[\tau'\subseteq\tau\right].\qedhere
\]
\end{thm}

In particular, Theorem~\ref{13} identifies the $\ZZ^d$-graded prime ideals containing $I_A=I^\boldzero_A$ with the elements of $\Phi^\boldzero_A$ from Definition~\ref{35}.

\begin{arxiv}
\begin{rmk}
Geometrically, $\Spec(S^L_A)$ reflects certain asymptotics of $\Spec(S_A)$.
We make this precise in the case where $L_{\del_j}>0$ for all $j$.
We may assume that the components of $L$ are coprime positive integers.
Then $X=\CC^n$ can be considered as the smooth chart defined by $\del_0\ne0$ in the weighted projective space $\PP^n_{(1,L)}$\index{PnL@$\PP^n_L$} where $(1,L):=(1,L_{\del_1},\dots,L_{\del_n})$, and $\Spec(S_A)\subseteq X$.
The closure $Z$ of $\Spec(S_A)$ in $\PP^n_{(1,L)}$ is defined by the $L$-homogenization $I_A^L(\del_0)$ of $I_A$ with respect to $\del_0$.
The ideal of the intersection of $Z$ with the weighted projective space $\PP^{n-1}_L=\Var(\del_0)\subseteq\PP^n_{(1,L)}$ is $I_A^L(\del_0)_{\vert\del_0=0}=I_A^L$.
Thus $Z\cap\PP^{n-1}_L=\Proj(S_A^L)\subseteq\PP^{n-1}_L$.
\end{rmk}
\end{arxiv}

\subsection{Index formula for multiplicities}

By Theorem~\ref{13}, there is a composition chain of $\ZZ^d$-graded $R$-modules
\begin{equation}\label{56}
0=N_0\subsetneq N_1\subsetneq\cdots\subsetneq N_{l-1}\subsetneq N_l=S_A^L
\end{equation}
with $N_i/N_{i-1}\cong S_{\tau_i}(-\boldu_i)$ for some $\tau_i\in\Phi_A^L$ and $\boldu_i\in\ZZ^d$.
As $\ZZ^d$-graded vector spaces $S^L_A=S_A$ and the $\ZZ^d$-graded Hilbert function of both rings has values in $\{0,1\}$.
Thus, the composition chain~\eqref{56} induces a partition of $\ZZ^d$-degrees
\begin{equation}\label{92}
\deg(S_A^L)=\deg(S_A)=\NN A=\bigsqcup_{i=1}^l(\boldu_i+\NN\tau_i)=\bigsqcup_{i=1}^l\deg(S_{\tau_i}(-\boldu_i)).
\end{equation}

\begin{dfn}
If $\tau\in\Phi_A^{L,d-1}$ is a facet then the number $\nu_A^{L,\tau}$\index{nLtA@$\nu_A^{L,\tau}$} of indices $i$ in the chain \eqref{56}, and hence in the partition \eqref{92}, with $\tau_i=\tau$ is the {\em multiplicity} of $S^L_A$ along $\bar O_A^\tau$. 
Note that $\nu^{L,\tau}_A$ is the length of the localization of $S^L_A$ at $I^\tau_A$ and hence independent of the particular composition chain.  
\end{dfn}

\begin{prp}\label{58}
For all $\tau\in\Phi^{L,d-1}_A$, $\nu_A^{L,\tau}=[\ZZ^d:\ZZ\tau]$. In
particular, the degree of $I_A^L$ equals $\sum_{\tau\in\Phi_A^{L,d-1}}\nu_A^{L,\tau}=\sum_{\tau\in\Phi_A^{L,d-1}}[\ZZ^d:\ZZ\tau]$.
\end{prp}

\begin{proof}
Fix $\tau\in\Phi_A^{L,d-1}$.
For disjoint $\boldu+\NN\tau$ and $\boldv+\NN\tau$, $\boldu-\boldv\not\in\ZZ\tau$ since $\NN\tau$ contains a shifted copy of its normalization $\QQ_+\tau\cap\ZZ\tau$.
This means that $\nu_A^{L,\tau}\le[\ZZ^d:\ZZ\tau]$ and it remains to show that $\bigsqcup_{\tau_i=\tau}(\boldu_i+\NN\tau_i)$ meets every coset of $\ZZ^d/\ZZ\tau$.
Pick $\boldu\in\NN\tau$ outside $\QQ_+\tau'$ for all
$\tau'\in\Phi_A^L$ with $\tau'\subsetneq\tau$.
Since $\NN A$ contains a shifted copy of its normalization $\QQ_+A\cap\ZZ^d$, for $k\gg 0$
\[
k\boldu+(\QQ_+\tau\cap\ZZ^d)\subseteq
k\boldu+(\QQ_+A\cap\ZZ^d)\subseteq\NN
A=\deg(S_A^L)=\bigsqcup_{i=1}^l\deg(S_{\tau_i}(-\boldu_i)).
\]
By Lemma~\ref{57} below, for $k\gg 0$ 
\[
k\boldu+(\QQ_+\tau\cap\ZZ^d)\subseteq\bigsqcup_{\tau_i=\tau}(\boldu_i+\NN\tau_i)
\]
and the left hand side meets every coset of $\ZZ^d/\ZZ\tau$. 
This yields $\nu_A^{L,\tau}\geq [\ZZ^d:\ZZ\tau]$.
\end{proof}

\begin{lem}\label{57}
Let $\tau\in\Phi_A^{L,d-1}$ and suppose
$\boldu\in\NN\tau\smallsetminus\QQ_+\tau'$ for all $\tau'\in\Phi_A^L$
with $\tau'\subsetneq\tau$. 
Fix $\boldu'\in\ZZ^d$ and $\tau'\in\Phi^L_A$ such that $\tau'\neq\tau$.
Then for $k\gg 0$
\[
(k\boldu+\QQ_+\tau)\cap (\boldu'+\QQ_+\tau')=\emptyset.
\]
\end{lem}

\begin{proof}
Suppose that $k\boldu+\boldv=\boldu'+\boldv'$ with
$\boldv\in\QQ_+\tau$ and $\boldv'\in\QQ_+\tau'$.  Pick a linear form
$0\neq h_{\tau,\tau'}\in\Hom_\QQ(\QQ^d,\QQ)$ with
$h_{\tau,\tau'}(\tau)\geq 0$, $h_{\tau,\tau'}(\tau')\leq 0$.  Then, by
hypothesis, $h_{\tau,\tau'}(\boldu)>0$.  But
$kh_{\tau,\tau'}(\boldu)+h_{\tau,\tau'}(\boldv)=h(\boldu')+h_{\tau,\tau'}(\boldv')$
implies that $kh_{\tau,\tau'}(\boldu)\le h(\boldu')$ which is
impossible for $k\gg0$.
\end{proof}

\subsection{Newton filtration and Cohen--Macaulayness}

Until the end of this section, we fix a weight vector $L$ with
$L_{\del_i}>0$ for all $i$.  Let $V=\bigoplus_{\bolda\in\NN A}
V_\bolda$ be a $\ZZ^d$-graded vector space.  Then the {\em Newton
filtration} $N=N^L_A$\index{NLA@$N^L_A$} on $V$ with respect to
$\Delta^L_A$ is defined by
\[
N_iV=\bigoplus_{\bolda\in i\cdot\Delta^L_A} V_\bolda.
\]
Note that $\gr^N(V)=V$ as $\ZZ^d$-graded vector spaces.
For $\tau\in\Phi_A^L$ or $\tau=A$, we denote $V(\tau)=\bigoplus_{\bolda\in\QQ_+\tau}V_\bolda$.
Generalizing \cite[Prop.~2.6]{Kou76}, there is a complex
\begin{equation}\label{54}
\xymat{
0\ar[r]&V_{d-1}\ar[r]&V_{d-2}\ar[r]&\cdots\ar[r]&V_0\ar[r]&0
}
\end{equation}
where the $\ZZ^d$-graded vector space $V_i$ is the direct sum of all $V(\tau)$ for which $\tau\in\Phi^{L,i}_A$ is not contained in the boundary of $\Gamma^L_A$ (see Definition~\ref{35}). 
The cohomology of this complex is concentrated in homological degree $d-1$ and equal to $V$. 

Now assume that $V$ is a $\ZZ^d$-graded $\CC$-algebra. 
Then $\gr^N(V)(\tau)=V(\tau)$ as $\CC$-algebras for $\tau\in\Phi_A^L$.
For $\tau,\tau'\in\Phi_A^L$ (or $\tau=A$) with $\tau\supseteq\tau'$, the maps
\[
\xymat{\gamma_{\tau,\tau'}\colon\gr^N(V)(\tau)\ar@{->>}[r]^-{}&\gr^N(V)(\tau')}
\]  
in the complex~\eqref{54} for $\gr^N(V)$ are natural projections of $\CC$-algebras.
In particular, \eqref{54} is a complex of $\gr^N(V)$-modules.

\begin{lem}\label{88}
For all $\tau\in\Phi_A^L$ and all $k>0$, $\CC[\NN\cdot k\cdot(\tau\cap\Phi_A^{L,0})]\subseteq S_A^L(\tau)$ is a module-finite ring extension for all $\tau\in\Phi_A^L$.
\end{lem}

\begin{proof}
One can construct a $\ZZ^d$-graded composition chain of $S_A^L$ as in \eqref{56} such that the resulting partition~\eqref{92} refines the partition $\QQ_+A=\QQ_+\tau\sqcup(\QQ_+A\smallsetminus\QQ_+\tau)$, and hence $\deg(S_{\tau_i}(-\boldu_i))$ meets $\QQ_+\tau$ only if $\tau_i\subseteq\tau$. 
Namely, given any chain as in \eqref{56} we refine it at any index $i\in\{1,\dots,l\}$ for which $\deg(S_{\tau_i}(-\boldu_i))=\boldu_i+\NN\tau_i$ meets both $\QQ_+\tau$ and $\QQ_+A\smallsetminus\QQ_+\tau$.
Let $\ell\in\Hom_\QQ(\QQ^d,\QQ)$ be a separating linear form, $\ell(\tau)\ge0$ and $\ell(\tau_i)\le0$.
If $\tau_i\subseteq\tau$ (and hence $\boldu_i\notin\QQ_+\tau$), we chose $\boldu'_i\in(\boldu_i+\NN\tau_i)\cap\QQ_+\tau$. 
The submodule $S_{\tau_i}(-\boldu'_i)$ of $S_{\tau_i}(-\boldu_i)$ has its $\ZZ^d$-degrees entirely in $\QQ_+\tau$, and the cokernel of the inclusion has smaller dimension. 
If conversely $\tau_i\not\subseteq\tau$, choose
$\ell\in\Hom_\QQ(\QQ^d,\QQ)$ above such that $\ell(\tau_i)\neq 0$ and pick $\boldu'_i\in \boldu_i+\NN\tau_i$ with $\ell(\boldu'_i)<0$. 
Then $S_{\tau_i}(-\boldu'_i)\subseteq S_{\tau_i}(-\boldu_i)$ and has degrees completely outside of $\QQ_+\tau$. 

Iterating this procedure we arrive at a composition chain as claimed. 
Then, however, the statement of the lemma is obvious:
for each $i\in\{1,\dots,l\}$ with $\deg(S_{\tau_i}(-\boldu_i))\cap\QQ_+\tau\ne\emptyset$, $\tau_i\subseteq\tau$, and $(S_{\tau_i}(-\boldu_i))(\tau)\cong S_{\tau_i}$ is even module-finite over the subring $\CC[\NN\cdot k\cdot(\tau_i\cap\Phi_A^{L,0})]$ of $\CC[\NN\cdot k\cdot(\tau\cap\Phi_A^{L,0})]$.
\end{proof}

\begin{thm}\label{36}
If $S_A^L(\tau)$ is Cohen--Macaulay for all $\tau\in\Phi_A^L$ not contained in the topological boundary of $\Gamma_A^L$ then $S_A^L$ is Cohen--Macaulay.
\end{thm}

\begin{proof}
There is an integer $k>0$ such that $k\cdot\Gamma_A^L$ has vertex set
$k\cdot\Gamma_A^{L,0}\subseteq\NN A$.  Within the space of sequences
$f_1,\ldots,f_d$ in $S^L_A$ for which each $f_i$ is a sum of terms
whose $\ZZ^d$-degrees are in $k\cdot \Gamma^{L,0}_A$, choose one
sequence that is generic.

Each $f_i$ is homogeneous with respect to $N=N^L_A$ and can hence be
identified with its $N$-symbol in $\gr^N(S_A^L)$.  By Lemma~\ref{88},
the vector space spanned by
$\gamma_{A,\tau}(f_1),\dots,\gamma_{A,\tau}(f_d)\in\CC[\NN\tau]$
contains a system of parameters on $S_A^L(\tau)$ for all
$\tau\in\Phi_A^L$.  This system of parameters is a regular sequence by
the Cohen--Macaulay hypothesis.  By the (spectral sequence) argument
in \cite[\S2.12]{Kou76}, the Koszul complex induced by $f_1,\dots,f_d$
on the complex~\eqref{54} with $V=S^L_A$ is a resolution.  It follows
that the Koszul complex induced by $f_1,\dots,f_d$ on $\gr^N(S_A^L)$,
and hence on $S_A^L$, is a resolution as well.  Therefore,
$f_1,\dots,f_d$ is a regular sequence on $S_A^L$ and the claim
follows.
\end{proof}

\section{Characteristic variety of the hypergeometric system}\label{charvar}

\subsection{Characteristic varieties and slopes}\label{18}

The {\em Weyl algebra} $D=\CC[x]\ideal{\del}$\index{D@$D$} in $n$
variables $x=x_1,\dots,x_n$ is the ring of $\CC$-linear differential
operators on $X=\CC^n$\index{X@$X$} and contains $R$ as a commutative
subring.  The $\ZZ^d$-grading on $R$ extends to $D$ by setting
$-\deg(x_i)=\bolda_i=\deg(\del_i)$\index{deg@$\deg$}.

With $L_x=(L_{x_1},\ldots,L_{x_n})$ and $L_\del=(L_{\del_1},\ldots,L_{\del_n})$, $L=(L_x,L_\del)\in\QQ^{2d}$ is a \emph{weight vector on $D$} if $L_x+L_\del\geq 0$. Fix any such weight vector $L$; it defines an increasing filtration $L$ on $D$ by $[x^\boldu\del^\boldv\in L_kD]\Leftrightarrow[L\cdot(\boldu,\boldv)\le k]$. 
Since the Gr\"obner fan of any $D$-ideal is defined over $\QQ$ (see \cite{ACG00}), the study of real weights can be reduced to the present rational case. 
If $P\in L_kD\smallsetminus L_{<k}D$ then $k=:\deg^L(P)$ is the \emph{$L$-degree} of $P$. 
The multiplicative, but not additive, {\em $L$-symbol map}
\[
\xymat{\sigma^L:D\ar[r]&\gr^L(D)=:W}\index{sigmaL@$\sigma^L$}\index{W@$W$}
\]
is defined by $\sigma^L(P)=P\modulo L_{<k}R$ if $\deg^L(P)=k$. 
An element of the form $P=\sum_{L\cdot(\boldu,\boldv)=k}P_{\boldu,\boldv}x^\boldu\del^\boldv\in D$ with $P_{\boldu,\boldv}\in\CC$ is \emph{$L$-homogeneous}. 
By abuse of notation, we identify $L$-homogeneous elements in $D$ with their image under $\sigma^L$ in $W$. 
We restrict ourselves to the case $L_x+L_\del>0$ in which case $W\cong \CC[x,\del]$ can be considered as the ring of polynomial functions on the cotangent space $T^*X$

The definition of the characteristic variety of a $D$-module is based on the concept of good filtrations discussed in more detail in Section~\ref{cycle} and \cite[Ch.~II, \S\S1.1--1.3]{Sch85}.
Let $F$ be a filtration on a ring $T$.
Then $G$ is called a \emph{good $F$-filtration} on a $T$-module $N$ (see \cite[Ch.~II, Def.~1.1.1]{Sch85}), if there are generators $n_1,\dots,n_m$ of $N$ and $\boldu\in\ZZ^m$ such that for all $k$ one has
\[
G_kN=\sum_{i=1}^m F_{k+u_i}T\cdot n_i.
\]
Note that good $F$-filtrations on $N$ exist if and only if $N$ is $T$-finite. 
From the definition follows the fact that all good $F$-filtrations on $N$ are equivalent in the sense that for all $k,l\in\NN$ there are $k_l,l_k\in\NN$ with $G_kN\subseteq G'_{l_k}N$ and $G'_l\subseteq G_{k_l}$.

\begin{dfn}\label{24}
The {\em $L$-characteristic variety} $\ch^L(M)$\index{chLM@$\ch^L(M)$} of a finite $D$-module $M$ on $X$ is the support of $\gr^L(M)$ in $T^*X$ for some good $L$-filtration on $M$.
A finite $D$-module $M$ is \emph{$L$-holonomic} if $\dim\ch^L(M)=n$.
\end{dfn}

The independence of $\ch^L(M)$ of the choice of the good $L$-filtration on $M$ follows from \cite[Ch.~II, Prop.~1.3.1.a]{Sch85}.
The following algebraic statement is a special case of a result by G.G.~Smith \cite[Thm.~1.1]{Smi01}.

\begin{thm}\label{64}
The dimension of any component of $\ch^L(M)$ is at least $n$.\qed
\end{thm}

Important special cases of filtrations on $D$ are the the {\em order filtration} $F=(F_x,F_\del)=(\boldzero,\boldone)$\index{F@$F$}\index{Fx@$F_x$}\index{Fd@$F_\del$} and the {\em $V$-filtration} $V=(V_x,V_\del)$\index{V@$V$}\index{Vx@$V_x$}\index{Vd@$V_\del$} along the coordinate variety
\begin{equation}\label{94}
Y:=\Var(x_\frakV)\subseteq X\index{Y@$Y$},\quad\frakV\subseteq\{1,\dots,n\}\index{V@$\frakV$},
\end{equation}
defined by $-V_{x_i}=1=V_{\del_i}$ for $i\in\frakV$ and $V_{x_i}=0=V_{\del_i}$ for $i\notin\frakV$.
The notion of slopes along $Y$ (see Definition~\ref{91}) involves the family of intermediate filtrations $L$ between $F$ and $V$ defined by the linear combination of weight vectors
\begin{equation}\label{73}
L=pF+qV,\quad p/q\in\QQ_{>0}\cup\{\infty\}\index{L@$L$}.
\end{equation}
Note that $L_x+L_\del>0$ since $p>0$.
If $p'/q'=p/q$ then the filtrations $L$ and $L'$ are identical, up to a dilation in the index: $L_{p'k}=L'_{pk}$.
By abuse of notation, we shall frequently identify the filtration $L$ with the number $p/q$. 
For $Y\subseteq X$ closed and reduced, let $T^*_YX$\index{TYX@$T^*_YX$} be (the closure of) the conormal bundle of (the smooth points of) $Y$ in $X$.
With notation as in \eqref{94} and \eqref{73}, the ring $W=\gr^L(D)$ can be considered as the ring of polynomial functions on the cotangent space $T^*T^*_YX$ of $T^*_YX$.
For $i\in\frakV$, $-x_i$ can be interpreted as the partial derivative with respect to the variable $\del_i$.
This sets up an explicit isomorphism between $T^*T^*_YX$ and $T^*X$ .

\begin{rmk}\label{63}
For $L=pF+qV$ as in \eqref{73}, $\ch^L(M)$ is the global algebraic version of Y.~Laurent's {\em microcharacteristic variety $\Sigma_\Lambda^r(\calM)$\index{SLr@$\Sigma_\Lambda^r$} of type $r=p/q+1$} in $T^*\Lambda$ with $\Lambda=T^*_YX$ (see \cite[\S3.2]{Lau87}).
Our algebraic $L$-filtration corresponds to the filtration on the sheaf of analytic differential operators on $X$ along $Y$ induced by Y.~Laurent's microlocal filtration $F_{\Lambda,r}$ along $\Lambda$ (see \cite[Def.~3.2.1]{Lau87}).
By a flatness argument, A.~Assi et al.~\cite[Lem.~1.1.2]{ACG96} show that the analytification functor commutes with the grading by these two corresponding filtrations. 
Thus, the components of our $\ch^L(M)$ which meet the preimage of $Y$ in $T^*X$ correspond to the components of $\Sigma_\Lambda^r(\calM)$.
\end{rmk}

The preceding remark motivated the following algebraic version of Y.~Laurent's {\em critical indices} (see \cite[\S3.4]{Lau87}) which are also called slopes.

\begin{dfn}\label{91}
For $L=pF+qV$ as in \eqref{73}, we mean by \emph{$f(L')$ jumps at $L'=L$} that the
set-valued function
\[
\xymat{
\QQ\ni p'/q'\ar@{|->}[r]&L'=p'F+q'V\mapsto f(L')
}
\]
is not locally constant at $p/q$.
A {\em slope} of a finite $D$-module $M$ \emph{at $y\in Y$ along $Y$} is a value $L=p/q$ such that the set of components of $\ch^{L'}(M)$ which meet $T_y^*X$ jumps at $L'=L$.
\end{dfn}

It follows from the existence of the Gr\"obner fan (see \cite{ACG00}) that a fixed $D$-module has only finitely many and only rational slopes along all coordinate varieties. 
Y.~Laurent proved this finiteness and rationality along general varieties in the microlocal setting in \cite[Thm.~3.4.1]{Lau87}.
He also showed that $\Sigma_\Lambda^r(\calM)$ is involutive (see \cite[Prop.~3.5.2]{Lau87}) and Lagrangian for holonomic $\calM$ (see \cite[Cor.~4.1.2.(ii)]{Lau87}).
In view of Remark~\ref{63}, this implies for holonomic $M$ that $\ch^L(M)$ is purely $n$-dimensional for $L=pF+qV$ as in \eqref{73}.

\subsection{Hypergeometric system and candidate components}

We now define our main object of interest, the hypergeometric $D$-module $M_A(\beta)$, introduced in \cite{GGZ87,GKZ89}.

\begin{dfn}\label{81}
The {\em Euler vector fields} $E=(E_1,\dots,E_d)$ of $A$ are defined by $E_i:=\sum_j a_{i,j}x_j\del_j$\index{E@$E$} for $i=1,\dots,d$.
The {\em $A$-hypergeometric} (or \emph{GKZ}) \emph{system} defined by
$A$ and a complex parameter vector $\beta\in\CC^d$\index{b@$\beta$} is
the $D$-module
\[
M_A(\beta):=D/D\ideal{I_A,E-\beta}\index{MAb@$M_A^\beta$}
\]
on $X$ defined by the toric ideal $I_A$ and the {\em Euler operators} $E-\beta$.
\end{dfn}

The $A$-hypergeometric system is always holonomic (see \cite{Ado94}).
Our goal is to understand its $L$-characteristic varieties and slopes under the assumption that
\begin{equation}\label{96}
L_x+L_\del=(c,\ldots,c)=:\boldc\text{ for some rational }c>0.
\end{equation}
This guarantees that $W$ is a polynomial ring and $E$ is $L$-homogeneous of positive $L$-degree. 
Note that $M_A(\beta)$ is $L$-homogeneous if $A$ is $L$-homogeneous as defined in Definition~\ref{35}. 
The following statement is a consequence of Buchberger's algorithm.

\begin{lem}\label{37}
One has the identity $\gr^L(DI_A)=WI_A^L$.\qed
\end{lem}

The vector fields $t_1\del_{t_1},\dots,t_d\del_{t_d}$ span the tangent space at any point of the torus $\TT=(\CC^*)^d$. 
Hence, for any $\xi\in T^*_0X$, the tangent space of the orbit $\Orb(\xi)$ is spanned by the pushforwards $E_i^T:=\sum_ja_{i,j}\del_j(-x_j)$ of $t_i\del_{t_i}$ under the map $\TT\to T^*_0X$, $t\mapsto t\cdot\xi$ from \eqref{16}. 
Since $E_i$ is $L$-homogeneous of positive degree with $\sigma^L(E_i)=-\sigma^L(E^T_i)$, the equations
$\sigma^L(E-\beta)=0$ impose the conormal condition to the orbit $\Orb(\xi)\subseteq T^*_0X$ in $T^*T^*_0X=T^*X$. 
In what follows we abuse notation by writing $\sigma^L(E)=E$.

\begin{dfn}\label{74}
For a subset $\tau$ of $\{1,\ldots,n\}$, we denote by
$C_A^\tau\subseteq T^*X$\index{Ct@$C_A^\tau$} the conormal space to the
orbit $O_A^\tau\subseteq T^*_0X$ from Definition~\ref{17}.  We denote
by $P_C\subseteq W$\index{PC@$P_C$} the defining ideal of any
irreducible variety $C\subseteq T^*X$ and abbreviate $P_{C_A^\tau}$ by
$P_\tau$\index{Pt@$P_\tau$}.
\end{dfn}

\begin{prp}\label{55}
The $L$-characteristic variety of $M_A(\beta)$ is 
\[
\ch^L(M_A(\beta))=\bigcup_{\tau\in\phi^L_A(\beta)}\bar
C_A^\tau=\bigsqcup_{\tau\in\phi^L_A(\beta)}C_A^\tau
\]
for some subset $\phi_A^L(\beta)\subseteq\Phi^L_A$\index{pALb@$\phi_A^L(\beta)$}.
In particular, the hypergeometric system $M_A(\beta)$ is $L$-holonomic for any $L$.
\end{prp}

\begin{proof}
By definition, $\ch^L(M_A(\beta))\subseteq\Var(W\ideal{I_A^L,E})$. 
By Theorem~\ref{13} and the preceding arguments,
\begin{equation}\label{67}
\Var(W\ideal{I_A^L,E})
=\bigsqcup_{\tau\in\Phi_A^L}C_A^\tau
=\bigcup_{\tau\in\Phi_A^L}\bar C_A^\tau.
\end{equation}
Then Theorem~\ref{64} assures that $\ch^L(M_A(\beta))$ is purely $n$-dimensional and hence a union of closures of conormals $C_A^\tau$ for certain $\tau\in\Phi_A^L$.
\end{proof}

\subsection{Existence of facet components}

In this subsection, we use an elementary localization argument to give
an index formula for the multiplicity of components in the characteristic cycle $\cc^L(M_A(\beta))$ (see Definition~\ref{27}) corresponding to facets in the $(A,L)$-umbrella.
It shows in particular that these multiplicities are independent of $\beta$ and positive, and that all facet components occur in $\ch^L(M_A(\beta))$.
We deduce from this a converse to R.~Hotta's Theorem: regular $A$-hypergeometric systems are homogeneous in the usual sense.

\begin{dfn}\label{32}
We denote by $\mu_{A,0}^{L,C}(\beta)$\index{mLCAb@$\mu_{A,0}^{L,C}(\beta)$} the multiplicity of $\gr^L(M_A(\beta))$ along $C$. 
(The reason for the appearance of the subscript ``$0$'' will become apparent in Section \ref{cycle}.)
This is the length of the $W_{P_C}$-module $\gr^L(M_A(\beta))\otimes_WW_{P_C}$.
We write $\mu_{A,0}^{L,\tau}(\beta)$\index{mLtAb@$\mu_{A,0}^{L,\tau}(\beta)$} if $C=\bar C_A^\tau$.
\end{dfn}

\begin{thm}\label{10}
For all $\beta\in\CC^d$ and all $\tau\in\Phi_A^{L,d-1}$, 
\[
\mu_{A,0}^{L,\tau}(\beta)=\nu^{L,\tau}_A=[\ZZ^d:\ZZ\tau].
\]
In particular, $\Phi^{L,d-1}_A\subseteq\phi_A^L(\beta)$ for all $\beta\in\CC^d$.
\end{thm}

\begin{proof}
Let $\tau\in\Phi^{L,d-1}$ and relabel columns such that
$\bolda_1^L,\dots,\bolda_d^L\in\tau$ and such that
$\bolda_1,\dots,\bolda_d$ are linearly independent.  We have to show
that $C_A^\tau\subseteq\ch^L(M_A(\beta))$.  First, we verify that $E$
is a regular sequence on
$W[\partial_\tau^{-1}]/W[\partial_\tau^{-1}]I_A^L$.  After Gauss
reduction on $E$ and multiplying $E_i$ by $\partial_i^{-1}$,
$E_i\equiv x_i$ modulo terms independent of $x_j$ for all $j\le d$.
After a change of the coordinates $x_1,\dots,x_d$ in
$W[\partial_\tau^{-1}]$, leaving $I_A^L$ invariant, $E_i=x_i$.  As
$I_A^L$ does not involve the variables $x_1,\dots,x_d$, $E$ is a
regular sequence modulo $W[\partial_\tau^{-1}]I_A^L$ as claimed.

Since $W[\partial_\tau^{-1}]I_A^L=\gr^L(D[\partial_\tau^{-1}]I_A)$ by Lemma~\ref{11} below, $\gr^L(D[\partial_\tau^{-1}]\ideal{I_A,E-\beta})=\gr^L(D[\partial_\tau^{-1}]I_A)+\gr^L(D[\partial_\tau^{-1}]\ideal{E})$ by the argument in \cite[Thm.~4.3.5]{SST00}. 
Again by Lemma~\ref{11},
\[
\gr^L(M_A(\beta))[\partial_\tau^{-1}]
=W[\partial_\tau^{-1}]/W[\partial_\tau^{-1}]\ideal{I_A^L,E}.
\]
Since $\partial_j\notin P_\tau$ for all $j\in\tau$, this yields
\[
\gr^L(M_A(\beta))_{P_\tau}
=W_{P_\tau}/W_{P_\tau}\ideal{I_A^L,E}
\cong(R_{I_A^\tau}/I_A^L)(x_{d+1},\dots,x_n)
\]
with the ring isomorphism defined by the above coordinate change, and hence
\[
\mu_{A,0}^{L,\tau}(\beta)
=\ell(\gr^L(M_A(\beta))_{P_\tau})
=\ell(R_{I_A^\tau}/I_A^L)
=\nu^{L,\tau}_A.
\]
Then Proposition~\ref{58} finishes the proof.
\end{proof}

\begin{lem}\label{11}
For $I\subseteq D$ and $\tau\subseteq\{1,\dots,n\}$, $\gr^L(D[\partial_\tau^{-1}]I)=\gr^L(D[\partial_\tau^{-1}])\gr^L(I)$.
In particular, $\gr^L((D/I)[\partial_\tau^{-1}])=(\gr^L(D/I))[\partial_\tau^{-1}]$.
\end{lem}

\begin{proof}
The inclusion $\gr^L(D[\del_\tau^{-1}]I)\supseteq\gr^L(D[\del_\tau^{-1}])\gr^L(I)$ holds trivially.
For $Q\in D[\del_\tau^{-1}]$, $\del_\tau^\boldu Q\in D$ for some
$\boldu\in\NN^{|\tau|}$. If $\{P_i\}$ is a finite set of generators
for $I$ 
then the opposite inclusion follows from 
\[
\sigma^L\bigl(\sum_iQ_iP_i\bigr)=\del_\tau^{-\boldu}\sigma^L\bigl(\sum_i\del_\tau^{\boldu}Q_iP_i\bigr)\in\gr^L(D[\del_\tau^{-1}])\gr^L(I)
\]
with a common $\boldu$ for the finitely many $Q_i\in D[\del^{-1}_\tau]$.
In particular,
\begin{gather*}\pushQED{\qed}
\gr^L((D/I)[\partial_\tau^{-1}])
=\gr^L(D[\partial_\tau^{-1}]/D[\partial_\tau^{-1}]I)
=\gr^L(D[\partial_\tau^{-1}])/\gr^L(D[\partial_\tau^{-1}])\gr^L(I)\\
=W[\partial_\tau^{-1}]/W[\partial_\tau^{-1}]\gr^L(I)
=(W/\gr^L(I))[\partial_\tau^{-1}]
=(\gr^L(D/I))[\partial_\tau^{-1}].\qedhere
\end{gather*}
\end{proof}

\begin{exa}\label{97}
Let $A$ be as in Figure~\ref{70}; we consider $L=F+tV$ where $V$ is
the $V$-filtration along $\Var(x_4)$, induced by the weight vector
$(0,0,0,-1,0,0,0,1)$. We consider specifically $t\in\{0,1,4\}$. As
these weights are generic, for all $\tau\subseteq\{1,\ldots,4\}$ the
conormal closures $\bar C^\tau_A$ are coordinate subspaces
$\Var(x_\tau, \del_{\ol\tau})$.  The facets in $\Phi^L_A$ for these
values of $t$ are shown in Figure~\ref{70}.

For $t=0$ one finds $\mu^{L,\{1,4\}}_{A,0}(\beta)=12$ and
$\mu^{L,\{2,4\}}_{A,0}(\beta)=1$ with sum
$13=\vol_{\ZZ^d}(\Delta^F_A)$ (for relevant notation and more information, see Definition~\ref{98} and the continuation in Example~\ref{100}). 
For $t=1$ one finds three facets, with $\mu^{L,\{1,3\}}_{A,0}(\beta)=3$,
$\mu^{L,\{3,4\}}_{A,0}(\beta)=7$ and $\mu^{L,\{2,4\}}_{A,0}(\beta)=1$. 
On the other hand, $t=4$ yields two facets only, with
$\mu^{L,\{1,3\}}_{A,0}(\beta)=3$ and $\mu^{L,\{2,3\}}_{A,0}(\beta)=2$.
The respective sums, $3+7+1$ and $3+2$, are the degrees of the $L$-graded toric ideals (see Proposition \ref{58}).
\end{exa}

\begin{lem}\label{12}
The map $L'\mapsto\Phi^{L',d-1}_A$ jumps at $L'=L$ if and only if the map $L'\mapsto\Phi^{L'}_A$ jumps at $L'=L$.
\end{lem}

\begin{proof}
The ``only if" part in the statement follows trivially from $\Phi^{L',d-1}_A\subseteq\Phi^{L'}_A$.
Assume that $\Phi^{L',d-1}_A$ is locally constant at $L$ and let $\tau\in\Phi^L_A$.
Then there is a $\tau'\in\Phi^{L,d-1}_A$ such that $\tau$ is a face of $\tau'$.
By assumption, $\tau'\in\Phi^{L',d-1}_A$ for $L'$ close to $L$.
Since $\tau'$ is not contained in a hyperplane through origin and since the $\bolda_j^{L'}$ depend on $L'$ only by scaling, $\tau$ remains a face of $\tau'$ for $L'$ close to $L$.
In particular, $\tau\in\Phi^{L'}_A$ for $L'$ close to $L$ and hence $\Phi^{L'}_A$ is locally constant at $L$.
\end{proof}

Theorem~\ref{10} and Lemma~\ref{12} yield the existence of slopes of $M_A(\beta)$ corresponding to jumps of the $(A,L)$-umbrella with $L$ as in \eqref{73}. 
The following result on the candidate components $\bar C_A^\tau$ of $\ch^L(M_A(\beta))$ allows us to show that all these slopes occur at the origin.

\begin{lem}\label{78}
For all $\tau\in\Phi_A^L$, the candidate component $\bar C_A^\tau$
meets $T_0^*X$.
\end{lem}

\begin{proof}
Since the ideal $W\ideal{I_A^L,E}$ is homogeneous in the $x$-variables, this holds also for its associated prime ideals $P_\tau$ where $\tau\in\Phi_A^L$.
Thus $\bar C_A^\tau=\Var(P_\tau)$ meets $T_0^*X$ as claimed.
\end{proof}

\begin{cor}\label{66}
For $L=pF+qV$ as in \eqref{73}, if $\Phi^{L}_A$ jumps at $L=L'$ then $L'=p'/q'$ is a slope of $M_A(\beta)$ along $Y$ at $0\in Y$.\qed
\end{cor}

In particular, we obtain the following converse to a Theorem by R.~Hotta \cite[Ch.~II, \S6.2, Thm.]{Hot98}.

\begin{cor}\label{38}
Regular $A$-hypergeometric systems are homogeneous with respect to the order filtration $F$.
\end{cor}

\begin{proof}
The polytope $\conv(A)$ is an intersection of closed half-spaces and $0\notin\conv(A)$ by positivity of $\NN A$.
Thus there is a facet $\tau$ of $\conv(A)$ and a linear form $\ell\in\Hom_\QQ(\QQ^d,\QQ)$ such that $\ell(\bolda_i)\ge1$ with equality equivalent to $i\in\tau$.
If $\tau\in\Phi_A^{F,d-1}$ then $\tau$ is a facet of $\Delta^F_A=\conv(\{0,\bolda_1,\dots,\bolda_n\})$.
But then $0\in\Delta^F_A$ and $\ell(0)=0<1$ implies that $\ell(\bolda_i)\le1$ and hence $\ell(\bolda_i)=1$ for $i=1,\dots,n$.
So in this case $A$ is $F$-homogeneous.
Assume now that $\tau\not\in\Phi_A^{F,d-1}$ and let $L=F+\frac{q}{p}V$
where $\frakV=\ol\tau$ in \eqref{94}.
For $i\in\tau$, $\bolda_i^L=\bolda_i$ and hence $\ell(\bolda_i^L)=1$.
But for $i\notin\tau$ and $p/q\to0$, $\bolda_i^L\to0$ and hence eventually $\ell(\bolda_i^L)<1$.
Thus, $\tau\in\Phi_A^{L,d-1}$ for $p/q\to0$ while $\tau\not\in\Phi_A^{F,d-1}$ and hence $L\mapsto\Phi_A^L$ can not be constant.
Then, by Corollary~\ref{66}, $M_A(\beta)$ has slopes. By
\cite[Thm.~2.4.2]{LM99},  $M_A(\beta)$  is hence irregular.
\end{proof}

We will see later in Corollary~\ref{39} that Corollary~\ref{66} actually gives a complete list of all slopes along coordinate varieties at the origin.

\begin{rmk}
For generic $\beta$, the equivalence of the regularity of $M_A(\beta)$
with the homogeneity of $A$ was already obtained in \cite[Thm.~2.4.11]{SST00}.
\end{rmk}

\section{Characteristic cycle of the Euler--Koszul complex}\label{cycle}

In this section we assume that $L_x+L_\del=\boldc$ is a constant positive vector as in \eqref{96}. 
By way of a dilation we further may, and do, assume that the index set of the $L$-filtration is $\ZZ$.

By Proposition~\ref{55}, all components of the $L$-characteristic
variety of $M_A(\beta)$ are of the form $\bar C_A^\tau$ where
$\tau\in\Phi_A^L$.  In this section we prepare the way for
Corollary~\ref{23} which ascertains the presence of every such
candidate component.  The approach is to consider $M_A(\beta)$ as the
$0$-th homology of a Koszul type complex by operators $E-\beta$ on the
$D$-module $D\otimes_R N$ for the $\ZZ^d$-graded $R$-module $N=S_A$.
For modules $N$ having a composition series with quotients of type
$S_\theta$, $\theta\in\Phi^\boldzero_A$, we can apply results in Sections \ref{spectoric} and \ref{charvar} combined with homological algebra. 
For basic results on filtered rings and modules we refer to \cite[Ch.~II, \S\S1.1-1.3]{Sch85}.

\subsection{Good filtrations and toric modules}\label{14}

In order to combine homological methods with good $L$-filtrations on $D$-modules, we need the $L$-filtration on $D$ to be Noetherian, which means the following. 

The \emph{Rees ring} of a ring $T$ with a $\ZZ$-indexed filtration $F$
is the graded ring
\[
\Rees^F(T)=\bigoplus_{i\in\ZZ}F_i(T)t^i\subseteq T[t,t^{-1}]
\]
and $F$ is called \emph{Noetherian} if $\Rees^F(T)$ is Noetherian.

\begin{lem}\label{31}
The filtrations $L$ on $D$ and $L_\del$ on $R$ are Noetherian.
\end{lem}

\begin{proof}
Let us write $(\boldu,\boldv):=(L_x,L_\del)$, then $\boldu+\boldv=\boldc>0$ by hypothesis. 
The elements of $\Rees^L(D)$ are of the form $\sum_{k\in\ZZ}\sum_{\boldu\cdot\boldu'+\boldv\cdot\boldv'\le k}c_{\boldu',\boldv',k}x^{\boldu'}\del^{\boldv'} t^k$. 
Apply the change of variables $x_i\mapsto x_it^{-u_i}$, $\del_i\mapsto\del_it^{u_i}$ and set $m=k+\boldu\cdot(\boldv'-\boldu')$. 
Then the elements are transformed to $\sum_{m\in\ZZ}\sum_{(\boldv+\boldu)\cdot\boldv'\le m}x^{\boldu'}\del^{\boldv'}t^m$.  Thus, $\Rees^L(D)\cong\Rees^{F'}(D)$ where $F':=(0,\boldv+\boldu)=(0,\boldc)$, which is Noetherian by \cite[Ch.~II, Prop.~1.1.8]{Sch85} since $F'_0D=\CC[x]$ is Noetherian.
Similarly, a change of variables $\del_i\mapsto\del_it^{1-v_i}$ shows that $\Rees^{L_\del}(R)\cong\Rees^{F_\del}(R)$ is Noetherian.
\end{proof}

The \emph{Rees module} of a $T$-module $N$ with an $F$-filtration $G$ is the graded $\Rees^F(T)$-module 
\[
\Rees^G(N)=\bigoplus_{i\in\ZZ}G_i(N)t^i\subseteq N[t,t^{-1}].
\]
For Noetherian $F$, $G$ is a good $F$-filtration on $N$ in the sense of Subsection~\ref{18} if and only if $\Rees^G(N)$ is Noetherian over $\Rees^F(T)$, (see \cite[Ch.~II, Prop.~1.1.7]{Sch85}). 
The advantage of this new definition is the following. 
Consider a short exact sequence of $T$-modules
\begin{equation}\label{72}
\xymat{0\ar[r]&N'\ar[r]&N\ar[r]&N''\ar[r]&0},
\end{equation}
a good $F$-filtration $G$ on $N$, and the induced $F$-filtrations $G'$ on $N'$ and $G''$ on $N''$. 
Then the associated Rees sequence 
\[
\xymat{
0\ar[r]&\Rees^{G'}(N')\ar[r]&\Rees^{G}(N)\ar[r]&\Rees^{G''}(N'')\ar[r]&0
}
\]
is exact and hence $G'$ and $G''$ are good according to the new definition.
By definition of $G'$ and $G''$, the maps in the sequence~\eqref{72} are strict and hence 
\[
\xymat{
0\ar[r]&\gr^{G'}(N')\ar[r]&\gr^{G}(N)\ar[r]&\gr^{G''}(N'')\ar[r]&0
}
\]
is exact.

For $D$-modules of the type $M=D\otimes_RN$ where $N$ is an $R$-module, a good $L_\del$-filtration on $N$ induces a good $L$-filtration on $M$ by
\begin{equation}\label{71}
L_k(D\otimes_RN)=\sum_{i+j=k}(L_x)_i\CC[x]\otimes_\CC (L_\del)_jN.
\end{equation}
Here, $(L_x)_i\CC[x]$ denotes the level-$i$ piece of the filtration
$L_x$ on $\CC[x]$.
For many purposes, we can replace a given good $L$-filtration on $M$ by that in \eqref{71} and then Lemma~\ref{37} generalizes as follows.

\begin{lem}\label{29}
For any $R$-module $N$, $\gr^{L}(D\otimes_RN)=W\otimes_R\gr^{L_\del}(N)$ as $W$-modules.\qed
\end{lem}

The following definition (see \cite[Def.~4.5]{MMW05}) describes a class of $R$-modules that arise naturally in the study of $A$-hypergeometric systems. 
Recall that, by Theorem~\ref{13}, the $\ZZ^d$-graded prime ideals containing $I_A$ are of the form
$I_A^\theta$ for $\theta\in\Phi^\boldzero_A$. 

\begin{dfn}
A $\ZZ^d$-graded $R$-module $N$ is {\em toric} if it has a {\em toric filtration}
\[
0=N_0\subsetneq N_1\subsetneq\cdots\subsetneq N_{l-1}\subsetneq N_l=N
\]
which means that for all $i$, $N_i/N_{i-1}=S_{\theta_i}(\boldu_i)$ for some
$\theta_i\in\Phi^\boldzero_A$ and some $\boldu_i\in\ZZ^d$.  The minimal such $l$
is called the {\em toric length} of $N$.  A \emph{toric morphism} is a
$\ZZ^d$-graded $R$-linear map between toric modules.
\end{dfn}

\subsection{Characteristic cycle and Euler--Koszul homology}

Recall that by Theorem~\ref{64} any component of the $L$-characteristic variety $\ch^L(M)$ of a finite $D$-module $M$ is at least $n$-dimensional.

\begin{dfn}\label{27}
The {\em $L$-characteristic cycle} of an $L$-holonomic $D$-module $M$ is the formal sum of $n$-dimensional irreducible varieties in $T^*X$
\[
\cc^L(M):=\sum_C\mu^{L,C}(M)\cdot C\index{CCL@$\cc^L(M)$},\quad\mu^{L,C}(M):=\ell\bigr(\gr^L(M)_{P_C}\bigl)\index{mLC@$\mu^{L,C}(M)$},
\]
for some good $L$-filtration on $M$.

A bounded complex of $D$-modules $K_\bullet$ is called {\em homologically $L$-holonomic} if all its homology modules $H_i(K_\bullet)$ are $L$-holonomic.
The {\em $L$-characteristic cycle} of a homologically $L$-holonomic complex of finite $D$-modules $K_\bullet$ is 
\begin{gather*}
\cc^L(K_\bullet):=\sum_{i\in\ZZ}(-1)^i\cc^L(H_i(K_\bullet))=\sum_C\mu^{L,C}(K_\bullet)\cdot
C\index{CCL@$\cc^L(K_\bullet)$},
\end{gather*}
where
\begin{gather*}
\mu^{L,C}(K_\bullet):=\sum_{i\in\ZZ}(-1)^i\mu^{L,C}_i(K_\bullet)\index{mLC@$\mu^{L,C}(K_\bullet)$},\quad
\mu^{L,C}_i(K_\bullet):=\ell\bigl(\gr^L(H_i(K_\bullet))_{P_C}\bigr)\index{mLCi@$\mu^{L,C}_i(K_\bullet)$}.
\end{gather*}
For an $L$-holonomic $D$-module $M$ considered as a complex with trivial differential concentrated in homological degree zero, the two definitions of $\cc^L(M)$ coincide.
\end{dfn}

A $D$-module is holonomic in the usual sense if it is $L$-holonomic for $L=F$.
The independence of $\cc^L(M)$ of the particular choice of the good $L$-filtration follows from \cite[Ch.~II, Prop.~1.3.1.a]{Sch85}. 
Essentially by definition, the $L$-characteristic variety $\ch^L(M)$ of an $L$-holonomic $M$ in Definition~\ref{24} is the union
\[
\ch^L(M)=\bigcup_{\mu^{L,C}(M)>0}C.
\]
For $L=pF+qV$ as in \eqref{73}, $\cc^L(M)$ is the global algebraic version of Y.~Laurent's {\em microcharacteristic cycle $\widetilde\Sigma_\Lambda^r(\calM)$\index{SLr@$\widetilde\Sigma_\Lambda^r$} of type $r=p/q+1$} in $T^*\Lambda$ with $\Lambda=T^*_YX$ whose support is $\Sigma_\Lambda^r(\calM)$ (see Remark~\ref{63}).
Lemma~\ref{31} and \cite[Ch.~II, Prop.~1.3.1.b]{Sch85} yield the following statement.

\begin{lem}\label{34}
The $L$-characteristic cycle $\cc^L$ is additive.\qed
\end{lem}

Our main technical tool for the study of $\cc(M_A(\beta))$ is the Euler--Koszul functor from \cite{MMW05}.
For a $\ZZ^d$-graded left $D$-module $M$ and a $\ZZ^d$-homogeneous $y\in M$,
\begin{equation}\label{86}
(E_i-\beta_i)\circ y:=(E_i-\beta_i+\deg_i(y))y
\end{equation}
defines $d$ commuting $D$-linear endomorphisms $E-\beta$ of $M$.

\begin{dfn}\label{19}
Let $N$ be a $\ZZ^d$-graded $R$-module and $\beta\in\CC^d$.
The \emph{Euler--Koszul complex} $K_{A,\bullet}(N;\beta)$\index{KANb@$K_{A,\bullet}(N;\beta)$} is the Koszul complex of the endomorphisms $E-\beta$ on the left $D$-module $D\otimes_RN$.
Its $i$-th homology $H_{A,i}(N;\beta):=H_i(K_{A,\bullet}(N;\beta))$\index{HAiNb@$H_{A,i}(N;\beta)$} is the $i$-th \emph{Euler--Koszul homology} of $N$.
\end{dfn}

To a good $L_\del$-filtration on a $\ZZ^d$-graded $R$-module $N$ we
associate a good $L$-filtration on $K_{A,\bullet}(N;\beta)$ by setting
for $i_1<\dots<i_k$ the $L$-degree of the $(i_1,\dots,i_k)$-th unit
vector in $ K_{A,k}(R;\beta) = D^{\binom{n}{k}} $ equal to
$\sum_{t=1}^k\deg^L(E_{i_t})$, and using equation~\eqref{71}. 
Then Lemma~\ref{29} implies that
\begin{equation}\label{30}
\gr^L(K_{A,\bullet}(N;\beta))=W\otimes_RK_\bullet(\gr^L(N);E),
\end{equation}
by which, abusing notation, we mean the usual Koszul complex induced by the
$L$-symbols of $E$ on $W\otimes_R\gr^L(N)$.

By \cite[Thm.~5.1]{MMW05}, $K_{A,\bullet}(N;\beta)$ is homologically $F$-holonomic for all toric modules $N$. 
Using Proposition~\ref{55} and Lemma~\ref{34}, one shows by induction on the toric length of $N$ combined with a spectral sequence argument as in the proof of Theorem~\ref{50} that $K_{A,\bullet}(N;\beta)$ is homologically $L$-holonomic for all toric modules $N$ and for all $L$.

\begin{dfn}
For a toric module $N$, we call
\[
\cc^L(N):=\cc^L(K_{A,\bullet}(N;\beta))=\sum_{i=0}^d (-1)^i\cc^L(H_{A,i}(N;\beta))
\] 
the {\em $L$-characteristic Euler--Koszul cycle} of $N$. 
Recall Definition~\ref{32}.
We denote the {\em $i$-th Euler--Koszul $L$-multiplicity} and the {\em Euler--Koszul $L$-characteristic} of $N$ along $C$ by
\[
\mu^{L,C}_{A,i}(N;\beta):=\mu^{L,C}_i(K_{A,\bullet}(N;\beta))\index{mLCAi@$\mu^{L,C}_{A,i}(N;\beta)$},\quad
\mu^{L,C}_A(N):=\mu^{L,C}(K_{A,\bullet}(N;\beta))\index{mLCA@$\mu^{L,C}_A(N)$}.
\]
If $N=S_A$ we skip the argument $N$, while if $C=\bar C_A^\tau$ we write $\tau$ rather than $\bar C_A^\tau$ as upper index.\index{mLtAi@$\mu^{L,\tau}_{A,i}(N;\beta)$}\index{mLtA@$\mu^{L,\tau}_A(N)$}\index{mLCA@$\mu^{L,C}_A(N;\beta)$}\index{mLtA@$\mu^{L,\tau}_A(N;\beta)$}
\end{dfn}

The independence of $\beta$ suggested by the notation $\mu^{L,C}_A(N)$ will be established in Theorem~\ref{50} below. 
Lemma~\ref{34} and the long exact Euler--Koszul homology sequence imply the following statement.

\begin{lem}\label{28}
The Euler--Koszul $L$-characteristic $\mu^{L,C}_A$ along $C$ is additive.\qed
\end{lem}

\subsection{Euler--Koszul characteristic and intersection multiplicity}

We now interpret $\mu^{L,\tau}_A$ as an intersection multiplicity (see \cite[Ch.~V, \S B]{Ser65}).
That will lead to an explicit formula in Theorem~\ref{26}.

\begin{dfn}
For two $W$-modules $M$ and $M'$ and a variety $C\subseteq\Spec(W)$
not strictly contained in an irreducible component of
$\supp(M\otimes_WM')$, the {\em intersection multiplicity of $M$ and
  $M'$ along $C$} is the alternating sum
\[
\chi^C(M,M'):=\sum_i(-1)^i\ell\left(\Tor_i^{W_{P_C}}(M_{P_C},M'_{P_C})\right)\index{cCMM'@$\chi^C(M,M')$}
\]
of lengths of $W_{P_C}$-modules.
If $C=\bar C_A^\tau$ we write $\tau$ rather than $\bar C_A^\tau$ as upper index. \index{ctMM@$\chi^\tau(M,M')$}
\end{dfn}

Lemmas~\ref{31} and \ref{29} and \cite[Ch.~II, Prop.~1.3.1.b]{Sch85} yield the following statements.

\begin{lem}\label{65}
Fix a $W$-module $M$ and a variety $C\subseteq\Spec(W)$.
Whenever defined, the quantity $\chi^C(M,W\otimes_R\gr^L(N))$ is independent of the choice of a good $L$-filtration on the $R$-module $N$ and additive in $N$.\qed
\end{lem}

\begin{thm}\label{50}
For any toric module $N$ and all $\beta\in\CC^d$, 
\[
\mu^{L,C}_A(N)=\chi^C(W/W\ideal{E},W\otimes_R\gr^L(N))
\]
and is hence independent of $\beta$.
For any $i\in\NN$, $\mu^{L,C}_{A,i}(N;\beta)>0$ only if $C=\bar C_A^\tau$ for some $\tau\in\Phi^L_A$, while $\mu^{L,C}_A(N)>0$ if and only if $\dim(N)=d$ and $C=\bar C_A^\tau$ for some $\tau\in\Phi^L_A$.
\end{thm}

\begin{proof}
We consider $K_\bullet=K_{A,\bullet}(N;\beta)$ and fix a variety $C$
for which $\mu^{L,C}_A(N)$ is defined. 
By Corollary~\ref{41} and equation~\eqref{30}, there is a sequence of homologically graded
$W_{P_C}$-modules $K^1,\ldots,K^r$ with $W_{P_C}$-linear endomorphisms
$d^i\colon K^i\to K^i$ such that
\[
K^1=W_{P_C}\otimes_RK_\bullet(\gr^L(N);E),
\]
$K^{i+1}=H(K^i,d^i)$ for $i=1,\ldots,r-1$, and $K^r=\gr^L(H_{A,\bullet}(N;\beta))_{P_C}$ as homologically graded modules.
Thus for all $\beta\in\CC^d$, $\mu^{L,C}_A(N)=\chi(K^r)=\cdots=\chi(K^1)$ is independent of $\beta$ by additivity of $\chi(K):=\sum_i(-1)^i\ell(H_i(K,d))$.
Since $E$ is a regular sequence in $W$, Koszul homology against $E$ computes $\Tor$ over both $W$ and $W_{P_C}$ and hence
\begin{align*}
\chi(K^1)
&=\sum_{i=1}^d(-1)^i\ell\left(\Tor_i^{W_{P_C}}(W_{P_C}/W_{P_C}\ideal{E},W_{P_C}\otimes_R\gr^L(N))\right)\\
&=\chi^C(W/W\ideal{E},W\otimes_R\gr^L(N)).
\end{align*}

By Lemmas~\ref{28} and \ref{65} we may assume now that $N=S_\theta$
for some $\theta\in\Phi^\boldzero_A$ with the standard good
$L_\del$-filtration induced by the single generator $1\in S_\theta$.
Then the complex $K^1$ above is the Koszul complex of $E$ on
$W_{P_C}\otimes_R S_\theta^L$. 
By \cite[Prop.~1.6.5]{BH93}, its homology is annihilated by both $\gr^L(I_A^\theta)=R\cdot I_\theta^L+J_\theta\supseteq I_A^L$ and $E$. 
Then equation~\eqref{67} shows that $\mu^{L,C}_{A,i}(S_\theta;\beta)>0$ for some $i$ implies that $C=\bar C_A^\tau$ is an irreducible component of $\Spec(W/W\ideal{E}\otimes_RS_A^L)$ with $\theta\supseteq\tau\in\Phi^L_A$. 
By Theorem~\ref{13},
\[
\codim(W/W\ideal{E})+\codim(W\otimes_RS_\theta^L)=d+n-\dim(S_\theta)\ge
n=\codim(\bar C_A^\tau).
\]
So Serre's Intersection Theorem~\cite[Ch.~V, \S C.1, Thm.~1]{Ser65} shows that $\mu^{L,\tau}_A(S_\theta)>0$ if and only if $\dim(S_\theta)=d$, which is to say $\theta=A$.
\end{proof}

In the special case $L=F$, Theorem~\ref{50} yields the holonomicity of Euler--Koszul homology of toric modules proved in \cite[Thm.~5.1]{MMW05}.

In order to discuss $\cc^L(M_A(\beta))$ we need more information about
$\mu^{L,\tau}_{A,0}(\beta)$:
Theorem~\ref{50} does not state that $\mu^{L,\tau}_{A,0}(\beta)>0$.
We show next that for generic $\beta$ this inequality does indeed hold. 
In Subsections~\ref{rigidity} and \ref{generic} we will show that the
genericity assumption is not necessary and in fact identify a global
combinatorial lower bound for $\mu^{L,\tau}_{A,0}(\beta)$, similar to
classical statements about the holonomic rank of $M_A(\beta)$ in \cite{GKZ89}. 

\begin{cor}\label{33}
For generic (more precisely, not rank-jumping) $\beta\in\CC^d$
\[
\mu^{L,\tau}_{A,0}(\beta)=\mu^{L,\tau}_A=\chi^\tau(W/W\ideal{E},W\otimes_RS_A^L)>0
\]
for all $\tau\in\Phi_A^L$.
In particular, $\ch^L(M_A(\beta))=\bigcup_{\tau\in\Phi^L_A}\bar C_A^\tau$.
\end{cor}

\begin{proof}
By assumption on $\beta$ and \cite[Thm.~6.6]{MMW05}, the Euler--Koszul complex
$K_{A,\bullet}(S_A;\beta)$ is a resolution of $M_A(\beta)$.
Hence $\cc^L(K_{A,\bullet}(S_A;\beta))=\cc^L(M_A(\beta))$ and $\mu^{L,\tau}_{A,0}(\beta)=\mu^{L,\tau}_A>0$ by Theorem~\ref{50}.
So the last claim follows from Proposition~\ref{55}.
\end{proof}

Recall that $\tilde S_A$ is the normalization of $S_A$. 
In Section~\ref{rigidity}, we show that $\mu^{L,\tau}_{A,0}(\beta)>0$ for all $\beta$.
This is done by reducing to $\tilde S_A$ and applying the following statement.

\begin{cor}\label{49}
For all $\tau\in\Phi_A^L$,
\[
\mu^{L,\tau}_{A,0}(\tilde S_A;\beta)=\mu^{L,\tau}_A(\tilde S_A)=\mu^{L,\tau}_A=\chi^\tau(W/W\ideal{E},W\otimes_RS_A^L)>0.
\]
\end{cor}

\begin{proof}
Since $\tilde S_A$ is a Cohen--Macaulay module of full dimension over $S_A$, its Euler--Koszul complex is a resolution of $H_{A,0}(\tilde S_A;\beta)$ by \cite[Thm.~6.6]{MMW05}.
Hence $\mu^{L,\tau}_{A,i}(\tilde S_A;\beta)=0$ for all $i>0$ which proves the first equality.
There is a short exact sequence 
\begin{equation}\label{51}
\xymat{
0\ar[r]&S_A\ar[r]&\tilde S_A\ar[r]&Q\ar[r]&0
}
\end{equation}
where $\dim(Q)<d$ and the rest of the claim follows from
Theorem~\ref{50} and Corollary~\ref{33}.
\end{proof}

\subsection{Rigidity of Euler--Koszul multiplicities}\label{rigidity}

In this section we show in Theorem~\ref{60} that the strict positivity of $\mu^{L,\tau}_{A,0}(\beta)$ in Corollary~\ref{33} holds without the genericity assumption.
As a consequence, we obtain a complete description of the $L$-characteristic variety and the slopes of the $A$-hypergeometric system along coordinate varieties at the origin.

\begin{dfn}
Let $\tau$ be a subset of the column set of $A$.
We set $W_\tau=\CC[x_\tau,\del_\tau]$\index{Wt@$W_\tau$} so that
$W/WI_A^\tau=\CC[x_{\ol\tau}]\otimes_\CC (W_\tau\otimes_{R_\tau}S_\tau)$.
\end{dfn}

Let $g\in\GL(\ZZ,d)$ and put $A'=gA$ and $\beta'=g\beta$. 
Then note that $S_A=S_{A'}$ and $E'-\beta'=g(E-\beta)$. 
Hence, the Euler--Koszul complexes $K_{A,\bullet}(S_A;\beta)$ and $K_{A',\bullet}(S_{A'};\beta')$ are homotopy equivalent. 
In particular, when investigating homotopy-invariant properties of $K_{A,\bullet}(S_A;\beta)$ one may replace $A$ by $gA$.  

\begin{lem}\label{61}
Let $\theta\in\Phi^\boldzero_A$ and assume that the top $\dim(\theta)$ rows $A'$ of $\theta$ are a lattice basis for the $\ZZ$-row span of $\theta$.
Then, for all $\tau\in\Phi_A^L$, 
\[
\mu^{L,\tau}_{A,i}(S_\theta;\beta)=
\begin{cases}
\mu^{L',\tau'}_{A',i}(S_{\theta'};\beta')&\text{if }\tau\subseteq\theta;\\
0&\text{if }\tau\not\subseteq\theta;
\end{cases}
\]
where $L'=(L_{x_{A'}},L_{\del_{A'}})$, $\tau=\tau'$, $\theta'=A'$, and $\beta':=(\beta_1,\dots,\beta_{\dim(\theta)})$.
\end{lem}

\begin{proof}
By \cite[Lem.~4.8]{MMW05}, $H_{A,i}(S_\theta;\beta)=\CC[x_{\ol\theta}]\otimes_\CC H_{A',i}(S_{\theta'};\beta')$ and hence 
\[
\gr^L(H_{A,i}(S_\theta;\beta))=\CC[x_{\ol\theta}]\otimes_\CC\gr^{L'}(H_{A',i}(S_{\theta'};\beta')).
\]
For $j\in\tau\smallsetminus\theta$, $\del_j$ is both annihilator and unit on $\gr^L(H_{A,i}(S_\theta;\beta))_{P_\tau}$ and hence $\tau\subseteq\theta$ if $\mu^{L,\tau}_{A,i}(S_\theta;\beta)\ne0$.
Assume now that $\tau\subseteq\theta$. 
Then 
\[
W_{P_\tau}\otimes_W\gr^L(H_{A,i}(S_\theta;\beta))=\CC(x_{\ol\theta})\otimes_\CC\left((W_{A'})_{P_\tau}\otimes_{W_{A'}}\gr^{L'}(H_{A',i}(S_{\theta'};\beta'))\right)
\]
where the symbol $P_{\tau}$ is used for the two ideals in $W_A$ and $W_{A'}$ from Definition~\ref{74} induced by $\tau\in\Phi_A^L$ and $\tau\in\Phi_{A'}^{L'}$. 
In particular, the length of $W_{P_\tau}\otimes_W\gr^L(H_{A,i}(S_\theta;\beta))$ over $W_{P_\tau}$ is the same as the length of $(W_{A'})_{P_\tau}\otimes_{W_{A'}}\gr^{L'}(H_{A',i}(S_{\theta'};\beta'))$ over $(W_{A'})_{P_\tau}$. 
\end{proof}

\begin{thm}\label{60}
For any toric module $N$ and for every parameter $\beta\in\CC^d$, if $\mu^{L,\tau}_{A,i}(N;\beta)>0$ for some $i>0$ then $\mu^{L,\tau}_{A,0}(N;\beta)>0$.
\end{thm}

\begin{proof}
Suppose there is a counterexample $(A,\beta,L,\tau,N)$ to the theorem.
This means that $\mu^{L,\tau}_{A,i}(N;\beta)\neq0$ for some $i>0$ while $\mu^{L,\tau}_{A,0}(N;\beta)=0$.
We choose a minimal counterexample in the sense that
\begin{asparaenum}[1.]
\item $A$ has minimal number of columns,
\item for a fixed such $A$, $\dim(N)$ is minimal, and
\item no quotient of $N$ is part of such a counterexample.
\end{asparaenum}
This last choice can be made since toric modules are Noetherian. 
Let $\theta\in\Phi^\boldzero_A$ such that $I_A^\theta$ is an associated prime of $N$.
Then there is a short exact sequence 
\[
\xymat{
0\ar[r]&S_\theta(-\bolda)\ar[r]&N\ar[r]&N'\ar[r]&0
}
\]
where $\bolda\in\ZZ^d$ is the degree of the image of $1\in S_\theta$ within $N$.
The long exact Euler--Koszul homology sequence shows that $H_{A,0}(N';\beta)$ is a quotient of $H_{A,0}(N;\beta)$. 
Recall that by Lemma~\ref{31} the $L$-characteristic cycle is additive.
Since by hypothesis $\mu^{L,\tau}_{A,0}(N;\beta)=0$, we conclude that $\mu^{L,\tau}_{A,0}(N';\beta)=0$ as well. 

By assumption, $(A,\beta,L,\tau,N')$ is not a counterexample and so $\mu^{L,\tau}_{A,i}(N';\beta)=0$ for all $i$.
It follows that the $\mu^{L,\tau}_{A,i}$-vanishing patterns of $N$ and of $S_\theta(-\bolda)$ coincide. 
In particular, $S_\theta(-\bolda)$ is a counterexample. 
Since $\dim(S_\theta)\le\dim(N)$, and since $S_\theta$ is a domain, $S_\theta(-\bolda)$ is actually a minimal counterexample.
Since $H_{A,i}(S_\theta(-\bolda);\beta)=H_{A,i}(S_\theta;\beta-\bolda)$ up to a $\ZZ^d$-shift we may assume that $N=S_\theta$. 
But if $\theta\ne A$ then Lemma~\ref{61} yields a counterexample whose matrix $A'$ has strictly fewer columns than $A$.
Therefore $\theta=A$ by minimality of $A$ and $(A,\beta,L,\tau,S_A)$ is a minimal counterexample.

By Corollary~\ref{49}, $\mu^{L,\tau}_{A,0}(\tilde S_A;\beta)>0$ where $\tilde S_A$ is the normalization of $S_A$. 
For appropriate $\bolda\in\ZZ^d$, there is a short exact sequence
\[
\xymat{
0\ar[r]&\tilde S_A(-\bolda)\ar[r]&S_A\ar[r]&Q\ar[r]&0
}
\]
where $Q$ is not a counterexample as $\dim(Q)<d$.
By additivity of the $L$-characteristic cycle,
$\mu^{L,\tau}_{A,0}(S_A;\beta)=0$ implies
$\mu^{L,\tau}_{A,0}(Q;\beta)=0$ and so, as $Q$ is not a counterexample,
$\mu^{L,\tau}_{A,i}(Q;\beta)=0$ for all $i$. 
This yields $\mu^{L,\tau}_{A,i}(S_A;\beta)=\mu^{L,\tau}_{A,i}(\tilde S_A;\beta)$ for all $i$ which is a contradiction for $i=0$.
\end{proof}

Finally, from Theorem~\ref{60}, Corollary~\ref{33} and
Theorem~\ref{50}, we obtain the $L$-characteristic variety of the
$A$-hypergeometric system and, as a consequence, the slopes of the
$A$-hypergeometric system along coordinate varieties at the origin.

\begin{cor}\label{23}
For all $\tau\in\Phi_A^L$ and $\beta\in\CC^d$,
$\mu_{A,0}^{L,\tau}(\beta)>0$. In consequence, the $L$-characteristic
variety of $M_A(\beta)$ consists for every $\beta$ precisely of
all conormal closures of the torus orbits indexed by the faces in the
$(A,L)$-umbrella:
\[\pushQED{\qed}
\ch^L(M_A(\beta))=\bigcup_{\tau\in\Phi^L_A}\bar C_A^\tau.\qedhere
\]
\end{cor}

\begin{cor}\label{39}
For $L=pF+qV$ as in \eqref{73}, $L'=p'/q'$ is a slope of $M_A(\beta)$
along $Y$ at $0\in Y$ if and only if $\Phi^{L}_A$ jumps at $L=L'$.
\qed
\end{cor}

\subsection{Euler--Koszul characteristic and volume}\label{volume}

In this section we develop an explicit combinatorial formula for the Euler--Koszul characteristic $\mu_A^{L,\tau}$ for all $\tau\in\Phi_A^L$, Theorem~\ref{26} below. 
For generic $\beta$, this formula determines the characteristic cycle
of $M_A(\beta)$ as stated in Corollary~\ref{23} below, but in general
it provides only a lower bound.

Note that $\rk(M_A(\beta))=\mu^{F,\emptyset}_{A,0}(\beta)$, and that $\mu^{L,\tau}_{A,0}(\beta)=\mu^{L,\tau}_A$ if $\beta$ is not a rank-jumping parameter by \cite[Thm.~6.6]{MMW05}. 
It is a classical result that $\rk(M_A(\beta))=\vol_{\ZZ^d}(\Delta^F_A)$ for generic $\beta$ (see \cite{GKZ89,Ado94}). 
Theorem~\ref{26} below contains this result as a special case when $\tau=\emptyset$ and $L=F$. 
Our generalization of the fact that generic rank equals volume reads
\[
\mu_A^{L,\emptyset}=\vol_{\ZZ^d}\left(\bigcup_{\tau\in\Phi_A^{L,d-1}}(\Delta^F_0\smallsetminus\conv(\tau))\right)\,.
\]

\medskip

This paragraph outlines our strategy towards Theorem~\ref{26}.  By
Theorem~\ref{50}, $\mu^{L,\tau}_A$ is an intersection multiplicity and
hence additive. 
Applying $\mu_A^{L,\tau}$ to the composition chain~\eqref{56} of $S_A^L$ yields
\begin{align}\label{75}
\mu_A^{L,\tau}=\chi^{\tau}(W/\ideal{E},S_A^L)
&=\sum_{i=1}^l\chi^{\tau}(W/\ideal{E},S_{\tau_i})\\
\nonumber&=\sum_{\tau\subseteq\tau'\in\Phi^{L,d-1}_A} \nu_A^{L,\tau'}\cdot
\chi^{\tau}(W/\ideal{E},S_{\tau'})\\
\nonumber&=\sum_{\tau\subseteq\tau'\in\Phi^{L,d-1}_A}[\ZZ^d:\ZZ\tau']\cdot
\mu_A^{L,\tau}(S_{\tau'}).
\end{align}
For the third equality note that $(R/I^{\tau'}_A)_{I_A^{\tau}}=0$
unless $\tau\subseteq\tau'$, and then use $\dim(S_{\tau'})=\dim(\tau')+1$ in Theorem~\ref{50}. 
Proposition~\ref{58} gives the last equality.
By equation~\eqref{75} it is sufficient to know the multiplicities $\mu_A^{L,\tau}(S_{\tau'})$ for $\tau\subseteq\tau'\in\Phi^{L,d-1}_A$.
To determine these we first reduce to $\tau'=A$ and then to $\tau=\emptyset$ while controlling the appropriate intersection multiplicity. 
Finally, we give a combinatorial description in that case. 
Some of our arguments are similar to \cite[\S\S2.3--2.4]{GKZ89} and \cite[\S5]{Ado94}.

Let $\Phi_A^L\ni\tau\subseteq\tau'\in\Phi_A^{L,d-1}$ and abbreviate
$\mu:=\mu_A^{L,\tau}(S_{\tau'})=\chi^\tau(W/\ideal{E},S_{\tau'})$.
Since $E$ is a regular sequence in $W$, $\mu$ can be computed from
the Koszul complex of $E$ on $W_{P_{\tau}}\otimes_RS_{\tau'}$.
Since $\del_{\ol\tau'}$ is zero on $S_{\tau'}$ one can erase the
$\del_{\ol\tau}$-terms in $E$.
Now the (canonical basis of the) ambient lattice $\ZZ^d=\ZZ A$ 
can be replaced by (a lattice basis of) $\ZZ\tau'$; this
leaves the quantities $W$, $\ideal{E}$, $I_\tau$ and $S_{\tau'}$ invariant.
We have thus reduced the problem to the case $\tau'=A$.

\medskip 

We now reduce to the case $\tau=\emptyset$. 
By the Fourier transform we may use $\del$ as base variables and identify $-x_j\in W$ with the (symbol of the) partial derivation along $\del_j$. 
Since $\chi^{\tau}$ is determined at a generic point of $C^\tau_A$, we may replace $W$ by $W'=W[\del_{\tau}^{-1}]$ whose spectrum $T^*U$ is the cotangent space of $U=\Spec(R[\del_{\tau}^{-1}])$.

We now modify the pair $(\tau',\tau)$ into a more convenient one, while keeping track of $\mu$.
A row operation $(\hat \tau',\hat\tau)=g(\tau',\tau)$ defined by $g\in\GL(\ZZ,d)$ amounts to a basis change in the $d$-torus $\TT$, so $\mu=\mu_{\hat\tau'}^{L,\hat\tau}(S_{\hat\tau'})$.
On the other hand, for a triple $j\in\tau'$, $i\in\tau$, $m\in\ZZ$, the elementary column operation $\hat\bolda_j=\bolda_j+m\bolda_i$ transforming $\tau'$ into $\hat\tau'$ corresponds to the automorphism $\hat\del_j=\del_j\del_i^m$ of $U$. 
The induced map on $W'$ obviously transforms $I_{\tau'}$ into $I_{\hat\tau'}$. 
Note that $x_i\del_i=\hat x_i\hat\del_i+m\hat x_j\hat\del_j$ and $x_k\del_k=\hat x_k\hat \del_k$ for $k\ne i$.
Thus, $\mu=\mu_{\hat\tau'}^{L,\hat\tau}(S_{\hat\tau'})$. 
In other words, all elementary column operations $\bolda_j\mapsto\bolda_j+m\bolda_i$ with $i\in\tau$, all column switches $\bolda_i\leftrightarrow\bolda_j$ with $i,j\in\tau$ or $i,j\not\in\tau$, and all $\ZZ$-invertible row operations leave $\mu$ invariant.

After a suitable sequence of such transformations, we may assume that 
\[
\tau=\begin{pmatrix}I&0\\0&0\end{pmatrix},
\quad
I=\begin{pmatrix}k_1&&0\\&\ddots&\\0&&k_{e}\end{pmatrix}
\]
where $k_i\in\NN\smallsetminus\{0\}$. Note that 
the product $k_1\cdots k_{e}$ is the index
$[(\ZZ \tau'\cap\QQ\tau):\ZZ\tau]$ of $\ZZ\tau$ in
$\ZZ \tau'\cap\QQ\tau$. 

We next consider the $k_1$-fold covering space $\kappa_1\colon U\to U$ induced by $\hat\del_1^{k_1}=\del_1$. 
Let $(\hat\tau,\hat\tau')$ be the matrices obtained from $(\tau,\tau')$ by dividing the first column by $k_1$ and let $\hat E$ be the Euler vector fields of $\hat\tau$. 
As one checks, $\hat x_1\hat\del_1=k_1x_1\del_1$ and so $\kappa_1(\Var(\hat E))=\Var(E)$, $\kappa_1(\Var(I_{\hat\tau'}))=\Var(I_{\tau'})$, and $\kappa_1(C_A^{\hat\tau})=C_A^\tau$. 
Moreover, the degree of $\kappa_1$ is equal to $k_1$ on $C_A^\tau$, and equal to $1$ on both $\Var(E)$ and $\Var(I_{\tau'})$. 
By Example 8.2.5 in \cite{Ful98},
\[
\mu=k_1\cdot\chi^{\hat\tau}(W/\ideal{\hat E},S_{\hat\tau'}).
\]
Analogous transformations for $k_2,\dots,k_{e}$ followed by suitable
elementary column operations as above yield
\[
\hat\tau=\begin{pmatrix}I&0\\0&0\end{pmatrix},
\quad
\hat\tau'=\begin{pmatrix}I&0&0\\0&0&\check \tau\end{pmatrix},
\quad
I=\begin{pmatrix}1&&0\\&\ddots&\\0&&1\end{pmatrix},
\]
where now 
\[
\mu=k_1\cdots k_e\cdot\chi^{\hat\tau}(W/E_{\hat\tau'},S_{\hat\tau'})
=[(\ZZ\tau'\cap\QQ\tau):\ZZ\tau]\cdot\chi^{\hat\tau}(W/E_{\hat\tau'},S_{\hat\tau'}).
\]
Note that $\check\tau$ is pointed since $\tau$ is a face of $\tau'$.
Let $\check E$ be the Euler vector fields of $\check\tau$.  Then on
$T^*U$, $\ideal{E}=\ideal{\hat x_1\hat\del_1,\dots,\hat
x_{e}\hat\del_{e},\check E}=\ideal{\hat x_1,\dots,\hat x_{e},\check
E}$ as well as $I_{\hat
\tau'}=\ideal{\hat\del_{e+1},\dots,\hat\del_{|\tau|}}+I_{\check\tau}$.
This establishes the announced reduction to the case $\tau=\emptyset$,
\[
\mu=[(\ZZ\tau'\cap\QQ\tau):\ZZ\tau]\cdot\chi^{\emptyset}(W/\ideal{\check E},S_{\check\tau})
=[(\ZZ\tau'\cap\QQ\tau):\ZZ\tau]\cdot\mu_{\check\tau}^{L,\emptyset}(S_{\check\tau}).
\]
The constructed lattice $\ZZ\check \tau\subseteq\ZZ\tau'\subseteq \ZZ^d$ corresponds to a splitting of the natural projection
\[
\xymat{
\pi_{\tau,\tau'}\colon \ZZ\tau'\ar@{->>}[r]^-{}&\ZZ\tau'/(\ZZ\tau'\cap\QQ\tau)\index{pt@$\pi_{\tau,\tau'}$}}
\]  
which identifies the positive semigroups $\NN\check\tau$ and
$\pi_{\tau,\tau'}(\NN\tau')$.

\begin{dfn}\label{98}
In a lattice $\Lambda$, the volume function $\vol_\Lambda$\index{volL@$\vol_\Lambda$} is normalized so that the unit simplex of $\Lambda$  has volume $1$.
We abbreviate $\vol_{\tau,\tau'}:=\vol_{\pi_{\tau,\tau'}(\ZZ\tau')}$\index{voltt@$\vol_{\tau,\tau'}$}.
\end{dfn}

We continue under the assumption that $\tau=\emptyset$, $\check\tau=\tau'=A$, and compute the intersection multiplicity of $\Var(E)$ with $\Var(I_A)$ along $C_A^\emptyset=X$ inside $T^*U=\Spec(W')$,
\[
\check\mu:=\mu_{\check\tau}^{L,\emptyset}(S_{\check\tau})=\chi^{\emptyset}(W/\ideal{\check E},S_{\check\tau}).
\]
Since for $y\in(\CC^*)^{|\check\tau|}$ the linear polynomials $x_1-y_1,\ldots,x_{|\check\tau|}-y_{|\check\tau|}$ form a regular sequence on $W/\ideal{E}$, on $S_{\check\tau}$ and on $W/P_\emptyset$, we may replace $W$ by $R$, $E_i$ by $\bar E_i:=\sum_j a_{i,j}y_j\del_j$, and $C_A^\emptyset$ by $\{0\}=\Var(\del_{\check\tau})$. 
By the sequence~\eqref{51}, we may further replace $S_{\check\tau}$ by its normalization $\tilde S_{\check\tau}$. 
Thus, we have reduced to computing
\[
\check\mu=\chi^{\{0\}}(R/\ideal{\bar E},\tilde S_{\check\tau}).
\]
For generic $y$, the function on $\Spec(\tilde S_{\check\tau})$ 
\[
f:=\sum_{j=1}^{|\check \tau|}y_j\del_j=\sum_{j=1}^{|\check \tau|}y_jt^{\bolda_j}
\]
is by \cite[Thm.~6.1]{Kou76} and its proof Newton nondegenerate at the origin. 
We can interpret $\bar E$ as functions on $\Spec(\tilde S_{\check\tau})$,
\[
\bar E_i=t_i\frac{\del f}{\del t_i}.
\]
Let $\frakm_{\check\tau}\subseteq S_{\check\tau}$ and
$\tilde\frakm_{\check\tau}\subseteq\tilde S_{\check\tau}$ be the
maximal ideals at $0$. 
In the special case where $\check\tau$ is the unit matrix, and where hence $\tilde S_{\check\tau}$ is a polynomial ring, \cite[Thm.~2.8]{Kou76} states that the Koszul complex of $\bar E$ on $\tilde S_{\check\tau}$ is acyclic in positive dimension when completed at $\tilde\frakm_{\check\tau}$. 
The $\frakm_{\check\tau}$- and $\tilde\frakm_{\check\tau}$-adic topologies on $\tilde S_{\check\tau}$ are equivalent and completion at $\tilde\frakm_{\check\tau}$ can be replaced by completion at $\frakm_{\check\tau}$. 
By faithful flatness of completion this implies that the Koszul complex of $\bar E$ on $(\tilde S_{\check\tau})_{\frakm_{\check\tau}}$ is acyclic in positive dimension and hence
\[
\check\mu=\dim_\CC((\tilde S_{\check\tau})_{\frakm_{\check\tau}}/\ideal{\bar E}).
\]
Since this dimension is finite, localization can be replaced by completion at $\frakm_{\check\tau}$ or equivalently at $\tilde\frakm_{\check\tau}$.
Then \cite[Thm.~A.I]{Kou76} states that this dimension equals
\[
\check\mu=\vol_{\ZZ\check\tau}(\Gamma_-(f)),\quad\Gamma_-(f)=\conv(\check\tau\cup\{0\})\smallsetminus\conv(\check\tau).
\] 
One checks that Kouchnirenko's result applies to a general normal semigroup ring $\tilde S_{\check\tau}$ for which $\check\tau$ is pointed (see \cite[\S2.12]{Kou76}).
The following definition serves to interpret this volume in terms of the original matrix $A$, and mirrors the one given in \cite[\S2.1]{GKZ89}.

\begin{dfn}
For $\Phi_A^L\ni\tau\subseteq\tau'\in\Phi_A^{L,d-1}$, define the polyhedra
\[
P_{\tau,\tau'}:=\conv(\pi_{\tau,\tau'}(\tau'\cup\{0\}))\index{Ptt@$P_{\tau,\tau'}$},\quad Q_{\tau,\tau'}:=\conv(\pi_{\tau,\tau'}(\tau'\smallsetminus\tau))\index{Qtt@$Q_{\tau,\tau'}$}.
\]
\end{dfn}

We are now ready to give the promised multiplicity formula.

\begin{thm}\label{26}
For all $\tau\in\Phi_A^L$, the multiplicity of $\bar C_A^\tau$ in $\cc^L(K_{A,\bullet}(S_A;\beta))$ is 
\[\pushQED{\qed}
\mu_A^{L,\tau}=
\sum_{\tau\subseteq\tau'\in\Phi^{L,d-1}_A}[\ZZ^d:\ZZ\tau']\cdot[(\ZZ\tau'\cap\QQ\tau):\ZZ\tau]\cdot\vol_{\tau,\tau'}(P_{\tau,\tau'}\smallsetminus Q_{\tau,\tau'}).\qedhere
\]
\end{thm}

\begin{cor}\label{20}
If $\tau=\emptyset$, then the local degree of $S_A^L$ at the origin equals 
\[
\mu_A^{L,\emptyset}=
\vol_{\ZZ^d}\left(\bigcup_{\tau'\in\Phi_A^{L,d-1}}\left(\Delta^F_{\tau'}\smallsetminus\conv(\tau')\right)\right)\le\vol_{\ZZ^d}(\Delta^F_A).
\]
\end{cor}

\begin{proof}
Let $\tau'\in\Phi_A^{L,d-1}$.
For $\tau=\emptyset$, $\vol_{\tau,\tau'}=\vol_{\ZZ\tau'}$, $P_{\tau,\tau'}=\Delta^F_{\tau'}$, and $Q_{\tau,\tau'}=\conv(\tau')$.
Clearly, $[(\QQ\tau\cap\ZZ\tau'):\ZZ\tau]=1$ and $\vol_{\ZZ\tau'}(\Delta^F_{\tau'})=\vol_{\ZZ^d}(\Delta^F_{\tau'})/[\ZZ^d:\ZZ\tau']$.
Then Theorem~\ref{26} implies the equality while the inequality is obvious.
\end{proof}

\begin{rmk}\label{99}
If all facets $\tau'$ of the $(A,L)$-umbrella are $F$-homogeneous then
the volumes $\vol_{\ZZ^d}(\conv(\tau'))$ are zero and one obtains 
\[
\mu_A^{L,\emptyset}=
\vol_{\ZZ^d}\left(\bigcup_{\tau'\in\Phi_A^{L,d-1}}\Delta^F_{\tau'}\right).
\]
In particular, $\mu^{F,\emptyset}_A=\vol_{\ZZ^d}(\Delta^F_A)$. 
\end{rmk}

\begin{exa}\label{100}
We continue Example~\ref{97} and investigate characteristic
cycles. The following table, whose rows are indexed by the three
weight vectors $L=F+tV$ considered in Example~\ref{97}, lists
in its columns the nonzero
multiplicities $\mu^{L,\tau}_{A,0}(\beta)$ for generic
$\beta$. Genericity implies that $\mu^{L,\tau}_{A,i}(\beta)=0$ for
$i>0$, so that $\mu^{L,\tau}_{A,0}(\beta)=\mu^{L,\tau}_{A}$.
\begin{figure}[ht]
\caption{Euler--Koszul multiplicities in Example~\ref{100}}
\begin{tabular}{|c|c|c|c|c|c|c|c|c|c|c|}
\hline
\backslashbox{$t$}{$\tau$}&$\emptyset$&\{1\}&\{2\}&\{3\}&\{4\}&\{1,3\}&\{1,4\}&\{2,3\}&\{2,4\}&\{3,4\}\\
\hline
  0&13           &12   &1    &     &13   &       &12     &       &1
  &\\
\hline
  1&11           &3    &1    &10   &8    &3      &       &       &1
  &7\\
\hline
  4&5            &3    &2    &5    &     &3      &       &2      &
  &\\
\hline
\end{tabular}
\end{figure}

The matrix $A$ permits two rank-jumping parameters, $\calE(A)=\{(2,1),
(3,1)\}$. In both cases, $H_{A,1}(S_A;\beta)\cong
D/\ideal{\del_1,\ldots,\del_n}$ while $H_{A,2}(S_A;\beta)\cong0$. As
$\ideal{\del}=P_{\emptyset}$, 
with $\beta\in\calE(A)$ one has
$\mu^{L,\tau}_{A,1}(\beta)=1$ if $\tau=\emptyset$, and
$\mu^{L,\tau}_{A,1}(\beta)=0$ otherwise. (This holds for \emph{all}
$L\in\QQ^{2n}$ as long as $L_x+L_\del>0$ since
$\gr^L(\ideal{\del_1,\ldots,\del_n})$ is independent of $L$).  In
accordance with Corollary~\ref{49}, the multiplicities
$\mu^{L,\tau}_{A,0}(\beta)$ for $\beta\in\calE(A)$ are given by the
data in the table above, with $\mu^{L,\emptyset}_{A,0}(\beta)$
incremented by one. This behavior is typical in dimension two as the
following proposition shows.
\end{exa}

\begin{prp}\label{101}
If $d=2$ then for every rank-jumping $\beta$ and for all $L$
with $L_x+L_\del>0$ one has
$\mu^{L,\tau}_{A,1}(\beta)=1$ if $\tau=\emptyset$ and
$\mu^{L,\tau}_{A,1}(\beta)=0$ otherwise. 
\end{prp}

\begin{proof}
Recall that $\tilde S_A$ is the Cohen--Macaulay $S_A$-module $\bigoplus_{\bolda \in \QQ_+A\cap\ZZ^d}t^\bolda$. 
We first construct the minimal toric submodule of $\tilde S_A$ that contains $S_A$ and satisfies Serre's condition $S_2$. 
To that end, apply the local cohomology functor $H^\bullet_\frakm(-)$ to the sequence~\eqref{51}; we obtain $H^0_\frakm(Q)=H^1_\frakm(S_A)$. 
Let $C$ be the (toric) preimage of $H^0_\frakm(Q)$ under the projection $\tilde S_A\to Q$; as $C\subseteq \tilde S_A$ we have $H^0_\frakm(C)=0$. 
Applying $H^\bullet_\frakm(-)$ to
\begin{equation}\label{102}
\xymat{
0\ar[r]&S_A\ar[r]&C\ar[r]&H^0_\frakm(Q)\ar[r]&0
}
\end{equation}
we find an exact piece $0\to H^0_\frakm(H^0_\frakm(Q))\to H^1_\frakm(S_A)\to H^1_\frakm(C)\to H^1_\frakm(H^0_\frakm(Q))=0$. 
As the second arrow is an isomorphism by construction, the two-dimensional module $C$ has depth two and is hence Cohen--Macaulay.

Application of the Euler--Koszul functor to the short exact sequence~\eqref{102} shows that $H_{A,1}(S_A;\beta)\cong H_{A,2}(H^0_\frakm(Q);\beta)$. 
Since $H^0_\frakm(Q)$ is a toric Artinian quotient of $C$, it has a toric filtration whose quotients
are of the form $R/\ideal{\del_1,\ldots,\del_n}\cdot t^\beta$ where $\beta\in\calE(A)$ and each such $\beta$ occurs exactly once (see \cite[Thm.~6.6,Thm~.9.1]{MMW05}). 
Hence, $H_{A,1}(S_A;\beta)\cong D/\ideal{\del_1,\ldots,\del_n}$ if $\beta\in\calE(A)$ and zero
otherwise. 
It follows as in the example above that $\gr^L(H_{A,1}(S_A;\beta))\cong W/\ideal{\del_1,\ldots,\del_n}$ in the nonvanishing case.
\end{proof}

\begin{rmk}
\begin{asparaenum}
\item The module $C$ constructed in the above proof is actually a ring, the
\emph{ideal transform} of $S_A$ relative to $\frakm$ from
\cite[Thm.~2.2.4]{BS98}.
\item It is suggestive, particularly in the light of
Example~\ref{100}, to view $\mu_A^{L,\emptyset}$ as
``$L$-rank''. Specifically, one might speculate whether the $L$-rank
detects rank jumps in the sense that
$\mu_{A,i}^{L,\emptyset}(\beta)=0$ for all $i>0$ implies that $\beta$
is not rank-jumping. We think that this is plausible.
\end{asparaenum} 
\end{rmk}

\subsection{Generic Euler--Koszul homology}\label{generic}

In \cite{GKZ89}, I.M.~Gel'fand, M.M.~Kapranov and A.V.~Zelevinski{\u\i} proved that in the projective case the rank of $M_A(\beta)$ is always at least equal to $\vol_{\ZZ^d}(\Delta^F_A)=\mu^{F,\emptyset}_A$, cf.~Remark~\ref{99}; this was generalized by A.~Adolphson in \cite{Ado94}. 
In this section, we prove that every Euler--Koszul multiplicity $\mu^{L,\tau}_{A,0}(\beta)$ is bounded from below by the Euler--Koszul characteristic $\mu^{L,\tau}_A$ for all $\beta$. 
This yields some evidence for the following conjecture.

\begin{cnj}
For fixed $A$, $L$ and $\tau$, $\mu^{L,\tau}_{A,0}(\beta)$ is upper
semicontinuous in $\beta$.
\end{cnj}

A proof of this would generalize \cite[Thm.~2.6, Thm.~7.5]{MMW05} which correspond to $\tau=\emptyset$ and $L=F$. 
However, we do not know how to approach this question since it involves a flat deformation in combination with a specialization. 
The case where $\tau=\emptyset$ and $L=F$ is much easier since there one may skip the computation of the graded by the Cauchy--Kovalevskaya--Kashiwara Theorem (see \cite[Thm.~1.4.14, Thm.~1.4.19]{SST00}).

\begin{thm}\label{47}
For any $\tau\in\Phi_A^L$ and $\beta\in\CC^d$, $\mu^{L,\tau}_{A,0}(\beta)\ge\mu^{L,\tau}_A$. 
Equality holds if $\beta$ is generic (more precisely, not rank-jumping).
\end{thm}

\begin{proof}
The sequence~\eqref{51} yields a four-step exact sequence
\[
\xymat@C=16pt{0\ar[r]&H_{A,1}(Q;\beta)\ar[r]&H_{A,0}(S_A;\beta)\ar[r]&H_{A,0}(\tilde S_A;\beta)\ar[r]&H_{A,0}(Q;\beta)\ar[r]&0}.
\]
By Corollary~\ref{49} and Lemma~\ref{48} below,
\[
\mu^{L,\tau}_{A,0}(S_A;\beta)\geq\mu^{L,\tau}_{A,0}(\tilde S_A;\beta)=\mu^{L,\tau}_A.
\]
The second claim is part of Corollary~\ref{33}.
\end{proof}

\begin{lem}\label{48}
For any toric module $N$ with $\dim(N)<d$, $\tau\in\Phi_A^L$, and $\beta\in\CC^d$,
\[
\mu^{L,\tau}_{A,0}(N;\beta)\le\mu^{L,\tau}_{A,1}(N;\beta).
\]
\end{lem}

\begin{proof}
Let $(A',\beta')=(gA,g\beta)$ for some $g\in\GL(\ZZ,d)$. 
Then the Euler--Koszul complexes $K_{A,\bullet}(N;\beta)$ and $K_{A',\bullet}(N;\beta')$ are homotopic for any toric module $N$.
Hence we may replace $(A,\beta)$ by $(A',\beta')$ and choose $g$ generic. 
Fix the toric module $N$ of dimension less than $d$. 
Then the Koszul complex $K_{A,\bullet}(N;\beta')$ of the endomorphisms $E'-\beta':=(E_1-\beta_1,\dots,E_{d-1}-\beta_{d-1})$ on the $\ZZ^{d}$-graded left $D$-module $D\otimes_RN$, defined as in
Definition~\ref{19}, is homologically $L$-holonomic. 
To see this it is by the long Euler--Koszul homology sequence enough to consider a length-one toric module $N=S_\theta$ for some $\theta\in\Phi^\boldzero_A$. 
Then, using Proposition \ref{55}, genericity of $g$ assures that $W\ideal{I_\theta^L,E'}$ is an $n$-dimensional ideal and hence $\dim\ch^L(H_{A,i}(N;\beta'))\le n$ by the spectral sequence
\[
H_i(\gr^L(K_{A,\bullet}(N;\beta')))\Longrightarrow \gr^L(H_{A,i}(N;\beta')).
\]

We interpret $K_{A,\bullet}(N;\beta)$ as the complex induced by
$E_d-\beta_d$ on $K_{A,\bullet}(N;\beta')$ and abbreviate
$H_i:=H_{A,i}(N;\beta)$ and $H'_i:=H_{A,i}(N;\beta')$.  There is a
double complex spectral sequence abutting to $H_\bullet$ with
$E_1$-term $H'_\bullet$ and differential $d_1$ induced by
$E_d-\beta_d$ which collapses at the $E_2$-term.  In particular,
$H_0=H'_0/(E_d-\beta_d)H'_0$ and there is a short exact sequence of
holonomic $D$-modules
\[
\xymat{0\ar[r]&\ker_{H'_0}(E_d-\beta_d)\ar[r]&H_1\to H'_1/(E_d-\beta_d)H'_1\ar[r]&0}.
\]
Using the exact sequence  
\[
\xymat{
0\ar[r]&
\ker_{H'_0}(E_d-\beta_d)\ar[r]&
H'_0\ar[r]^-{E_d-\beta_d}&
H'_0\ar[r]&H'_0/(E_d-\beta_d)H'_0\ar[r]&0
}
\]
it follows that
\begin{align*}
\cc^L(H_0)&=\cc^L(H'_0/(E_d-\beta_d)H'_0)\\
&=\cc^L(\ker_{H'_0}(E_d-\beta_d))\\
&=\cc^L(H_1)-\cc^L(H'_1/(E_d-\beta_d)H'_1).
\end{align*}
We conclude that $\cc^L(H_1)-\cc^L(H_0)=\cc^L(H'_1/(E_d-\beta_d)H'_1)$ is nonnegative which implies the claim.
\end{proof}

\begin{arxiv}

\section{Projectivized hypergeometric systems}\label{infty}

Here we introduce and study a natural extension of the $A$-hypergeometric system $M_A(\beta)$ on $X=\CC^n$ to a sheaf on $(\PP^1_\CC)^n$ where $\PP^1_\CC=\CC\sqcup\{\infty\}$ in each factor.
We call the resulting extension $\calM_A(\beta)$ of $M_A(\beta)$ the \emph{projectivized $A$-hypergeometric system} defined by $A$ and $\beta$.
We use the results of the previous sections to discuss, in this sequence, $\calM_A(\beta)$ from the points of view of the previous sections:
$(A,L)$-umbrella, $L$-characteristic variety, $L$-characteristic cycle, and, as a special case, its slopes along all projective coordinate subspaces.
Once again the $(A,L)$-umbrella will play a pivotal role in our investigation of $\calM_A(\beta)$.

\begin{ntn}\label{85}
We use the canonical embedding  
\[
\xymat{
\CC^1=\Spec(\CC[x_j])\ar@{^(->}[r]&\Proj(\CC[y_j,y'_j])=\PP^1_\CC=\CC^1\sqcup\{\infty\}
}
\]
induced by $x_j=y_j/y'_j$ to embed $X=\Spec(\CC[x_1,\ldots,x_n])$ into $(\PP^1_\CC)^n=:\bar X$\index{X@$\bar X$} as the complement of the variety of the ideal sheaf generated by $\prod_{j=1}^ny_j'$.
Recall from Notation~\ref{80}, that the symbol $\ol\tau$ denotes the complement of a set $\tau$ in the set $\{1,\ldots,n\}$.
On the other hand, $\bar C$ will always refer to the Zariski closure of a set $C$ in $T^*\bar X$.
We denote by $x'_j:=y'_j/y_j$ the coordinate on $\PP^1_\CC\smallsetminus\{0\}$. 
A subset $\frakP\subseteq\{1,\ldots,n\}$ defines an affine patch $X_\frakP\cong\CC^n$ of $\bar X$ by
\[
X_\frakP=\Spec(\CC[\{x_j\mid j\not\in\frakP\}\cup\{x'_j\mid j\in\frakP\}])\index{XP@$X_\frakP$}.
\]
In particular, $X=X_\emptyset$\index{X@$X$} in this notation.
If $\calF$ is a sheaf on $\bar X$ then we denote its sections over $X_\frakP$ by $F_{(\frakP)}:=\calF(X_\frakP)$\index{FP@$F_{(\frakP)}$}.
On the other hand, if $g$ is an $\bar X$-global section then we denote by $g_{(\frakP)}$\index{gP@$g_{(\frakP)}$} the restriction of $g$ to $X_\frakP$.

Denote by $\calD$\index{D@$\calD$} the sheaf of algebraic $\CC$-linear differential operators on $\bar X$.
Then with the above notation the ring $D_{(\frakP)}=\calD(X_\frakP)$\index{DP@$D_{(\frakP)}$} is the Weyl algebra in the variables $\{x_j\mid j\notin\frakP\}\cup\{x'_j\mid j\in\frakP\}$.
Let $\del'_j$ denote the derivation on $\PP^1_\CC$ relative to $x'_j$.
The sheaf $\calD$ has a special subsheaf of rings $\calR$\index{R@$\calR$}, consisting entirely of the global differential operators in $R=\CC[\del]$.
There are ring isomorphisms 
\begin{equation}\label{87}
\xymat{\rho_\frakP\colon R_{(\emptyset)}=\CC[\del]\ar[r]&\CC[\del_{(\frakP)}]=\CC[\del_{\ol\frakP}][\{{x'_j}^2\del'_j\mid j\in\frakP\}]=R_{(\frakP)}}\index{rP@$\rho_\frakP$}\index{RP@$R_{(\frakP)}$}
\end{equation}
defined by $\del_j\mapsto-{x'_j}^2\del'_j$ for $j\in\frakP$.
Any $\ZZ^d$-graded $R$-module induces a $\ZZ^d$-graded $\calR$-module via equation~\eqref{87}.

Note that $x'_j$ naturally inherits the $\ZZ^d$-degree $\deg(x_j'):=-\deg(x_j)=\bolda_j$\index{deg@$\deg$}.
Since $x_j\del_j=-x'_j\del_j'$, we put $\deg(\del'_j):=\deg(x_j\del_j)-\deg(x'_j)=-\bolda_j$\index{deg@$\deg$}.
This makes $\calD$ a $\ZZ^d$-graded sheaf of right $\calR$-modules.

Similarly, let $L=(L_x,L_\del)$ be a weight vector on $D=D_\emptyset$, so $L_\del+L_x\geq 0$. 
We extend $L$ to the variables $x'_j$ and $\del'_j$ by setting $L_{x'_j}:=-L_{x_j}$\index{Lx@$L_{x'}$} and $L_{\del'_j}:=L_{x_j}+L_{\del_j}-L_{x'_j}$\index{Ld@$L_{\del'}$}.
This extends the $L$-filtration on $D=D_{(\emptyset)}$ to a \emph{global $L$-filtration} on the sheaf $\calD$\index{L@$L$}.
Let $\calW:=\gr^L(\calD)$\index{W@$\calW$} be the sheaf of associated graded algebras.

By definition, $L_{\del}+L_{x}=L_{\del'}+L_{x'}$.
Suppose $L_x+L_\del>0$, then $W_{(\frakP)}=\gr^L(D_{(\frakP)})$\index{WP@$W_{(\frakP)}$} is the polynomial ring $\CC[x_\frakP,\del_\frakP,x'_{\ol\frakP},\del'_{\ol\frakP}]$ and $\calW$ can be interpreted as the sheaf of regular functions on $T^*\bar X$.

Definitions~\ref{24} and \ref{27} of the \emph{$L$-characteristic variety} $\ch^L$ and \emph{$L$-characteristic cycle} $\ch^L$ generalize to finite $\calD$-modules $\calM$\index{chL@$\ch^L(\calM)$} and homologically holonomic complexes $\calK_\bullet$ of $\calD$-modules\index{ccL@$\ch^L(\calM)$}\index{ccL@$\cc^L(\calK_\bullet)$}.
Both $L$-characteristic variety and cycle behave well under restriction to charts as $L$ is global.
\end{ntn}

Since both $\del_j, x_j\del_j\in D=D_{(\emptyset)}$ extend to global differential operators on $\bar X$, the same is true for the generators $\square_\boldu$ with $A\cdot\boldu=0$ and $E-\beta$ with $1\le i\le n$ in Definition~\ref{81}. 
We will denote these extended operators by the same symbols as the ones we used on $X_\emptyset$.
This provides the motive for introducing a projectivized version of $M_A(\beta)$ as follows.

\begin{dfn}
Define the $\calR$-module $\calI_A:=\calR\cdot I_A$\index{IA@$\calI_A$} as the extension of $I_A$ to $\bar X$.
The \emph{projectivized $A$-hypergeometric system} $\calM_A(\beta)$ is the sheaf of $\calD$-modules 
\[
\calM_A(\beta):=\frac{\calD}{\calD\ideal{I_A,E-\beta}}\index{M_Ab@$\calM_A(\beta)$}.
\]
On the affine patch $X_\frakP$, using notation introduced above, one has
\[
M_{A,(\frakP)}(\beta)=D_{(\frakP)}/D_{(\frakP)}\ideal{I_{A,(\frakP)},E_{(\frakP)}-\beta}\index{M_APb@$\calM_{A,(\frakP)}(\beta)$}.
\]
\end{dfn}

Note, how $I_{A,(\frakP)}=\rho_\frakP(I_A)$\index{IAP@$I_{A,(\frakP)}$} and $E_{(\frakP)}=\rho_\frakP(E)$\index{EP@$E_{(\frakP)}$} behave under change of $\frakP$: 
moving the index $j$ into $\frakP$ results in the change $a_{i,j}x_j\del_j\mapsto-a_{i,j}x'_j\del'_j$ in $E_{(\frakP)}$, and in the substitution $\del_j\mapsto-{x'_j}^2\del'_j$ in $I_{A,(\frakP)}$.

\medskip

We begin our study of $\calM_A(\beta)$ with a basic observation. 
Since the $L$-filtration is compatible with restriction to charts, we have
$(\sigma^L(P))_{(\frakP)}=\sigma^L(P_{(\frakP)})$ for all global sections $P$ of $\calD$. 
Thus the local sections of the ideal sheaf
\[
\calI^L_A:=\gr^L(\calI_A)=\calR\ideal{\sigma^L(\square_\boldu)\mid A\cdot\boldu=0}\index{ILA@$\calI^L_A$}
\]
are given by $I^L_{A,(\frakP)}=\gr^L(I_{A,(\frakP)})$\index{ILAP@$I^L_{A,(\frakP)}$}.
More interestingly, Buchberger's algorithm yields the following analog to Lemma~\ref{7} regarding the $L$-graded ideal of $\calD\cdot I_A$.

\begin{lem}\label{8}
As ideal sheaves in $\calW$, 
\[
\gr^L(\calD\cdot\calI_A)=\calW\cdot\calI^L_A
\]
so that for all $\frakP$ one has
\[
\gr^L(D_{(\frakP)}\cdot I_{A,(\frakP)})=W_{(\frakP)}\cdot I^L_{A,(\frakP)}=W_{(\frakP)}\ideal{\sigma^L(\square_\boldu)_{(\frakP)}\mid\boldu\in\ZZ^n,A\cdot\boldu=0}.\qed
\]
\end{lem}

We now study the variety of $T^*\bar X$ defined by $\calW\cdot\calI^L_{A}$.
Since $\calI^L_A$ is $\ZZ^d$-graded, the associated prime ideals of $\calW\cdot\calI^L_{A}$ will be $\ZZ^d$-graded as well and we hence search for $\ZZ^d$-graded sheaves of prime ideals $\calP\subseteq\calW$ containing $\calW\cdot\calI^L_{A}$.
By Theorem~\ref{13} the $\ZZ^d$-graded prime ideals of $R_{(\emptyset)}$ containing $I^L_A$ are of the form $I^\tau_A$ with $\tau\in\Phi_A^L$.
Since their sections are $\bar X$-global, every $\calP$ as above must contain one of the ideals $\calI_A^\tau:=\calR\cdot I^\tau_A$\index{ItA@$\calI_A^\tau$} with $\tau\in\Phi_A^L$.
Thus we shall look for the components of the variety 
\[
\Upsilon^\tau_A:=\Var(\calW\cdot\calI^\tau_A)\subseteq T^*\bar X\index{YtA@$\Upsilon^\tau_A$}.
\]
The sections $I^\tau_{A,(\frakP)}$\index{ItAP@$I^\tau_{A,(\frakP)}$} of $\calI^\tau_A$ over a chart $X_\frakP$ are the elements of $I^\tau_A$ with $\del_j$ replaced by $-{x'_j}^2\del'_j$ for all $j\in\frakP$.
On $T^*X_\emptyset$, $I^\tau_{A,(\emptyset)}=I^\tau_A$ defines an irreducible variety whose (again irreducible) closure in $T^*\bar X$ we denote by $\Upsilon^{\tau,\emptyset}_A$\index{YteA@$\bar\Upsilon^{\tau,\emptyset}_A$}. 

Assume first that $\frakP$ is contained in $\tau$.
In this case, $W_{(\frakP)}\cdot I^\tau_{A,(\frakP)}$ is generated by $I_{\tau,(\frakP)}$ and $\del_{\ol{\tau\cup\frakP}}$ which reduces the problem to understanding $W_{(\frakP)}\cdot I_{\tau,(\frakP)}$.

\begin{lem}\label{84}
If $\frakP\subseteq\tau$ then the ideal $W_{(\frakP)}\cdot I_{\tau,(\frakP)}$ is a prime ideal.
\end{lem}

\begin{proof}
By equation~\eqref{87}, $I_{\tau,(\frakP)}=\rho_\frakP(I_\tau)$ is
prime in $R_{(\frakP)}$ and hence so is $I_{\tau,(\frakP)}\cdot R_{(\frakP)}[x'_{\frakP},x_{\ol\frakP}]$.
Now $W_{(\frakP)}=R_{(\frakP)}[x'_\frakP,x_{\ol\frakP}][\del'_\frakP]$ and $\del'_j=-\rho_\frakP(\del_j)/{x'_j}^2$ for $j\in\frakP$.
Since for $j\in\frakP$, $x'_j$ is not in $R_{(\frakP)}[x'_\frakP,x_{\ol\frakP}]\cdot I_{\tau,(\frakP)}$ (and, by primeness, not in any associated prime either), it follows that under the ring extension $R_\frakP\hookrightarrow W_\frakP$ the ideal $I_{\tau,(\frakP)}$ remains prime.
\end{proof}

It follows that for $\frakP\subseteq\tau$, $\Upsilon^\tau_A$ is irreducible on $X_\frakP$ and $W_{(\frakP)}\cdot I^\tau_{A,(\frakP)}$ is precisely the defining ideal of $\Upsilon^{\tau,\emptyset}_A\cap T^*X_\frakP$.

We now view $\frakP=(\frakP\cap\tau)\sqcup(\frakP\smallsetminus\tau)$.
By the above, $W_{(\frakP)}\cdot I^\tau_{A,(\frakP)}$ is the ideal of $W_{(\frakP)}$ generated by $I_{\tau,(\frakP)}$, $\del_{\ol{\tau\cup\frakP}}$ and $\{{x'_j}^2\del'_j\mid j\in\frakP\smallsetminus\tau\}$.
By Lemma~\ref{84}, the radical of this ideal has a decomposition into prime ideals
\[
\sqrt{I^\tau_{A,(\frakP)}}
=\bigcap_{\frakT\subseteq\frakP\smallsetminus\tau}W_{(\frakP)}\ideal{I_{\tau,(\frakP)},\del_{\ol{\tau\cup\frakP}},x'_\frakT,\del'_{\frakP\smallsetminus(\frakT\cup\tau)}}.
\]
We are prompted to make the following definitions generalizing Definitions~\ref{17} and \ref{35}.

\begin{dfn}
Let
\[
\bar\Phi_A^L:=\{(\tau,\frakT)\mid\tau\in\Phi_A^L,\,\frakT\subseteq\{1,\ldots,n\},\,\frakT\cap\tau=\emptyset\}\index{PLA@$\bar\Phi_A^L$}
\]
and $\bar\Phi_A^{L,k}:=\{(\tau,\frakT)\in\bar\Phi_A^L\mid\tau\in\Phi_A^{L,k}\}$\index{PLkA@$\bar\Phi_A^{L,k}$}.
For $(\tau,\frakT)\in\bar\Phi_A^L$, put
\[
J_{\tau,\frakT,(\frakP)}:=R_{(\frakP)}\ideal{x'_\frakT,\partial_{\ol{\tau\cup\frakT\cup\frakP}},\del'_{\ol{\tau\cup\frakT}\cap\frakP}}\index{JtTP@$J_{\tau,\frakT,\frakP}$}
\] 
and define the prime ideal sheaf $\calI^{\tau,\frakT}_A\supseteq\calW\cdot I^\tau_A\supseteq\calW\cdot I^L_A$\index{ItTA@$\calI_A^{\tau,\frakT}$} on each chart $X_\frakP$ by
\[
I_{A,(\frakP)}^{\tau,\frakT}:=
\begin{cases}
W_{(\frakP)}\ideal{I_{\tau,(\frakP)},J_{\tau,\frakT,(\frakP)}},&\text{if }\frakP\supseteq\frakT,\\ 
W_{(\frakP)},&\text{otherwise.}
\end{cases}\index{ItTAP@$I_{A,(\frakP)}^{\tau,\frakT}$}
\]
The irreducible varieties $\Upsilon^{\tau,\frakT}_A:=\Var(\calI^{\tau,\frakT}_A)\subseteq T^*\bar X$\index{YtTA@$\Upsilon_A^{\tau,\frakT}$} define quasi-affine varieties
\[
\dot\Upsilon_A^{\tau,\frakT}
:=\Var(I_{A,(\frakT)}^{\tau,\frakT})\smallsetminus\bigcup_{\substack{(\tau',\frakT)\in\bar\Phi_A^L\\\tau'\subseteq\tau}}\Var(I_{A,(\frakT)}^{\tau',\frakT})\subseteq T^*X_\frakT\index{YtTA@$\dot\Upsilon_A^{\tau,\frakT}$}.
\]
Note that $\Upsilon^{\tau,\frakT}_A$ (and hence $\dot\Upsilon_A^{\tau,\frakT}$)  has dimension $\dim(\tau)+n$ and meets the chart $T^*X_\frakP$ exactly if $\frakP\supseteq\frakT$.
Let finally $\nu_A^{L,\tau,\frakT}$\index{nLtTA@$\nu_A^{L,\tau,\frakT}$} be the multiplicity of $\calI^L_A$ along $\Upsilon^{\tau,\frakT}_A$.
\end{dfn}

The quasi-affine varieties $\dot\Upsilon_A^{\tau,\frakT}$ will play the role of the orbits in Section~\ref{spectoric}.
In fact, the action of the $d$-torus $\TT$\index{T@$\TT$} on $T_0^*X$ can be extended to
$\bar X$ via the extension of the coordinate ring $R$ of $T_0^*X$:
Denote
\[
Z_\frakP:=\Var(x'_\frakP)\subseteq X_\frakP\index{ZP@$Z_\frakP$}.
\]
For $p\in Z_\frakP$, $\{\del_j\mid j\notin\frakP\}\cup\{\del'_j\mid j\in \frakP\}$ are coordinates on $T_p^*\bar X$. 
Similar to \eqref{16}, $t\in\TT$ applies componentwise to $\xi=(\xi, \xi')\in T_p^*\bar X$ by $t\cdot\xi_j=t^{\bolda_j}\xi_j$ but $t\cdot \xi'_j=0$ because $\del_j=-{x'_j}^2\del'_j$ vanishes on $Z_\frakP$. 
When restricted to $T_p^*Z_\frakP$, this becomes a group action uniform in $p\in Z_\frakP$. 
To generalize Definition~\ref{17}, let $(\tau,\frakT)\in\bar\Phi_A^L$, denote by $\boldone^{\tau,\frakT}_A$\index{1tT@$\boldone^{\tau,\frakT}_A$} the point $\boldone^\tau_A$ in $T_0^*Z_\frakT$, and set
\[
O_A^{\tau,\frakT}:=\Orb(\boldone^{\tau,\frakT}_A)\subseteq T_0^*Z_\frakT\index{OtTA@$O_A^{\tau,\frakT}$}.
\]
By Lemma~\ref{62}, $O_A^{\tau,\frakT}=\Var(I_{A\smallsetminus\frakT}^\tau)$ so that $O_A^{\tau,\frakT}\times\CC^{|\frakT|}$ identifies with $O_A^\tau$.
The inclusion $Z_\frakT=\Var(x'_\frakT)\subseteq X_\frakT$ induces 
natural maps of cotangent spaces
\[
\xymat{T^*X_\frakT\,&\ar@{_(->}_-{\iota_\frakT}[l]\,\,\,T^*X_{\frakT\vert Z_\frakT}\,\ar@{->>}^-{\pi_\frakT}[r]&\,\,T^*Z_\frakT}\index{iT@$\iota_\frakT
$}\index{pT@$\pi_\frakT$}.
\]
and $\dot\Upsilon_A^{\tau,\frakT}=\iota_\frakT(\pi_\frakT^{-1}(O_A^{\tau,\frakT}\times Z_\frakT))\index{YAtT@$\dot\Upsilon_A^{\tau,\frakT}$}$.
The partial Fourier transform
\begin{equation}\label{83}
\xymat{W_{(\frakT)}\ar[r]&W_{(\emptyset)}},\quad
\xymat{x'_\frakT\ar@{|->}[r]&\del_\frakT},\quad
\xymat{\del'_\frakT\ar@{|->}[r]&-x_\frakT},
\end{equation}
identifies $\dot\Upsilon_A^{\tau,\frakT}$ with the product $O_A^\tau\times\CC^n$

By the preceding discussion, Theorem~\ref{13} and Proposition~\ref{58} can be generalized as follows.

\begin{thm}\label{82}
A decomposition of $\Upsilon^L_A$ into irreducible components is given by
\[
\Upsilon^L_A=\bigcup_{(\tau,\frakT)\in\bar\Phi^{L,d-1}_A}\Upsilon^{\tau,\frakT}_A.
\]
More generally, there is a stratification
\[
\Upsilon_A^L=\bigsqcup_{(\tau,\frakT)\in\bar\Phi_A^L}\dot
\Upsilon_A^{\tau,\frakT}
\]
such that $\dot\Upsilon_A^{\tau',\frakT'}$ is in the closure of $\dot\Upsilon_A^{\tau,\frakT}$ if and only if $\tau'\subseteq\tau$ and $\tau
\smallsetminus\tau'\subseteq\frakT'\supseteq\frakT$.
For $(\tau,\frakT)\in\bar\Phi^{L,d-1}_A$, the multiplicity of $\Upsilon^{\tau,\frakT}_A$ in $\Upsilon^L_A$ is
\[
\nu^{L,\tau,\frakT}_A=2^{|\frakT|}\cdot\nu^{L,\tau}_A=2^{|\frakT|}\cdot[\ZZ^d:\ZZ\tau].
\]
\end{thm}

\begin{proof}
The irreducible decomposition is a consequence of Lemma~\ref{84}.
The existence of the stratification follows from Theorem~\ref{13} and the preceding discussion. 
For $j\in\tau\smallsetminus\frakP$, $\rho_{\frakP\cup\{j\}}(\del_j)=-{x'_j}^2\del'_j$ and the adjacencies in Theorem~\ref{13} imply that
\[
\sqrt{\ideal{I_{A,(\frakT\cup\{j\})}^{\tau,\frakT},x_j'}}
=\bigcap_{\tau'\subseteq\tau\smallsetminus\{j\}}I_{A,(\frakT\cup\{j\})}^{\tau',\frakT\cup\{j\}}
\]
which leads to the adjacencies of the $\dot\Upsilon_A^{\tau,\frakT}$.
The multiplicity $\nu^{A,\tau,\frakT}_A$ can be measured on the chart $X_\frakT$. 
Then the statement is a consequence of the fact that the ideal generated by
${x_j'}^2\del'_j$ has multiplicity two along $x'_j=0$, combined with the formula for $\nu^{L,\tau}_A$ in  Proposition~\ref{58}.
\end{proof}

This ends our discussion of the globalization of the toric ideal $I_A$.
From now on we assume that $L_{x_j}+L_{\del_j}=c>0$ for all $j$ and investigate the projectivized $A$-hypergeometric system $\calM_A(\beta)$.
In that case, $E_{(\frakP)}$ is $L$-homogeneous of $L$-degree $c$ for every chart index $\frakP$. 
We imitate Definition~\ref{74} as follows.

\begin{dfn}
For $(\tau,\frakT)\in\bar\Phi^{\tau,\frakT}_A$, set $C_A^{\tau,\frakT}:=\dot\Upsilon^{\tau,\frakT}_A\cap\Var(E_{(\frakT)})\subseteq T^*X_\frakT$\index{CtTA@$C_A^{\tau,\frakT}$} and let $\calP_{\tau,\frakT}\subseteq\calW$\index{PtT@$\calP_{\tau,\frakT}$} be the prime ideal sheaf of $\bar C_A^{\tau,\frakT}$.
\end{dfn}

The isomorphism of $\dot\Upsilon^{\tau,\frakT}_A$ with $O_A^\tau\times\CC^n$ via the partial Fourier transform~\eqref{83} mentioned above identifies $E_{(\frakT)}$ with $E_{(\emptyset)}$.
Thus $C_A^{\tau,\frakT}$ is a quasi-affine variety isomorphic to $C^\tau_A$.
This yields the following estimate for the $L$-characteristic variety of $\calM_A(\beta)$ in the spirit of Proposition~\ref{55}.

\begin{prp}
The characteristic variety of $\calM_A(\beta)$ is contained in the union of the closures of the Fourier twisted orbit conormals:
\[\pushQED{\qed}
\ch(\calM_A(\beta))\subseteq\bigcup_{(\tau,\frakT)\in\Phi^L_A}\bar C_A^{\tau,\frakT}.\qedhere
\]
\end{prp}

\medskip

We now wish to compute the actual components and their multiplicities in the characteristic cycle of $\calM_A(\beta)$.
We proceed exactly as in Section~\ref{cycle}. 
Multiplicities along a candidate component $\bar C^{\tau,\frakT}_A$ where $(\tau,\frakT)\in\bar\Phi^L_A$ can be computed on the chart $X_\frakT$. 
This means the $\calD$-module $\calM_A(\beta)$ can be replaced by the $D_{(\frakT)}$-module $M_{A,(\frakP)}(\beta)$.
Parallel to the affine case in Definition~\ref{19}, one can introduce a global Euler--Koszul complex on $\bar X$.

\begin{dfn}
Let $N$ be a $\ZZ^d$-graded $R$-module and $\beta\in\CC^d$.
We consider $\calD$ as a $\ZZ^d$-graded right $\calR$-module and $N$ as a $\ZZ^d$-graded $\calR$-module as in Notation $\ref{85}$.
Then $\calD\otimes_\calR N$ is a $\ZZ^d$-graded left $\calD$-module.
On a $\ZZ^d$-graded left $\calD$-module $\calM$, the assignment~\eqref{86} and $\calD$-linear extension yield $d$ commuting $\calD$-linear endomorphisms $E-\beta$.
We define the \emph{global Euler--Koszul complex} $\calK_{A,\bullet}(N;E-\beta)$\index{KAiNEb@$\calK_{A,i}(N;E-\beta)$} and the \emph{global Euler--Koszul homology} $\calH_{A,\bullet}(N;E-\beta)$\index{HAiNEb@$\calH_{A,i}(N;E-\beta)$} as the sheafified version of the objects in Definition~\ref{19}.
Note that $\calH_{A,0}(S_A;E-\beta)=\calM_A(\beta)$.
\end{dfn}

A good behavior of this global Euler--Koszul functor such as a long exact sequence of Euler--Koszul homology requires flatness of $\calD$ over $\calR$.

\begin{lem}\label{68}
The sheaves of $\calR$-modules $\calD$ and $\calW$ are locally free and hence flat.
\end{lem}

\begin{proof}
It suffices to consider the single variable case on a chart $X_\frakP$.
For $j\in\frakP$, $\CC[x'_j,\del'_j]$ is a free $\CC[-{x'_j}^2\del'_j]$-module with basis ${x'_j}^k{\del'_j}^l$ where $k\le1$ or $l=0$. 
Since $\CC[x'_j,\del'_j]=\gr^F(\CC[x'_j]\ideal{\del'_j})$ the same is true for $\CC[x'_j]\ideal{\del'_j}$.
\end{proof}

The Euler--Koszul homology sheaves of a toric $R$-module are holonomic for reasons similar to those given in \cite{MMW05} for Euler--Koszul homology modules.
Their $L$-characteristic variety is hence $n$-dimensional by Proposition~\ref{55} and \cite{Smi01}.
A spectral sequence argument like in Theorem~\ref{50} shows the possible components of $\ch^L(\calH_{A,\bullet}(N;\beta))$ are of the form $\bar C^{\tau,\frakT}_A$ with $(\tau,\frakT)\in\bar\Phi_A^L$.

\begin{dfn}
For $(\tau,\frakT)\in\bar\Phi_A^L$ and $i\in\NN$, we denote by $\mu^{L,\tau,\frakT}_{A,i}(\beta)$\index{mLtTAib@$\mu^{L,\tau,\frakT}_{A,i}(\beta)$} the multiplicity of $\gr^L(\calH_{A,i}(S_A;\beta))$ along $\bar C_A^{\tau,\frakT}$ and put $\mu^{L,\tau,\frakT}_A:=\sum_{i=1}^d(-1)^i\mu^{L,\tau,\frakT}_{A,i}(\beta)$.
\end{dfn}

By looking at the chart $X_\frakT$ it follows from the appendix that, as in the affine case, the alternating sum of the $\mu^{L,\tau,\frakT}_{A,i}(\beta)$ agrees with the intersection multiplicity of $I^L_{A,(\frakT)}$ with $E_{(\frakT)}$ in the local ring at $\bar C_A^{\tau,\frakT}$.
By Theorem~\ref{82}, this reduces to the intersection multiplicity of $I^{\tau,\frakT}_{A,(\frakT)}$ with $E_{(\frakT)}$, multiplied by $\nu^{L,\tau,\frakT}_A=2^{|\frakT|}\nu^{L,\tau}_A$.
From Theorem~\ref{50} we conclude the following statement which shows in particular that $\mu^{L,\tau,\frakT}_A$ is independent of $\beta$ as suggested by the notation.
For the sake of simplicity we concentrate on the case of interest $N=S_A$.

\begin{thm}
For all $(\tau,\frakT)\in\bar\Phi_A^L$ and all $\beta\in\CC^d$, 
\[
\mu^{L,\tau,\frakT}_A
=\chi^C(\calW/\calW\ideal{E},\calW\otimes_\calR S_A^L)
=2^{|\frakT|}\mu^{L,\tau}_A,
\]
where $C=\bar C_A^{\tau,\frakT}$, is independent of $\beta$ and positive.\qed
\end{thm}

In \cite[Thm.~6.3]{MMW05}, a spectral sequence is constructed that is used to show good vanishing properties of Euler--Koszul homology. 
We leave it to the reader to verify that it generalizes to our global context.
It follows that if $\beta$ is generic then $\calK_{A,\bullet}(S_A,E-\beta)$ has only homology in degree zero.
In particular, $\mu^{L,\tau,\frakT}_{A,0}(\beta)=\mu^{L,\tau,\frakT}_{A}>0$ for generic $\beta$, similar to Corollary~\ref{33}.

In order to show nonvanishing for all $\beta$ one can generalize Theorems~\ref{10} and \ref{47} by following exactly the steps of our proof in the chart $X_\frakT$ corresponding to $\bar C_A^{\tau,\frakT}$, always twisting the toric input module by $\rho_\frakP$ from equation~\eqref{87}.

\begin{thm}\label{76}\pushQED{\qed}
The $L$-characteristic variety of $\calM_A(\beta)$ is given by
\[
\ch^L(\calM_A(\beta))=\bigcup_{(\tau,\frakT)\in\bar\Phi^L_A}\bar C_A^{\tau,\frakT}.
\]
For all $(\tau,\frakT)\in\bar\Phi_A^{L}$ and all $\beta\in\CC^d$,
\[
\mu^{L,\tau,\frakT}_{A,0}(\beta)\geq
\mu^{L,\tau,\frakT}_A=2^{|\frakT|}\mu^{L,\tau}_A
\]
is positive, and equality holds for generic $\beta$.
In particular, the $L$-characteristic variety of $\calM_A(\beta)$ is entirely determined by the $(A,L)$-umbrella $\Phi^L_A$.
For $\tau\in\Phi_A^{L,n-1}$,
\[
\mu^{L,\tau,\frakT}_A=2^{|\frakT|}\nu^{L,\tau}_A=2^{|\frakT|}[\ZZ^d:\ZZ\tau].\qedhere
\]
\end{thm}

\medskip

In order to study slopes of $\calM_{A}(\beta)$ along coordinate varieties, we consider $\frakV$ as a partitioned set $\frakV=\frakV_0\sqcup\frakV_\infty\subseteq\{1,\ldots,n\}$\index{V@$\frakV$}\index{V@$\frakV_0$}\index{V@$\frakV_\infty$}.
As in \eqref{94}, $\frakV$ defines the coordinate variety
\[
Y:=\Var(x_{\frakV_0},x'_{\frakV_\infty})\subseteq\bar X\index{Y@$Y$}
\]
which meets the chart $X_\frakP$ exactly if $\frakV_0\cap\frakP=
\emptyset$ and $\frakV_\infty\subseteq\frakP$.
On any such $X_\frakP$, $Y$ induces the family of filtrations
\[
L=pF+qV,\quad p/q\in\QQ_{>0}\cup\{\infty\}\index{L@$L$},
\]
on $D_\frakP$ as in \eqref{73} where $V$ is the $V$-filtration along $Y\cap X_
\frakP$ defined by the assignment $-V_{x_i}=1=V_{\del_i}$ for $i\in
\frakV_0$, $V_{x'_i}=1=-V_{\del'_i}$ for $i\in\frakV_\infty$, and $V_{x'_i}=0=-V_{\del'_i}$ for $i\not\in\frakV$\index{V@$V$}\index{Vx@$V_{x'}$}\index{Vd@$V_{d'}$}.
One may view the patch $W_\frakP=\gr^L(D_\frakP)$ of $\calW$ as the ring of polynomial functions on $T^*(T^*_YX_\frakP)$ and identify $T^*T^*_Y\bar X=T^*(\bar X\smallsetminus\Var(x'_{\frakV_0},x_{\frakV_\infty}))$.
By Definition \ref{91}, $L=p/q$ a slope of a finite $\calD$-module $\calM$ on $\bar X$ at $y\in Y$ along $Y$ if the set of components of $\ch^L(\calM)$ which meet $T^*_y\bar X$ jumps at $L$. 

Recall that all slopes of $M_A(\beta)$ along coordinate varieties are slopes at the origin in $X$ by Lemma \ref{78}.
But the following example shows that jumps of $\ch^L(\calM_A(\beta))$ can be irrelevant for the slopes along $Y$. 

\begin{exa}
Let $A=\begin{pmatrix}3&1&0\\0&1&3\end{pmatrix}$, pick $Y=\Var(x_1',x_2')$ and $L=pF+qV$ with $p>0,q\geq 0$ in $\QQ$.
Since $I^L_A$ is Cohen--Macaulay, $\gr^L(\calW(I_A,E-\beta))$ equals $\ideal{\del_2^3,3x_1\del_1+x_2\del_2,x_2\del_2+3x_3\del_3}$ if $p/q>2$ but $\ideal{\del_1\del_3,3x_1\del_1+x_2\del_2,x_2\del_2+3x_3\del_3}$ if $p/q<2$.
It follows that $\ch^L(\calM_A(\beta))$ is defined by $\ideal{\del_2,x_1\del_1,x_3\del_3}$ for $p/q>2$ and by $\ideal{x_2\del_2,\del_1\del_3,x_3\del_3,x_1\del_1}$ if $p/q<2$.
Thus, $2$ is the only possible slope value and, as $p/q$ passes through $2$ from above, the component of $\ch^L(\calM_A(\beta))$ to $\ideal{x_1,\del_2,x_3}$ is replaced by the components to $\ideal{x_2,\del_1\del_3,x_3\del_3,x_1\del_1}$.
However, as $Y$ is defined by $x_1=x_2=\infty$, none of these components meet $Y$ and hence $\calM_A(\beta)$ is regular along $Y$.
\end{exa}

It follows that Theorem~\ref{76} does not give directly a combinatorial description of the slopes of $\calM_A(\beta)$ along $Y$.
We must select the components of $\ch^L(\calM_A(\beta))$ whose projection meets $Y$.
For this, we look at a chart $X_\frakP$ which meets $Y$ and hence $\frakP\supseteq\frakV_\infty$ and $\frakP\cap\frakV_0=\emptyset$.
Consider the component $\bar C^{\tau,\frakT}_A$ of $\ch^L(\calM_A(\beta))$ for some $(\tau,\frakT)\in\bar\Phi^L_A$.
Since $I^{\tau,\frakT}_{A,(\frakP)}$ contains $x'_\frakT$, $\bar C^{\tau,\frakT}_A$ can be visible in $X_\frakP$ only if $\frakP\supseteq\frakT$. 
On the other hand, if $\bar C^{\tau,\frakT}_A$ is not visible in $X_\frakP$ then (since $\bar C^{\tau,\frakT}_A$ is irreducible) its defining ideal contains an $x_i$ with $i \in\frakP$.
This may be the case even if $\frakP\supseteq\frakT$.

\begin{dfn}
$\tau\in\Phi^L_A$ is a \emph{pyramid} with vertex $i\in\tau$ if $\dim(\tau\setminus\{i\})<\dim(\tau)$.
\end{dfn}

\begin{lem}
Let $(\tau,\frakT)\in\bar\Phi^L_A$.
There is no chart $X_\frakP$ that meets both $Y$ and $\bar C^{\tau,\frakT}_A$ if and only if either $\frakT\cap\frakV_0\ne\emptyset$ or $\tau$ is a pyramid with vertex in $\frakV_\infty$.
Moreover, if $\bar C^{\tau,\frakT}_A$ is a component of $\cc^{L'}(\calM_A(\beta))$ for $L'$ near $L$ and meets some chart $X_\frakP$ which meets $Y$ then its projection also meets $Y$.
\end{lem}

\begin{proof}
To prove the first statement, assume that $\frakT\cap\frakV_0=\emptyset$ and $\bar C^{\tau,\frakT}_A$ does not meet the chart $X_\frakP$ with $\frakP=\frakV_\infty\cup\frakT$.
Then $C^{\tau,\frakT}_A\subseteq\Var(x_i)$ for some $i\in\frakV_\infty$ by the preceding discussion.
Since $C^{\tau,\frakT}_A$ is the conormal to $O^{\tau,\frakT}_A$, the ideal $I_A^{\tau,\frakT}$ is hence independent of of $\del_i$.
By definition of $I_A^{\tau,\frakT}$ through the toric ideal $I_\tau$, this is equivalent to $\tau$ being a pyramid with vertex $i$.
The converse implication of the first statement is obvious.

Under the hypothesis of the second statement, the ideal $\calP_{\tau,\frakT}$ of $\bar C^{\tau,\frakT}_A$ is an associated prime of an $L$- and an $L'$-homogeneous ideal for some $L'$ near $L$.
As such, $\calP_{\tau,\frakT}$ is $(L,L')$- and hence $(F,V)$-bihomogeneous.
It follows that the nontrivial $W_{(\frakP)}$-ideal $P'=P_{\tau,\frakT,(\frakP)}+\ideal{\del_{\ol\frakP},\del'_\frakP}$ is $(F,V)$-bihomogeneous, generated by $\ideal{\del_{\ol\frakP},\del'_\frakP}$ and some $\del$-independent $V$-homogeneous elements.
In particular, $P'+\ideal{x_{\frakV_0},x'_{\frakV_\infty}}$ is nontrivial and so $\bar C^{\tau,\frakT}_A$ meets $Y$.
\end{proof}

From the dimension one case, it is clear that non-$(F,V)$-bihomogeneous components are relevant for slopes.
Indeed, the $L$-characteristic variety equals $\Var(x_1\del_1)$ for most $L$-homogeneous differential operators in one variable $x_1$.
However we believe that the non-$(F,V)$-bihomogeneous components are irrelevant in our case.
More precisely, it seems that the conormal space of an $L$-homogeneous non-$(F,V)$-bihomogeneous pointed toric ring can not meet any point $T^*_y\bar X$ with $y\in Y$ unless $\frakV_\infty=\emptyset$.  
This would imply the following combinatorial description of slopes of projectivized $A$-hypergeometric systems along coordinate varieties.

\begin{dfn}
Let $L=pF+qV$ with $p/q\in\QQ_{>0}\cup\{\infty\}$ where $V$ is the $V$-filtration along $Y=\Var(x_{\frakV_0},x'_{\frakV_\infty})$.
By $\Phi^Y_A(p/q)$\index{PYA@$\Phi^Y_A(p/q)$} we denote the subset of $\tau\in\Phi^L_A$ which are not a pyramid with vertex in $\frakV_\infty$.
\end{dfn}

\begin{cnj}\label{93}
The slopes of the projectivized $A$-hypergeometric system $\calM_A(\beta)$ along a coordinate variety $Y$ are the jump parameters $p/q\in\QQ_{>0}\cup\{\infty\}$ of $\Phi^Y_A(p/q)$.
\end{cnj}
	
\end{arxiv}

\appendix
\numberwithin{equation}{section}

\section{On spectral sequences}

Let $(K^\bullet,d)$ be a complex of groups equipped with a descending filtration $F$ of subgroups such that $d(F^pK^{\bullet})\subseteq F^pK^{\bullet +1}$.
We stress that there are no boundedness conditions on the filtration.
We abbreviate $\gr=\gr_F$ and freely use induced filtrations on subgroups and quotient groups.

By \cite[Ch.~I, Thm.~4.2.1, Thm.~4.2.2]{God58} there is an associated spectral sequence involving the following data:
\begin{align}\label{45}
\nonumber Z^{p,q}_r&=F^pK^{p+q}\cap d^{-1}(F^{p+r}K^{p+q+1}),\\
\nonumber B^{p,q}_r&=F^pK^{p+q}\cap d(F^{p-r}K^{p+q-1}),\\
E^{p,q}_r&=\gr^p\frac{d^{-1}(F^{p+r}K^{p+q+1})}{d(F^{p-r+1}K^{p+q-1})}
=\gr^p\frac{Z^{p,q}_r}{B^{p,q}_{r-1}}
=Z^{p,q}_r/\left(Z^{p+1,q-1}_{r-1}+B^{p,q}_{r-1}\right).
\end{align}
Note that $Z^{p,q}_r\cap F^{p+1}K^{p+q}=Z^{p+1,q-1}_{r-1}$, $d(Z^{p-r,q+r-1}_r)=B^{p,q}_r$, $E^{p,q}_0=\gr^pK^{p+q}$, and $E^{p,q}_1=H^{p+q}(\gr^p(K^\bullet),\gr^p(d))$.

For all $p,q$ and $r$, $d$ induces a group homomorphism
$d^{p,q}_r:E^{p,q}_r\to E^{p+r,q-r+1}_r$ such that $d^{p,q}_r\circ
d^{p-r,q+r+1}_r=0$.  More explicitly,
\begin{equation}\label{42}
d^{p,q}_r\left(z\modulo (Z^{p+1,q-1}_{r-1}+B^{p,q}_{r-1})\right)
=d(z)\modulo (Z^{p+1+r,q-r}_{r-1}+B^{p+r,q-r+1}_{r-1})
\end{equation}
and $E^{p,q}_r=(\ker d_r^{p,q})/(\image d_r^{p-r,q+r-1})$.
Let 
\begin{align}\label{46}
\nonumber Z^{p,q}_\infty&=F^pK^{p+q}\cap d^{-1}(0),\\
\nonumber B^{p,q}_\infty&=F^pK^{p+1}\cap d(K^{p+q-1}),\\
E^{p,q}_\infty&=\gr^pH^{p+q}(K^\bullet)
=\gr^p\frac{Z^{p,q}_\infty}{B^{p,q}_\infty}
=Z^{p,q}_\infty/(Z^{p+1,q-1}_\infty+B^{p,q}_\infty).
\end{align}

In general, there is no relation between the terms $E^{p,q}_r$ for finite $r$ and $E^{p,q}_\infty$. The purpose of this section is to prove the following theorem and its corollary below.

\begin{thm}\label{40}
Let $D$ be a ring and $(K^\bullet,d)$ be a complex of $D$-modules, both equipped with a descending filtration $F$ of subgroups subject to the hypotheses:
\begin{enumerate}[(a)]
\item\label{40a} $D$ is a filtered ring: $(F^pD)\cdot(F^{p'}D)\subseteq F^{p+p'}D$.
\item\label{40c} $K^\bullet$ is a filtered $D$-module: $(F^pD)\cdot(F^{p'}K^n)\subseteq F^{p+p'}K^n$.
\item\label{40b} $d$ is compatible with $F$: $d(F^pK^{p+q})\subseteq F^p K^{p+q+1}$.
\item\label{40e} $E_0=\gr(K^\bullet)$ is a Noetherian $\gr(D)$-module.
\item\label{40d} $F$ is exhaustive and separated on $K^\bullet$: $\bigcup_pF^pK^n=K^n$, $\bigcap_pF^pK^n=0$.
\end{enumerate}
Then one has the consequences:
\begin{enumerate}[(1)]
\item\label{40f} $d_r\colon E_r\to E_r$ are graded $\gr(D)$-linear maps of graded $\gr(D)$-modules with $\gr^{p'}(D)\cdot E^{p,q}_r\subseteq E^{p+p',q-p'}_r$.
\item\label{40g} $E_r$ is Noetherian and $d_r$ is the zero map for all sufficiently large $r$.
\item\label{40h} The stable value of $E_r$ is isomorphic to $E_\infty$.
\end{enumerate}
\end{thm}

\begin{proof}
We abbreviate $W=\gr(D)$.
\begin{asparaenum}[(1)]

\item follows from \eqref{40a}, \eqref{40c}, and \eqref{40b}:
$E_0$ is a graded $W$-module by \eqref{40a} and \eqref{40c}.
By induction, we may assume that $E_r$ is a graded $W$-module and $d_{r-1}$ is graded $W$-linear.
Let $\delta\in F^{p'}D$ and $c\in Z^{p,q}_r$ represent non zero elements $\bar\delta\in W^{p'}$ and $\bar c\in E^{p,q}_r$. 
By \eqref{42} and \eqref{40b},
\begin{align*}
d_r(\bar\delta\cdot\bar c)
&=d_r\left((\delta\modulo F^{p'+1}D)\cdot(c\modulo (Z^{p+1,q-1}_{r-1}+B^{p,q}_{r-1}))\right)\\
&=d_r\left(\delta c\modulo (Z^{p+p'+1,q-p'-1}_{r-1}+B^{p+p',q-p'}_{r-1})\right)\\
&=d(\delta c)\modulo (Z^{p+p'+r+1,q-p'-r-1}_{r-1}+B^{p+p'+r,q-r-p'}_{r-1})\\
&=\delta d(c)\modulo (Z^{p+p'+r+1,q-p'-r-1}_{r-1}+B^{p+p'+r,q-r-p'}_{r-1})\\
&=(\delta\modulo F^{p'+1}D)\cdot\left(d(c)\modulo (Z^{p+r+1,q-r-1}_{r-1}+B^{p+r,q-r}_{r-1})\right)\\
&=(\delta\modulo F^{p'+1}D)\cdot d_r\left(c\modulo (Z^{p+1,q-1}_{r-1}+B^{p,q}_{r-1})\right)\\
&=\bar\delta\cdot d_r(\bar c)\in E^{p+p',q-p'}_r.
\end{align*}
Thus, $d_r$ is graded $W$-linear and $\ker(d_r)$, $\image(d_r)$, and $E_r=\ker(d_r)/\image(d_r)$ are graded $W$-modules.

\item follows from \eqref{40e}: 
As a subquotient of $E_r$, $E_{r+1}$ is Noetherian by induction. 
Let $\pi_r\colon \ker(d_r)\to\ker(d_r)/\image(d_r)=E_{r+1}$ denote the natural projection. 
Note that $\pi_r^{-1}(\image(d_{r+1}))\supseteq\pi_r^{-1}(0)=\image(d_r)$. 
Then the chain
\[
\image(d_0)\subseteq
\pi_0^{-1}(\image(d_1))\subseteq
\pi_0^{-1}(\pi_1^{-1}\image(d_2))\subseteq
\cdots
\]
of submodules of $E_0$ must stabilize at some $r_0$:
$\pi_{r}^{-1}(\image (d_{r+1}))=\image(d_{r})=\pi_r^{-1}(0)$ and hence $\image(d_{r+1})=0$ for all $r\ge r_0$.

\item follows from \eqref{40d}:
\begin{align*}
B^{p,q}_\infty&=F^pK^{p+q}\cap d(K^{p+q-1})
=F^pK^{p+q}\cap d\left(\bigcup_rF^{p-r}K^{p+q-1}\right)\\
&=\bigcup_r(F^pK^{p+q}\cap d(F^{p-r}K^{p+q-1}))
=\bigcup_r B^{p,q}_r,\\
Z^{p,q}_\infty&=F^pK^{p+q}\cap d^{-1}(0)
=F^pK^{p+q}\cap d^{-1}\left(\bigcap_rF^{p+r}K^{p+q+1}\right)\\
&=\bigcap_r(F^pK^{p+q}\cap d^{-1}(F^{p+r}K^{p+q+1}))
=\bigcap_rZ^{p,q}_r.
\end{align*}
By \eqref{42}, $d_r=0$ means that
\begin{equation}\label{43}
B^{p,q}_r=d(Z^{p-r,q+r-1}_r)
\subseteq Z^{p+1,q-1}_{r-1}+B^{p,q}_{r-1}
\subseteq B^{p,q}_{r-1}+F^{p+1}K^{p+q}\\
\end{equation}
and also
\[
d(Z^{p,q}_r)\subseteq Z^{p+r+1,q-r}_{r-1}+B^{p+r,q-r+1}_{r-1}=Z^{p+r+1,q-r}_{r-1}+d(Z^{p+1,q-1}_{r-1})
\]
which implies that
\begin{equation}\label{44}
Z^{p,q}_r
\subseteq Z^{p,q}_{r+1}+Z^{p+1,q-1}_{r-1}
\subseteq Z^{p,q}_{r+1}+F^{p+1}K^{p+q}.
\end{equation}
From \eqref{43} and \eqref{44} follows $\gr^p B^{p,q}_r=\gr^p B^{p,q}_{r-1}$ and $\gr^p Z^{p,q}_r=\gr^p Z^{p,q}_{r+1}$ for $r\ge r_0$.
By \eqref{45} and \eqref{46} and exactness of $\gr$ for induced filtrations, we conclude that, for $r\ge r_0$,
\[
E^{p,q}_r
=\gr^p\frac{Z^{p,q}_r}{B^{p,q}_{r-1}}
=\frac{\gr^pZ^{p,q}_r}{\gr^pB^{p,q}_{r-1}}
=\frac{\gr^p\left(\bigcap_rZ^{p,q}_r\right)}{\gr^p\left(\bigcup_rB^{p,q}_{r-1}\right)}
=\frac{\gr^pZ^{p,q}_\infty}{\gr^pB^{p,q}_\infty}
=\gr^p\frac{Z^{p,q}_\infty}{B^{p,q}_\infty}
=E^{p,q}_\infty.
\]
\end{asparaenum}
\end{proof}

\begin{cor}\label{41}
Under the hypotheses of Theorem~\ref{40} assume that $W=\gr(D)$ is commutative, let $P\in\Spec(W)$ be a prime ideal and denote by $W_P$
the local ring of $\Spec(W)$ at $P$.
Then there is a spectral sequence 
\[
H^p(\gr(K^\bullet)\otimes_W W_P) = {'E}_1\Longrightarrow {'E}_\infty
=(\gr(H^p(K^\bullet)))\otimes_W W_P.
\]
This spectral sequence collapses eventually and converges:
${'E}_r={'E}_{r+1}={'E}_\infty$ for all large $r$. 
\end{cor}

\begin{proof}
The spectral sequence $(E_r,d_r)$ of Theorem~\ref{40} is a sequence of $W$-modules and $W$-morphisms.
Since localization commutes with the formation of kernels and cokernels, the objects and morphisms $({'E}_r,'d_r)=(E_r\otimes_W W_P, d_r\otimes_W W_P)$ form a spectral sequence.
Since $(E_r,d_r)$ collapses, so does $({'E}_r,'d_r)$, converging to the localization of $E_\infty$.
\end{proof}

\begin{rmk}
As far as we can see, the collapse of the sequence is essential for the corollary.
If in general $(E_r,d_r)$ is a spectral sequence of $W$-modules and $W$-morphisms converging to $E_\infty$ in the sense that $E_\infty=\bigcap_r Z_r/\bigcup_r B_r$ then while $(\bigcup_rB_r)\otimes_W W_P= \bigcup_r (B_r\otimes_W W_P)$, it is not necessarily the case that $(\bigcap_r Z_r)\otimes_W W_P= \bigcap_r(Z_r\otimes_W W_P)$.
The problem is that infinite intersections (in contrast to infinite unions) and homomorphisms (such as the localization morphism) are not necessarily interchangeable.

In the situation of the theorem it is probably not sufficient to replace \eqref{40e} by the condition that $K^\bullet$ consists of Noetherian $D$-modules.
For example, in \cite{ST04} some rings appear in the theory of hypergeometric systems that are Noetherian, but the order filtration leads to non-Noetherian associated graded rings.
\end{rmk}

\section*{Acknowledgments}

UW is grateful to Francisco Castro-Jim\'enez, Alicia Dickenstein, Ezra Miller, Anurag Singh and Will Traves for interesting and helpful conversations. 
MS wishes to thank Purdue University for their hospitality during much of the preparation of this article. 
\begin{arxiv}We thank Mar\'{\i}a-Cruz Fern\'andez-Fern\'andez for pointing out an error in the initial version of Lemma~\ref{7}.\end{arxiv}
\begin{duke}We are grateful to the referees for a very careful reading and for valuable suggestions.
We also thank Mar\'{\i}a-Cruz Fern\'andez-Fern\'andez for pointing out an error in the initial version of Lemma~\ref{7}.\end{duke}

\printindex

\bibliographystyle{amsalpha}
\bibliography{sgkz}

\def\cprime{$'$}
\providecommand{\bysame}{\leavevmode\hbox to3em{\hrulefill}\thinspace}
\providecommand{\MR}{\relax\ifhmode\unskip\space\fi MR }
\providecommand{\MRhref}[2]{%
  \href{http://www.ams.org/mathscinet-getitem?mr=#1}{#2}
}
\providecommand{\href}[2]{#2}
\begin{thebibliography}{MMW05}

\bibitem[ACJG96]{ACG96}
A.~Assi, F.~J. Castro-Jim{\'e}nez, and J.~M. Granger, \emph{How to calculate
  the slopes of a {$\mathcal D$}-module}, Compositio Math. \textbf{104} (1996),
  no.~2, 107--123. \MR{MR1421395 (98i:32010)}

\bibitem[ACJG00]{ACG00}
A.~Assi, F.~J. Castro-Jim{\'e}nez, and M.~Granger, \emph{The {G}r\"obner fan of
  an {$A\sb n$}-module}, J. Pure Appl. Algebra \textbf{150} (2000), no.~1,
  27--39. \MR{MR1762918 (2001j:16036)}

\bibitem[Ado94]{Ado94}
Alan Adolphson, \emph{Hypergeometric functions and rings generated by
  monomials}, Duke Math. J. \textbf{73} (1994), no.~2, 269--290. \MR{96c:33020}

\bibitem[BH93]{BH93}
Winfried Bruns and J{\"u}rgen Herzog, \emph{Cohen-{M}acaulay rings}, Cambridge
  Studies in Advanced Mathematics, vol.~39, Cambridge University Press,
  Cambridge, 1993. \MR{MR1251956 (95h:13020)}

\bibitem[BS98]{BS98}
M.~P. Brodmann and R.~Y. Sharp, \emph{Local cohomology: an algebraic
  introduction with geometric applications}, Cambridge Studies in Advanced
  Mathematics, vol.~60, Cambridge University Press, Cambridge, 1998.
  \MR{MR1613627 (99h:13020)}

\bibitem[CDS01]{CDS01}
Eduardo Cattani, Alicia Dickenstein, and Bernd Sturmfels, \emph{Rational
  hypergeometric functions}, Compositio Math. \textbf{128} (2001), no.~2,
  217--239. \MR{MR1850183 (2003f:33016)}

\bibitem[CJT03]{CT03}
Francisco~Jes{\'u}s Castro-Jim{\'e}nez and Nobuki Takayama, \emph{Singularities
  of the hypergeometric system associated with a monomial curve}, Trans. Amer.
  Math. Soc. \textbf{355} (2003), no.~9, 3761--3775 (electronic). \MR{MR1990172
  (2004j:32009)}

\bibitem[CK99]{CK99}
David~A. Cox and Sheldon Katz, \emph{Mirror symmetry and algebraic geometry},
  Mathematical Surveys and Monographs, vol.~68, American Mathematical Society,
  Providence, RI, 1999. \MR{MR1677117 (2000d:14048)}

\bibitem[Ful98]{Ful98}
William Fulton, \emph{Intersection theory}, second ed., Ergebnisse der
  Mathematik und ihrer Grenzgebiete. 3. Folge. A Series of Modern Surveys in
  Mathematics [Results in Mathematics and Related Areas. 3rd Series. A Series
  of Modern Surveys in Mathematics], vol.~2, Springer-Verlag, Berlin, 1998.
  \MR{MR1644323 (99d:14003)}

\bibitem[GGZ87]{GGZ87}
I.~M. Gel{\cprime}fand, M.~I. Graev, and A.~V. Zelevinski{\u\i},
  \emph{Holonomic systems of equations and series of hypergeometric type},
  Dokl. Akad. Nauk SSSR \textbf{295} (1987), no.~1, 14--19. \MR{MR902936
  (88j:58118)}

\bibitem[God58]{God58}
Roger Godement, \emph{Topologie alg\'ebrique et th\'eorie des faisceaux},
  Actualit'es Sci. Ind. No. 1252. Publ. Math. Univ. Strasbourg. No. 13,
  Hermann, Paris, 1958. \MR{MR0102797 (21 \#1583)}

\bibitem[GZK89]{GKZ89}
I.~M. Gel{\cprime}fand, A.~V. Zelevinski\u{\i}, and M.~M. Kapranov,
  \emph{Hypergeometric functions and toric varieties}, Funktsional. Anal. i
  Prilozhen. \textbf{23} (1989), no.~2, 12--26. \MR{90m:22025}

\bibitem[HH03]{Har03}
Mar{\'{\i}}a~Isabel Hartillo~Hermoso, \emph{Slopes of hypergeometric systems of
  codimension one}, Proceedings of the International Conference on Algebraic
  Geometry and Singularities (Spanish) (Sevilla, 2001), vol.~19, 2003,
  pp.~455--466. \MR{MR2023195 (2005d:16039)}

\bibitem[HH05]{Har05}
Mar{\'{\i}}a~Isabel Hartillo-Hermoso, \emph{Irregular hypergeometric systems
  associated with a singular monomial curve}, Trans. Amer. Math. Soc.
  \textbf{357} (2005), no.~11, 4633--4646 (electronic). \MR{MR2156724}

\bibitem[HK84]{HK84}
R.~Hotta and M.~Kashiwara, \emph{The invariant holonomic system on a semisimple
  {L}ie algebra}, Invent. Math. \textbf{75} (1984), no.~2, 327--358.
  \MR{MR732550 (87i:22041)}

\bibitem[Hot98]{Hot98}
R.~Hotta, \emph{Equivariant $d$-modules}, arXiv.org (1998), no.~RT/9805021.

\bibitem[Inc44]{Inc44}
E.~L. Ince, \emph{Ordinary {D}ifferential {E}quations}, Dover Publications, New
  York, 1944. \MR{MR0010757 (6,65f)}

\bibitem[Kou76]{Kou76}
A.~G. Kouchnirenko, \emph{Poly\`edres de {N}ewton et nombres de {M}ilnor},
  Invent. Math. \textbf{32} (1976), no.~1, 1--31. \MR{MR0419433 (54 \#7454)}

\bibitem[Lau85]{Lau85}
Yves Laurent, \emph{Th\'eorie de la deuxi\`eme microlocalisation dans le
  domaine complexe}, Progress in Mathematics, vol.~53, Birkh\"auser Boston
  Inc., Boston, MA, 1985. \MR{MR776973 (86k:58113)}

\bibitem[Lau87]{Lau87}
\bysame, \emph{Polyg\^one de {N}ewton et {$b$}-fonctions pour les modules
  microdiff\'erentiels}, Ann. Sci. \'Ecole Norm. Sup. (4) \textbf{20} (1987),
  no.~3, 391--441. \MR{MR925721 (89k:58282)}

\bibitem[LM99]{LM99}
Yves Laurent and Zoghman Mebkhout, \emph{Pentes alg\'ebriques et pentes
  analytiques d'un {$\mathcal D$}-module}, Ann. Sci. \'Ecole Norm. Sup. (4)
  \textbf{32} (1999), no.~1, 39--69. \MR{MR1670595 (2001b:32015)}

\bibitem[Mal74]{Mal74}
Bernard Malgrange, \emph{Sur les points singuliers des \'equations
  diff\'erentielles}, Enseignement Math. (2) \textbf{20} (1974), 147--176.
  \MR{MR0368074 (51 \#4316)}

\bibitem[Meb89]{Meb89}
Zoghman Mebkhout, \emph{Le th\'eor\`eme de comparaison entre cohomologies de de
  {R}ham d'une vari\'et\'e alg\'ebrique complexe et le th\'eor\`eme d'existence
  de {R}iemann}, Inst. Hautes \'Etudes Sci. Publ. Math. (1989), no.~69, 47--89.
  \MR{MR1019961 (91a:32016)}

\bibitem[MMW05]{MMW05}
Laura~Felicia Matusevich, Ezra Miller, and Uli Walther, \emph{Homological
  methods for hypergeometric families}, J. Amer. Math. Soc. \textbf{18} (2005),
  no.~4, 919--941 (electronic). \MR{MR2163866}

\bibitem[MS05]{MS05}
Ezra Miller and Bernd Sturmfels, \emph{Combinatorial commutative algebra},
  Graduate Texts in Mathematics, vol. 227, Springer-Verlag, New York, 2005.
  \MR{MR2110098 (2006d:13001)}

\bibitem[OT01]{OT01}
Peter Orlik and Hiroaki Terao, \emph{Arrangements and hypergeometric
  integrals}, MSJ Memoirs, vol.~9, Mathematical Society of Japan, Tokyo, 2001.
  \MR{MR1814008 (2003a:32048)}

\bibitem[Sch85]{Sch85}
Pierre Schapira, \emph{Microdifferential systems in the complex domain},
  Grundlehren der Mathematischen Wissenschaften [Fundamental Principles of
  Mathematical Sciences], vol. 269, Springer-Verlag, Berlin, 1985. \MR{MR774228
  (87k:58251)}

\bibitem[Ser65]{Ser65}
Jean-Pierre Serre, \emph{Alg\`ebre locale. {M}ultiplicit\'es}, Cours au
  Coll\`ege de France, 1957--1958, r\'edig\'e par Pierre Gabriel. Seconde
  \'edition, 1965. Lecture Notes in Mathematics, vol.~11, Springer-Verlag,
  Berlin, 1965. \MR{MR0201468 (34 \#1352)}

\bibitem[Smi01]{Smi01}
Gregory~G. Smith, \emph{Irreducible components of characteristic varieties}, J.
  Pure Appl. Algebra \textbf{165} (2001), no.~3, 291--306. \MR{MR1864474
  (2003a:16034)}

\bibitem[SST00]{SST00}
Mutsumi Saito, Bernd Sturmfels, and Nobuki Takayama, \emph{Gr\"obner
  deformations of hypergeometric differential equations}, Algorithms and
  Computation in Mathematics, vol.~6, Springer-Verlag, Berlin, 2000.
  \MR{2001i:13036}

\bibitem[ST04]{ST04}
Mutsumi Saito and William~N. Traves, \emph{Finite generation of rings of
  differential operators of semigroup algebras}, J. Algebra \textbf{278}
  (2004), no.~1, 76--103. \MR{MR2068067 (2005e:16041)}

\bibitem[Stu96]{Stu96}
Bernd Sturmfels, \emph{Gr\"obner bases and convex polytopes}, University
  Lecture Series, vol.~8, American Mathematical Society, Providence, RI, 1996.
  \MR{MR1363949 (97b:13034)}

\end{thebibliography}

\end{document}